\documentclass[12pt]{article}
\usepackage{amsmath,amsfonts,amssymb}
\usepackage{mathrsfs}
\usepackage{dsfont}
\usepackage[utf8x]{inputenc}
\usepackage[T1]{fontenc}
\usepackage{lmodern} \usepackage{ucs}
\usepackage{fullpage}

\usepackage{graphicx}
\usepackage{caption}
\usepackage{subcaption}

\usepackage{array}
\usepackage{multicol}
\usepackage{multirow}
\usepackage{color}
\usepackage{authblk}

\usepackage{tikz}

\newtheorem{lemma}{Lemma}[section]
\newtheorem{theo}[lemma]{Theorem}
\newtheorem{rmk}[lemma]{Remark}
\newtheorem{proposition}[lemma]{Proposition}
\newtheorem{defin}[lemma]{Definition}

    \newenvironment{ProofOf}[1]{\noindent
    \parindent=0pt\abovedisplayskip = 0.5\abovedisplayskip
    \belowdisplayskip=\abovedisplayskip{\bfseries Proof of  #1. }}{\QED\medskip}

\newcommand{\QED}{\mbox{}\hfill \raisebox{-0.2pt}{\rule{5.6pt}{6pt}\rule{0pt}{0pt}} \medskip\par}
\newcommand{\N}{\mathbb{N}}
\newcommand{\R}{\mathbb{R}}

\newcommand{\ds}{\displaystyle}
\newcommand{\ud}{\, {\mathrm{d}}}
\newcommand{\supp}{\mbox{supp}}

\DeclareMathOperator{\sgn}{sgn}

\title{Long time behaviour of interacting particles through a vibrating medium: comparison between the N-particle systems and the natural kinetic equation dynamics}

\author{Arthur Vavasseur\thanks{ {\tt avavasseur@bcamath.org}}}
\affil{Basque Center for Applied Mathematics (BCAM), Alameda Mazarredo 14, 48009 Bilbao, Spain}

\begin{document}
\maketitle

\abstract{We are interested in a kinetic equation intended to describe the interactions of particles with their environment. We focus on the long time behaviour. We prove that the time derivative of the spatial density goes to $0$ and exhibit the omega limit set for the distribution function. We then apply this result to the empirical density associated to a finite number of particles and prove that the speeds of all of them go to $0$. 
It also allows us to prove the convergence of all the positions with an increasing external potential and to get a precise description of their long time behaviour with a non decreasing  external potential. Those results allow us to prove that in very large time, the distribution function is not a good approximation of the $N$-particle system. From all those considerations, we get a very detailed description of the dynamic of the 
$N$-particle system.}

\vspace*{.5cm}
{\small
Inelastic Lorentz gas.} \\

\noindent{\bf Math.~Subject Classification.} 
82C70, 
70F45, 
37K05, 
74A25. 

\normalsize

\section{Introduction}

\subsection{Presentation of the model}

The paper is concerned with the long time behaviour of the solutions of the following Vlasov-like system 
\begin{eqnarray}
\label{kin}
\ds \partial_t f + v.\nabla_x f =\nabla_x ( V + \Phi).\nabla_v f,
\quad (x,v) \in \R^d\times\R^d, \quad t\geq 0,
\\ \ds
\label{Phi}
\Phi(t,x) = \Phi_0(t,x) - \int_0^t \int_{\R^d}
\rho(s,z) \Sigma(x-z) p(t-s) \ud z \ud s ,
\quad x \in \R^d, \quad t\geq 0,
\\ \ds
\label{rho}
\rho(t,x) = \int_{\R^d} f(t,x,v) \ud v,
\quad x \in \R^d, \quad t\geq 0,
\end{eqnarray} 
for a set of parameters $(p,\Sigma,\Phi_0)$ regular enough to ensure
$\Phi \in C(\R_+,W^{2,\infty}(\R^d))$ and a non negative external potential $V \in W^{2,\infty}_{loc}(\R^d)$. The system is completed with the initial data taken in the space of positive finite measures on $\R^d\times\R^d$:
\begin{equation}
\label{CI}
f\Big|_{t=0}=f_0.
\end{equation}
 Our interest for 
\eqref{kin}-\eqref{CI} come for the special case where $\Phi$ is defined by 
\begin{equation}
\label{VW}
\left\lbrace \begin{array}{ll} \ds
\Phi(t,x) 
= \int_{\R^d\times \R^n} \Psi(t,z,y) \sigma_1(x-z) \sigma_2(y) \ud z \ud y,
&  x \in \R^d, \quad t\geq 0,
\vspace{0.1cm}
\\ \ds
(\partial_t^2 \Psi -c^2 \Delta_y \Psi )(t,x,y)
= - (\sigma_1 \ast \rho)(x)\sigma_2(y) ,
& (x,y) \in  \R^d  \times \R^n, \quad t \geq 0,
\vspace{0.3cm}
\\ \ds
\Psi(0,x,y) =\Psi_0(x,y), 
\quad \partial_t \Psi(0,x,y) = \Psi_1(x,y),
& (x,y) \in \R^d \times \R^n .
\end{array} \right.
\end{equation}
By solving the linear wave equation in \eqref{VW}, it has been shown in \cite{GV} that 
\eqref{VW} can be recast as \eqref{Phi} for a certain set of parameters 
$(p,\Sigma,\Phi_0)$ determined by 
$c,\sigma_1,\sigma_2, \Psi_0,\Psi_1$ and $n$.

The system defined by \eqref{kin},\eqref{rho}-\eqref{VW} has been introduced in \cite{GV,AVPhD}. It describes the interactions of a large number of particles with their environement.
At any time $t\geq 0$, the particles are described by their distribution function $f_t$ acting on the phase space 
$\R^d_x\times\R^d_v$; the quantity 
$\int_{\mathcal{A}\times \mathcal{B}} \ud f_t(x,v)$ represents the mass of particles with position $x\in \mathcal{A}$ and speed 
$v\in \mathcal{B}$ at time $t$.
The environement is considered as a vibrating medium, its state is described by a scalar field $\Psi$. It can be thought of as an infinite set of membranes (one for each position $x\in \R^d$) which vibrate in a direction $y\in \R^n$ perpendicular to the particles' motion with a speed of wave $c>0$. The coupling
beteween the environement and the particles is described by the form factors $\sigma_1$ and $\sigma_2$ which are supposed to be smooth, non negative, compactly supported and radially symmetric. The particles are also submitted to the action of a certain external potential $V$.

The system defined by \eqref{kin}\eqref{rho}-\eqref{VW} is a kinetic version of an hamiltonian system introduced in \cite{BdB} which decribes the interactions between a single particle with its environement. By noting $q(t) \in \R^d_x$, the position of the particle at time $t\geq 0$, the dynamic is determined by the following equations
\begin{equation}
\label{bdb} 
\left\lbrace \begin{array}{l} 
\ddot{q}(t) = -\nabla V(q(t)) - \ds\iint_{\mathbb R^d\times \mathbb R^n} \sigma_1 (q(t)-z)\ \sigma_2(y) \ \nabla_x \Psi(t,z,y) \ud y \ud z,
 \\[.4cm]
(\partial_{tt}^2 \Psi- c^2 \Delta_z \Psi)(t,x,y) = -\sigma_1 (x-q(t))\sigma_2(y), \qquad x \in \mathbb R^d ,\  y \in \mathbb R^n. 
\end{array} \right. 
\end{equation}
In \cite{BdB}, the existence and uniqueness of the solution is proved but the main part of the paper is concerned with the large time behaviour. Under certain circonstances (roughly speaking when $n=3$, $c$ is large enough and $\hat{\sigma_2} \neq 0$ everywhere) it is proved for several kinds of external potentials
that the environement acts on the particle like a linear friction force (with a friction coefficient $\gamma$ explicit) in large time. The energy of the particle is evacuated in the menbranes while the energy of the whole system stay constant. In particular, the particle stops exponentially fast when $V=0$ and when $\lim_{|x| \rightarrow \infty} V(x)= + \infty$, in the second situation the limit of $q(t)$ is a singular point of $V$.
We refer the reader to  \cite{dBS,dBPS, LPdB, AdBLP,dBLP,dBP} for thorough  investigations on the model,
that contains both analytical treatments and numerical experiments.

The system \eqref{bdb} can easily be adapted in order to consider several particles $(q_k)_{1\leq k \leq N}$ interacting with the environement as a single one. By rescaling the interactions, we are led to the following system of equations, for all $k \in \{1,..,N\}$
\begin{equation}
\label{chmbdb} 
\left\lbrace \begin{array}{l} 
\ddot{q_k}(t) = -\nabla V(q_k(t)) - \ds\iint_{\mathbb R^d\times \mathbb R^n} \sigma_1 (q_k(t)-z)\ \sigma_2(y) \ \nabla_x \Psi(t,z,y) \ud y \ud z, 
 \\[.4cm]\ds
(\partial_{tt}^2 \Psi- c^2 \Delta_z \Psi)(t,x,y) = -\frac{1}{N} \sum_{k=1}^N \sigma_1 (x-q_k(t))\sigma_2(y) , \qquad x \in \mathbb R^d ,\  y \in \mathbb R^n. 
\end{array} \right. 
\end{equation}
We refer the reader to \cite{GV} for the precise motivations of this particular rescaling and its physical meaning. It is also proved in \cite{GV} that the empirical density 
$\hat{f}^N_t = \frac{1}{N} \sum_{k=1}^N \delta_{(q_k(t), \dot{q}_k(t))}$ converges to a solution of \eqref{kin},\eqref{rho}-\eqref{VW} when its initial data $\hat{f}^N_0$ converges to $f_0$. 
See e.g.\ \cite{FG,GMR} for a more general overview on the mean field regimes in statistical physic.

The existence of solution for \eqref{kin},\eqref{rho}-\eqref{VW}
has been established in \cite{BGV}. Certain asymptotic regimes has allowed us to connect this system with a family of Vlasov equations including the attractive gravitational Vlasov-Poisson system.  
The orbital stability of a family of equilibrium states has been investigated in \cite{BGV18}. More recently the Landau damping
has been proved in \cite{GVi} and some numerical simulation are in preparation in \cite{GVi2}. The relaxation to equilibrium with an additional dissipative Fokker-Planck operator has been established in \cite{AGV} thanks to hypocoercive methods.

However until now, the general long time behaviour of the solutions of \eqref{kin},\eqref{rho}-\eqref{VW} and even of the solutions
\eqref{chmbdb} have both remained unknown. It was not even clear if the dissipation of energy from the particles to the environement 
remained true in those new situations. The trend of $\dot{q}$ to $0$ established for one particle in \cite{BdB} had no equivalent either. The purpose of this paper is to fill that gap. We will see that the dissipation of energy holds. The trend of $\dot{q}$ to $0$ is replaced by the trend of $\partial_t \rho$ to $0$. The long time dynamic determined by \eqref{chmbdb} is similar to the dynamic of one particle but for the dynamic of the solutions of \eqref{kin},\eqref{rho}-\eqref{VW} with an initial data $f_0 \in L^1(\R^d\times\R^d)$, the consequences are quite different.
We have also cared about understanding why \eqref{kin}-\eqref{CI} led to such behaviour. The assumptions we make on the parameters are far less restrictive than the ones satisfied by the parameters 
$(p,\Sigma,\Phi_0)$ associated to \eqref{VW}. 

It is worth mentioning that the dynamic we prove for \eqref{chmbdb} has strong similarities with the one observed for a large family of $N$-particle systems with mean field interactions coresponding to Vlasov equations (see e.g.\ \cite{YBBDR04,BG13,AR95,JBM07}). Those system have been widely studied over the last decades by numerical and analytical works (see \cite{CDR,BGM10} for a review). However, to the author's knowledge, it is the first time that such kind of behaviour is proved rigorously.

\subsection{Statement of the results}

\subsubsection*{Long time behaviour with a measure initial data}

The main result of the present paper is a simple criterion leading to a general understanding of the long time dynamic of the solutions of \eqref{kin}-\eqref{CI}. Such criteria are not known for Vlasov equations. The situation is very different here due to the convolution with respect to the time standing in \eqref{Phi}. We now explain our criterion and its concequences without detailing the regularity hypothesis which are classical. 
We first extend $p$ on the whole line by assuming $p(-t) = -p(t)$. At least formally $\int_{\R} p(t) \ud t = 0$ and the values of the Fourier transform of any primitive function of $p$ are real. The main hypothesis is
\[ \nabla_x \Phi_0 \in L^1(\R_+, L^{\infty}(\R^d)),
\quad p = P',
\quad P \in L^1(\R), \quad \hat{P} \leq 0, \quad \hat{\Sigma} \geq 0, \quad V \geq 0,\]
where $\hat{P}$ and $\hat{\Sigma}$ are the Fourier transforms of $P$ and $\Sigma$.
It already allows us to prove that the total amount of energy transferred from the particles to the environement during a period $[0,T]$ is uniformly bounded with respect to $T$. By using this property we  get for all solution of \eqref{kin}-\eqref{CI} the following control
\begin{equation}
\label{int1}
\int_\R \left| \left|  \int_0^\infty P(t-s) \partial_t (\Sigma \ast \rho) (s) \ud s \right| \right|_{W^{1,\infty}(\R^d)}^2 \ud t \leq \mathcal{E}_2
\qquad \mbox{(Theorem \ref{existence}-3)} 
\end{equation}
where $\mathcal{E}_2$ is an explicit constant depending on the parameters.
Roughly speaking, if some oscillations of 
$t \mapsto \partial_t \Sigma \ast \rho(t)$ do not disappear in large time, then
the frequences of those oscillations belong to the zeros of 
$\hat{P}$. Under the additional restriction
$\hat{P}(\omega) \neq 0 \mbox{ for all } \omega \neq 0$, 
\eqref{int1} allows us to deduce
\begin{equation}
\label{int2}
\lim_{t\rightarrow \infty} \| \partial_t \Sigma \ast \rho(t) \|_{W^{1,\infty}(\R^d)} = 0
\qquad \mbox{(Lemma \ref{relax})} .
\end{equation}
Under our assumptions, it means that the force acting from the environement to the particles turns to be quasi-stationary in large time; setting $\kappa = -P(0)$, $f$ behaves as a solution of the following equation
\begin{equation}
\label{eq-lim}
\partial_t g + v.\nabla_x g 
= \nabla_x (V-\kappa \Sigma \ast \rho).\nabla_v g .
\end{equation}
By assuming $supp(\hat{\Sigma}) = \R^d$, 
\eqref{int2} also means that the spatial density 
$\rho$ is quasi-stationary in large time. 
$\rho$ being fixed, \eqref{eq-lim} is a linear transport equation. It has a unique global solution on $\R$ for all initial finite measure $g_0$ and it defines a group of operators $(\mathcal{S}^\rho_t)_{t\in\R}$. Finally our methods allows us to deduce that 
$f(t)$ turns to get closer an closer to the set of initial data 
$g$ satisfying
$\int_{\R^d_v} \ud (\mathcal{S}^{\rho_{g}}_t g )(v) \equiv \rho_{g}$ at any time $t$. We decided to call a spatially-stationary state such initial datas $g$ and to write 
$SEq$ the set of all the spatially-stationary states
(see definition \ref{semistat}). $SEq$ obviously contains the stationary states of \eqref{eq-lim} but the inclusion is strict in the general case as we will see in the appendix. We can now state our main result (see Theorem \ref{main} and Remark \ref{W-W1}):

\begin{theo}
\label{th-intro}
Under our assumptions, for any initial data 
$f_0 \in \mathcal{P}(\R^d\times\R^d)$ with 
finite energy,
if $\lim_{x\rightarrow \infty} V(x) = +\infty$ or more generally if $(f_t)_{t\geq 0}$ is tight, then for all $T>0$,
\[\begin{array}{llll}
(i) & \ds \lim_{t\rightarrow \infty} 
\|\partial_t \Sigma\ast \rho \|_{W^{1,\infty}(\R^d)} =0, &
(iii) & \ds
\lim_{T\rightarrow \infty} |\rho|_{C^{1/2}_{W_1}(T,+\infty)} = 0,
\vspace{0.2cm}
\\ (ii) & \ds
 \lim_{t\rightarrow \infty} W_1 (f(t),SEq) = 0, 
& (iv) & \ds 
\lim_{u \rightarrow \infty} \sup_{|t-u|\leq T} 
W_1(f(t), \mathcal{S}^{\rho(u)}_{t-u} f(u)) = 0.
\end{array}\]
\end{theo}

When $\Phi$ is also defined by \eqref{VW}, the only additionnal assumptions required to apply Theorem \ref{th-intro} are 
$n\geq 3$,
$\hat{\sigma_2} \neq 0$ on $\R^n$  and a stronger integrability condition on the initial datas $(\Psi_0,\Psi_1)$. Those hypothesis were already required in \cite{BdB} where the dimension $n=3$ is treated for one particle, it is natural to see them here as well. 
The counterexamples given in the appendix allows us to deduce that $SEq$ is the optimal attractive set of the dynamic. No damping can be expected for the whole distribution function in the general case (see Remark \ref{non-amortissement}). 

\subsubsection*{Long time behaviour for the $N$-particle system}

We then consider the finite particle system 
$(q_k,p_k)_{1\leq k \leq N}$ whose dynamic is determined by the following system of equations, for all $k \in \{1,..,N\}$
\begin{equation}
\label{chm}
\left\lbrace \begin{array}{l} \ds 
\dot{q_k}(t) = p_k(t), \\  \ds \dot{p_k}(t) = 
-\nabla_x (V + \Phi_0(t))(q_k(t))
+ \int_0^t p(t-s) 
\big(\frac{1}{N} \sum_{j=1}^N \nabla_x \Sigma(q_k(t)-q_j(s))\big) \ud s.
\end{array} \right. 
\end{equation}
Here again, it is an equivalent formulation of \eqref{chmbdb} for the same set of parameters $(p,\Sigma,\Phi_0)$. 
The system is completed with the initial datas
\begin{equation}
\label{chminit}
q_k(0)=q_{k,0} \quad \mbox{and} \quad p_k(0)) = p_{k,0}  \quad
\mbox{ for all } k \in \{1,..,N\}. 
\end{equation}
We define the empirical density
$\hat{f}^N(t) = \frac{1}{N} \sum_{k=1}^N \delta_{(q_k(t),p_k(t))}$.
It solves \eqref{kin}-\eqref{CI} for the initial data $\hat{f}^N(0)$ and it allows us to apply Theorem \ref{th-intro}. By considering such initial datas, we add some rigidity: the solution stay a sum of Dirac measures at any time. Any spatially-stationary state accessible by the dynamic $g$ can be written $g=\frac{1}{N} \sum_{k=1}^N \delta_{(x_k,v_k)}$
where $(x_k,v_k)_{1\leq k \leq N}$ belongs to $Eq^N$, the set of equilibrium of \eqref{chm}. If the trajectory of
$(q_k,p_k)_{1\leq k \leq N}$ is bounded, then 
$(q_k,p_k)_{1\leq k \leq N}$ is also quasi-stationary in large time:
\begin{equation}
\label{prop-intro}
\lim_{t\rightarrow \infty} 
D\big((q_k(t),p_k(t))_{1\leq k \leq N}, Eq^N \big) = 0. 
\qquad \mbox{(Proposition \ref{kin-part})}.
\end{equation}
In particular, the speeds of all the particles $(p_k(t))_{1\leq k \leq N}$ goes to $0$ when $t$ goes to infinity. In order to be more precise, we add the following restriction
\[ {\bf (M)} \quad \lim_{x\rightarrow \infty} V(x) = + \infty, 
\quad x.\nabla V(x) \geq 0, \quad x.\nabla \Sigma(x) \leq 0. \]

When $x.\nabla V(x) >0$ for all $x\neq 0$, it allows us to prove the large time convergence of the positions and velocity to $0$ (Theorem \ref{monotonie}-$1)$. 
It generalizes the convergence obtained for one particle when $n=3$ in \cite{BdB} to any finite number of particles and any dimension $n\geq 3$.

In order to understand better the effects of the interactions between the particles, it would be interesting to study the case $V=0$. 
A first step to describe this particular situation would be to prove that the system stay confined uniformly in time. It is wrong in the general case for the solutions of \eqref{kin}-\eqref{CI} with an initial data $f_0 \in \mathcal{M}_+^1(\R^d\times\R^d)$ 
(see Remark \ref{non-conf-simp}); it requires some specific technics for the $N-$particles systems that we do not wish to develop in the present paper. We only partially answer this question by considering an external potential $V$ satisfying {\bf (M)}. To fix the ideas one can think $V(x) = [|x|-R]_+^2$. In that case the distance between $(q_k(t))_{1\leq k \leq N}$ and  
$\overline{B(0,R)}^N$ goes to $0$, we prove that the particles turn to form a bounded number of clusters. Picking two particles, there is only two scenarios: or their distance goes to $0$ and they belong to the same cluster, or they do not interact in large time
(Theorem \ref{monotonie}-$2$). 

\subsubsection*{Comparison between the long time dynamics
and concequences}

In the end, we revisit the mean-field limit established in \cite{GV} by taking a random sequence of initial datas 
$(q_{k,0},p_{k,0})_{k\geq 1}$ identically distributed according to 
$f_0 \in L^1(\R^d\times\R^d)$. Almost surely, up to a subsequence $N_k$, the empirical density $\hat{f}^N$ converges weakly to a solution of \eqref{kin}-\eqref{CI} when $N$ goes to infinity. We prove it by two different methods. The first one is inspired by the mean field limit established for the Vlasov equation  by Dobrushin in \cite{Dob}. This method does not require any additionnal restriction but it does not provide any explicit estimates in that particular context. In the second proof inspired by \cite{Sz,SM} we add the restriction $\nabla^2 V \in L^\infty(\R^d)$ and get some explicit estimates thanks to the coupling methods. 

By taking a $L^1$ initial data, we also add some rigidity: the solution $f(t)$ is 
equimesurable to $f_0$  at any time. When $V$ is a confining potential, the measure $f(t)\ud x\ud v$ cannot concentrate on the support $\R^d_x \times \{0_{\R^d_v} \}$
while $\hat{f}^N$ allways does in large time by \eqref{prop-intro}. More precisely, we can find an explicit constant 
$\Theta>0$ such that almost surely,
\begin{equation}
\label{int3}
\sup_{N \geq 1} \liminf_{t \rightarrow \infty}
W_1(\hat{f}^N(t), f(t))
\geq \Theta >
\sup_{t\geq 0} \limsup_{N \rightarrow \infty} W_1(\hat{f}^N(t), f(t)) =0
\end{equation}
(see Theorem \ref{lim-chm}). It allows us to prove that the convergences exhibited in Theorem \ref{monotonie} turn to get slower and slower when $N$ goes to infinity almost surely. 
When $\nabla^2 V \in L^\infty$, we give some explicit lower bounds for those convergences (Proposition \ref{bornes-infs}). 
The instability given by \eqref{int3} does not invalidate the relevance of \eqref{kin}-\eqref{CI} with an initial data in $L^1$ to describe finite particle systems. It simply ensures that the long time behaviour exhibited for the distribution function $f(t)$ should be seen as a mean time behaviour for the empirical density 
$\hat{f}^N$ when $N$ is large enough. When $f(t)$ is uniquely defined and $N$ is large enough to allow $f(t)$ to get close enough to $SEq$ during the time of validity on the mean-field approximation, our results allows us to split the time in four different periods (to be compared with e.g.\ \cite{YBBDR04}):
\begin{itemize}
\item[{\bf(1)}] 
At the beginning $\hat{f}^N$ is close to $f$, it moves from its initial position to get in a neighbourhood of $SEq^{f_0}$, the set of spatially-stationary states accessible by $f$.
\item[{\bf(2)}] 
Then  $\hat{f}^N$ is still close to $f$ but also to 
$SEq^{f_0}$ . The empirical spatial density is quasi-stationary.
\item[{\bf(3)}]  After a certain time $f(t)$ is no longer a good approximation of $\hat{f}^N(t)$. $\hat{f}^N(t)$ leaves the neighbourhood of $SEq^{f_0}$ to get in a neighbourhood of the set of empirical density associated with the equilibrium of \eqref{chm}.
\item[{\bf(4)}]
In the end, the positions and speeds of all the particles are close to the equilibrium $Eq^N$, they are at least quasi-stationary and possibly converge.
\end{itemize}

\subsection{Organization of the paper}

The paper is organized as follows. Section \ref{sec1}  contains the preleminary discussion. We recall how \eqref{VW} can be recast as \eqref{Phi} and prove that the parameters 
($p,\Sigma,\Phi_0$) satisfy all the property required to study the large time behaviour of the solutions of \eqref{kin}-\eqref{CI} by our methods. It also explains how the solution of \eqref{kin} can be expressed by mean of characteristic curves and establishes some important estimates on those curves. The last part of section \ref{sec1} is devoted do define the spatially-stationary states and to exhibit some basic property of those distributions. 
In section \ref{sec2}, we improve the results of existence and uniqueness for the solution of \eqref{kin}-\eqref{CI} already established in \cite{BGV} and we prove the key estimate \eqref{int1}.
Section \ref{sec3} is devoted to prove Theorem \ref{th-intro}  and section \ref{sec4} is devoted to exhibit its concequences on the finite particle systems. 

\subsection{Notations and conventions}

Taking $S=\R^d$ or $S=\R^d\times\R^d$, we note 
$C_b(S),C_c(S)$ and $C_0(S)$, the spaces of continuous function on $S$ which are respectively bounded, compactly supported and vanishing at infinity. 
$\mathcal{M}^1_+(S)$ is the set of non negative finite measures on $S$. For all $\chi \in W^{1,\infty}(S)$, we set 
$\|\chi\|_{W^{1,\infty}(S)} = 
\|\chi\|_{L^\infty(S)} + \|\nabla \chi\|_{L^\infty(S)}$. The norm on the dual of $W^{1,\infty}(S)$ define a distance $W$ on 
$\mathcal{M}^1_+(S)$:
\[ W(\mu,\nu) := \sup_{\|\chi\|_{W^{1,\infty}(S)} \leq 1} 
\left( \int_S \chi(z) \ud \mu(z) - \int_S \chi(z) \ud \nu(z) \right). \]
The use of the strong norm in $W^{1,\infty}(\R^d\times\R^d)'$ appeared naturally in \cite{BGV},
it will also be convenient here. 
$W$ is related to the Kantorowich-Rubinstein distance $W_1$ defined on the probability measures on $S$ by 
\[ W_1(\mu,\nu) := \inf_{\pi} \int_{S^2} (|z_1-z_2|\wedge 1) \ud \pi(z_1,z_2)\]
where the infimum is taken on the set of measures $\pi$ acting on $S^2$ such that $\int_S \ud \pi(\cdot,z_2) \equiv \mu$ and 
$\int_S \ud \pi(z_1,\cdot) \equiv \nu$.
We refer the reader to  \cite{Dob,CV} for a detailed introduction to this notion.
We recall (see \cite{Dob} or \cite[Chapter~6]{CV}) that $W_1$ metrizes the tight convergence on the probability measures (the weak convergence on $C_b(S)$). It is also related to the strong convergence on the dual of $W^{1,\infty}(S)$ thanks to the
Kantorowich-Rubinstein identity for the distance 
$d(x,y)=|x-y|\wedge 1$:
\begin{equation}
\label{KR}
W_1(\mu,\nu) = 
\sup_{ \ds \begin{subarray}{c} 0\leq \chi \leq 1 \\
 \|\nabla \chi\|_{L^\infty} \leq 1 \end{subarray} }
\left( \int_S \chi(z) \ud \mu(z) - \int_S \chi(z) \ud \nu(z) \right)
\end{equation}
(see \cite[Theorem~5.10 or even more precisely Remark~6.5]{CV}).
It allows us to deduce that $W$ and $W_1$ are equivalent on the probability measures:
\[ \frac{1}{2} W(\mu,\nu) \leq W_1(\mu,\nu) \leq 2 W(\mu,\nu). \]
Since 
$W(\mu,\nu) \geq \left| \int_S \ud \mu - \int_S \ud \nu \right|$
and $W(\mu,\nu) = \mathfrak{m} 
W \big(\frac{\mu}{\mathfrak{m}},\frac{\nu}{\mathfrak{m}}\big)$, 
$W$ metrizes the tight convergence on the measures as well as $W_1$ does. Then, we recall that a sequence $(\mu^k)_{k\geq 0}$
converges for the distance $W$ if and only if it converges for the weak star topology on $C_0(S)$ and is tight:
$\lim_{R\rightarrow \infty} \sup_{k\geq 0} 
\int_{|z|\geq R} \ud \mu^k(z) =0$.
We note $C_{w\ast}(\R_+,\mathcal{M}^1_+(S))$ (resp. $C_{W}(\R_+,\mathcal{M}^1_+(S))$) the set of measured valued functions $t\mapsto \mu_t$ continuous for the weak star topology on $C_0(S)$ (resp. for the distance $W$).  In the end, for any connex 
set $I\subset \R$, we define an Hölder continuity estimate on
$C_{W}(\R_+,\mathcal{M}^1_+(S))$:
\[ |\mu|_{C^{1/2}_{W}(I)}:= \sup_{s,t\in I} 
\frac{W(\mu_s,\mu_t)}{|s-t|^{1/2}}. \] 

\section{Preliminary discussion}\label{sec1}

\subsection{Simplification of the Vlasov-Wave system}

We recall the derivation of \eqref{kin}-\eqref{CI} from 
\eqref{kin},\eqref{rho}-\eqref{VW} established 
in \cite{BGV,GV} 
under the hypothesis of \cite{BdB}:
\[{\bf (A1)}\qquad
\left\{
\begin{array}{l}
\sigma_1\in C^\infty_c(\mathbb R^d,\mathbb R), \ \sigma_2\in C^\infty_c(\mathbb R^d,\mathbb R),
\\
\sigma_1(x)\geq 0,\ \sigma_2(y)\geq 0\  \text{for any $x\in\mathbb R^d$, $y\in\mathbb R^n$},
\\
\text{$\sigma_1,\sigma_2$ are radially symmetric},
\end{array}
\right.\]
\[{\bf (A2)}\qquad
\Psi_0, \Psi_1\in 
 L^2(\mathbb R^d\times\mathbb R^n).
\]
Assuming {\bf (A1)-(A2)}, we define
\begin{equation}
\label{reduct-param}
\begin{array}{l} \ds
\Sigma = \sigma_1 \ast \sigma_1, 
\quad p(t)=\frac{1}{(2\pi)^n} \int_{\R^n}
\frac{\sin(c|\xi|t)}{c|\xi|} |\hat{\sigma_2}(\xi)|^2 \ud \xi,
\\ \ds \Phi_0(t,x)
= \int_{\R^d\times \R^n} \tilde{\Psi}(t,z,y) \sigma_1(x-z) \sigma_2(y) \ud z \ud y 
\end{array}
\end{equation}
where $\tilde{\Psi}$ is the unique solution of the homogeneous wave equation with initial datas $(\Psi_0,\Psi_1)$ in 
$C(\R_+;L^2(\R^d\times\R^n))$. By solving the wave equation in 
\eqref{VW}, we have already established previously

\begin{lemma}{\cite{BGV,GV}}
\label{resolution}
Assume {\bf(A1)-(A2)}, take $f $ in $C_{W}(\R_+;\mathcal{M}_+^1(\R^d\times\R^d))$ and $(p,\Sigma,\Phi_0)$ defined by \eqref{reduct-param}. Then, the self consistent potential $\Phi$ defined by \eqref{rho} and \eqref{VW}
can be recast as \eqref{Phi}.
\end{lemma}
We now exhibit some conditions to ensure that
$(p,\Sigma,\Phi_0)$ satisfy the required property for the rest of the paper.

\begin{proposition}
\label{prop-P}
Take $n\geq 3$ and assume {\bf (A1)}, then
\begin{enumerate}
\item $p=P'$ where $P \in C^\infty(\R)$, 
$P^{(k)} \in (L^1 \cap L^\infty)(\R)$ for all $k\geq 0$ and
$\hat{P} \geq 0$.
\item If $\hat{\sigma}_2(\xi) \neq 0$ for all $\xi \neq 0$, then 
$\hat{P}(\omega) \neq 0$ for all $\omega \neq 0$.
\end{enumerate}
\end{proposition}

\begin{proposition}
\label{intphi0}
Assume {\bf(A1)-(A2)}, then $\Phi_0 \in C(\R_+,W^{2,\infty}(\R^d))$ and under one of the following restrictions
\begin{itemize}
\item
$n=3$ and 
$(x,y) \mapsto (1+|y|)^{\epsilon}(|\Psi_0(x,y)| + |\Psi_1(x,y)|)$ lies in $L^{q}(\R^d_x;L^1(\R^n_y))$ for $q,\epsilon$ such that 
$1 \leq q \leq \infty$ and $\epsilon>0$, 
\item
$n\geq 4$ and $\Psi_0, \Psi_1 \in L^{q}(\R^d_x;L^r(\R^n_y))$ for 
$q,r$ such that 
$1 \leq q \leq \infty$ and $1\leq r <\frac{2n}{n+3}$,
\end{itemize}
$\Phi_0 \in L^1(\R_+,W^{1,\infty}(\R^d))$ and
$\nabla_x \Phi_0 \in L^\infty(\R_+\times\R^d)$.  
\end{proposition}
About $\Sigma$, we just point out that the condition $supp(\hat{\Sigma}) = \R^d$ mentioned in the introduction is satisfied since $\hat{\sigma_1}$ is an analytic function under {\bf (A1)} (see \cite{Mit}).
Before proving those claims, we summarize here the assumptions we need to add to {\bf(A1)-(A2)} in order to ensure that $p,\Sigma$ and $\Phi_0$ defined by 
\eqref{reduct-param} satisfy all the conditions required in
all the rest of the paper.
\begin{itemize}
\item $n\geq 3$ and $\hat{\sigma}_2(\xi) \neq 0$
for all $\xi \neq 0$,
\item The integrability assumption on $\Psi_0$,$\Psi_1$ given 
in Proposition \ref{intphi0}, 
\item $x. \nabla_x \sigma_1(x) \leq 0$ for all $x\in \R^d$.
\end{itemize}

The two first points generalize the assumptions made in \cite{BdB}
where the dissipative behaviour of the solutions of \eqref{bdb} is only exhibited when $n=3$ under the assumption 
$\hat{\sigma}_2 \neq 0$ and a more restrictive condition on 
$\Psi_0$ and $\Psi_1$. The dimension $n=3$ is particular anyway,
from \eqref{Pchap} below it is the only one where 
$\hat{P}(0)$ can be different from $0$. By using the identity 
$-2 \int_0^\infty  p(s) s \ud s = \hat{P}(0)$,
the friction coefficient $\gamma$ explicited by the computations 
 of \cite{BdB} is equal to $0$ for all dimension $n\geq 4$. 
The results of the current paper hold anyway when $\hat{P}(0) = 0$ (at the cost of a less direct proof, see Remark \ref{Pitt} below) but they are far less precise than the results established for a single particle in \cite{BdB}. We believe that deeper investigations would reveal some important
differences when $n\geq 4$.

The assumption on $x.\nabla_x \sigma_1(x)$ is new. It ensures the condition on $\Sigma$ in {\bf(M)} to be satisfied
(see Appendix B below). It is just a tool to simplify the expression of the equilibrium of \eqref{chm}. It is not a real constrained since $\sigma_1$ and $\sigma_2$ are introduced as regular approximations of the Dirac measures $\delta_0$ of $\R^d$ and $\R^n$ in \cite{BdB}.

\begin{ProofOf}{Proposition \ref{prop-P}}


From \cite{BGV}, we already know that $p\in L^1(\R)$ when 
$n\geq 3$, we start by computing its Fourier transform. 
Taking $\chi$ in 
$\mathscr S(\R)$, we have by definition
\[ \int_{\R}  \hat{p}(\omega) \chi(\omega) \ud \omega 
\begin{array}[t]{l} \ds
=\int_{\R} p(t) \hat{\chi}(t) \ud t \\ \ds
=\frac{1}{(2\pi)^n} \int_{\R^n} \frac{|\hat{\sigma_2}(\xi)|^2}{c|\xi|} 
\left( \int_{\R} \sin(c|\xi|t) \hat{\chi}(t) \ud t \right) \ud \xi \\ \ds
=\frac{1}{(2\pi)^n} \int_{\R^n} \frac{|\hat{\sigma_2}(\xi)|^2}{c|\xi|} 
\left\langle \widehat{\sin(c|\xi|\cdot)} ; \chi 
\right\rangle_{\mathscr{S}'\times \mathscr{S}} \ud \xi \\ \ds
=\frac{1}{(2\pi)^n} \int_{\R^n}  \frac{|\hat{\sigma_2}(\xi)|^2}{c|\xi|} 
i\pi (-\chi(c|\xi|)+\chi(-c|\xi|)) \ud \xi \\ \ds
=\frac{i\pi|\mathbb{S}^{n-1}|}{(2\pi)^n} \int_0^\infty 
\frac{|\hat{\sigma_2}(r e_1)|^2}{c r} 
 (-\chi(c r)+\chi(-c r)) r^{n-1}\ud r \\ \ds
=-\frac{i\pi|\mathbb{S}^{n-1}|}{(2\pi)^n} \int_{\R} \chi(s)
\mbox{sgn}(s) |\hat{\sigma_2}(s e_1/c)|^2 |s|^{n-2} \frac{\ud s}{c^n}.
\end{array}\]
Finally, 
$ \hat{p}(\omega) = -i \frac{\pi|\mathbb{S}^{n-1}|}{(2\pi c)^n}
\mbox{sgn}(\omega) |\omega|^{n-2} |\hat{\sigma_2}(\omega e_1/c)|^2  $.
Since $n\geq 3$, $P$ is well defined by the expression of its Fourier transform
\begin{equation}
\label{Pchap} 
\hat{P}(\omega) = - \frac{\pi|\mathbb{S}^{n-1}|}{(2\pi c)^n}
|\omega|^{n-3} |\hat{\sigma_2}(\omega e_1/c)|^2  
\end{equation}
and we can check that 
$\hat{p}(\omega) = i\omega \hat{P}(\omega) = \widehat{P'}(\omega)$. We set 
$a_n = \frac{(\pi |\mathbb{S}^{n-1}|)^{1/2}}{(2\pi )^{n/2}} $
and for all $t \in \R$, we define
\[\tilde{\sigma_2}(t):= \int_{\R^{n-1}} 
\sigma_2(t,y_2,\cdots,y_n) \ud y_2 \cdots \ud y_n \]
so that 
$ \widehat{\tilde{\sigma_2}}(\omega) = \hat{\sigma_2}(\omega e_1)$
and $\hat{P}(\omega) = 
-\left| \frac{a_n}{\sqrt{c}}  \frac{1}{c}\left|\frac{\omega}{c}\right|^{\frac{n-3}{2}}  
\widehat{\tilde{\sigma_2}}\left(\frac{\omega}{c}\right)\right|^2$. 
It allows us to set
\begin{equation}
\label{q1}
 P = -\check{q_1} \ast q_1 \quad \mbox{for} 
\quad q_1(s)
= \left\lbrace \begin{array}{ll}  \ds
\frac{a_n}{\sqrt{c}}  \tilde{\sigma_2}^{(m)}(c s) 
& \mbox{if}\quad n=2m+3,\vspace{0.2cm} \\  \ds
\frac{a_n}{\sqrt{c}}  (-\Delta)^{1/4} \left( \tilde{\sigma_2}^{(m)}\right)
(c s) 
& \mbox{if} \quad n=2m+4,
\end{array} \right.
\end{equation}
where the classical fractional laplacian is given by 
\[ (-\Delta)^{1/4}u (x) = c_{1,\frac{1}{4}} 
\int_\R \frac{u(x)-u(y)}{|x-y|^{3/2}} \ud y .\]
Since $\sigma_2 \in C^\infty_c(\mathbb R^n,\mathbb R)$, 
$\tilde{\sigma_2}^{(m)}$ belongs to $C^\infty_c(\mathbb R,\mathbb R)$ for all $m\geq 0$ and it is easy to check that its fractional laplacian also belongs to $(L^1\cap L^\infty)(\R)$, therefore so does $q_1$. The expression of 
$ q_1^{(k)}$ is directly given  by changing $m$ in $m+k$ in \eqref{q1}, therefore 
$q_1^{(k)} \in L^1(\R)$ for all $k\geq 0$. From the Young inequality we deduce that
$P^{(k)}= -\check{q_1} \ast q_1^{(k)} $ also belongs to 
$(L^1\cap L^\infty)(\R)$. 
We have proved Proposition \ref{prop-P}-$1)$. 
Proposition \ref{prop-P}-$2)$ follows directly from \eqref{Pchap}.

%
\end{ProofOf}

\begin{rmk}{The proof also allows us to state the following points:}
\begin{itemize}
\item
If $n \in 2\N+3$, then $P \in C_c^\infty(\R)$ by \eqref{q1}.
\item If one remove the hypothesis $\sigma_2 \geq 0$ from {\bf (A1)}, the computations leading to \eqref{q1} also make sens for $n=1,2$ when 
$\int_{\R^d} \sigma_2 \ud x = 0$. In that case, setting 
$S(t) = \int_{-\infty}^t \tilde{\sigma_2}(s) \ud s$ we get
$q_1(t) = \frac{1}{\sqrt{c}} S(ct) $ when $n=1$ and
$q_1(t) = \frac{1}{\sqrt{2c}} (-\Delta)^{1/4} (S)(ct)$ when $n=2$.
By adding some restrictions on $\Psi_0$ and $\Psi_1$ adapted to those dimensions in order to ensure $\Phi_0 \in L^1(\R_+,W^{1,\infty}(\R^d))$, all the results of the paper would also hold for 
$n =1,2$ in that particular situation. 
\end{itemize}
\end{rmk}

\begin{ProofOf}{Proposition \ref{intphi0}}

The solution of the homogeneous wave equation
$t\mapsto \tilde{\Psi}(t) \in L^2(\R^d\times\R^n)$ is continuous on $\R_+$ therefore $\Phi_0 \in C(\R_+,W^{2,\infty}(\R^d))$ directly from the Hölder and the Young inequalities. 

For any $n\geq 3$, we
use the representation formula 
$\tilde{\Psi}(t) = K^0_t \underset{y}{\ast} \Psi_0
+K^1_t \underset{y}{\ast} \Psi_1 $ where $K^0,K^1$ are the green functions associated with the wave equation
(see e.g.\ \cite{Wa}), it allows us to set
\begin{equation}
\label{intphi1}
\begin{array}{l} \ds
 \left| \int_{\R^n} \tilde{\Psi}(t,x,y) \sigma_2(y) \ud y \right|
 \\ \qquad \ds
= \left| \int_{\R^n} 
\big(\Psi_0(x,y) (K^0_t \ast \sigma_2 )(y)
+\Psi_1(x,y) (K^1_t \ast \sigma_2)(y)\big) \ud y \right|
\vspace{0.2cm}
 \\ \qquad  \ds
\leq \|\Psi_0(x,\cdot)\|_{L^r(\R^n)} 
\| K^0_t \ast \sigma_2 \|_{L^{r^\ast}(\R^n)}
+  \|\Psi_1(x,\cdot)\|_{L^r(\R^n)} 
\| K^1_t \ast \sigma_2 \|_{L^{r^\ast}(\R^n)}.
\end{array}
\end{equation}
From \cite{Wa}, since $\sigma_2 \in C_c^\infty(\R^n)$, there exists two constants $c_0(\sigma_2,\mu)$ and $c_1(\sigma_2,\mu)$
such that for all $t\geq 0$,
\begin{equation}
\label{Wahl}
\begin{array}{l} \ds 
\| K^0_t \ast \sigma_2 \|_{L^{1/(1-\lambda -\mu)}(\R^n)}
\leq  (1+t)^{n(1-\lambda)-\frac{n-1}{2}} c_0(\sigma_2,\mu),
\\ \ds
\| K^1_t \ast \sigma_2 \|_{L^{1/(1-\lambda -\mu)}(\R^n)}
\leq  (1+t)t^{n(1-\lambda)-\frac{n-1}{2}} c_1(\sigma_2,\mu),
\end{array}
\end{equation}
for all set of parameters $\lambda,\mu$ satisfying
$\mu+\lambda\leq 1$, $0\leq \mu<1$ and $\lambda>0$.

When $n\geq 4$, we apply \eqref{Wahl} with $\mu=0$ and $\lambda = \frac{1}{r}$,
since $r^\ast = \frac{1}{1-\lambda-\mu}$, we get
from \eqref{intphi1}:
\begin{equation}
\label{ngeq4}
 \left| \left| 
\int_{\R^n} \tilde{\Psi}(t,\cdot,y) \sigma_2(y) \ud y 
\right| \right|_{L^q(\R^d_x)} 
\leq 
C (1+t)^{n(1-\frac{1}{r})-\frac{n-1}{2}} 
\end{equation}
where we have set 
$C= c_0(\sigma_2,0) \| \Psi_0\|_{L^q(\R^d_x; L^r(\R^n_y))} 
+c_1(\sigma_2,0) \| \Psi_1 \|_{L^q(\R^d_x; L^r(\R^n_y))}$. The exponent can be rewriten 
$n(1-\frac{1}{r})-\frac{n-1}{2}
= -1-n(\frac{1}{r}-\frac{n+3}{2n})<-1$.

When $n=3$, we apply \eqref{Wahl} with $\mu=0$ and $\lambda = 1$
to get $\| K^0_t \ast \sigma_2 \|_{L^\infty(\R^3)}
\leq  (1+t)^{-1} c_0(\sigma_2,0)$ and 
$\| K^1_t \ast \sigma_2 \|_{L^\infty(\R^3)}
\leq  (1+t)^{-1} c_1(\sigma_2,0)$ but it is not enough to get the 
integrability with respect to $t$. 
From the Kirchhoff's formula 
(see e.~g.~\cite[Eq. (22), Chapter~2.4, p.~73]{Evans}), 
we also have
\[\begin{array}{l} \ds 
K^0_t \ast \sigma_2(y) = \frac{1}{4\pi c^2 t^2} \int_{|y-y'|=ct}
\big( \sigma_2(y') + \nabla_y \sigma_2(y').(y'-y) \big) \ud S(y'),
\\ \ds
K^1_t \ast \sigma_2(y) = \frac{1}{4\pi c^2 t} \int_{|y-y'|=ct}
\sigma_2(y') \ud S(y').
\end{array} \]
Taking $R>0$ such that $\sup(\sigma_2) \subset B(0,R)$, it allows us to deduce that $K^0_t \ast \sigma_2$ and $K^1_t \ast \sigma_2$
are both supported in $B(0,ct+R) \setminus B(0,ct-R)$ for any 
$t\geq R/c$. Coming back to \eqref{intphi1}, we get for all 
$t \geq \frac{R}{c}$
\begin{equation}
\label{negale3}
\begin{array}{l} \ds
 \left| \int_{\R^3} \tilde{\Psi}(t,x,y) \sigma_2(y) \ud y \right|
 \\ \qquad \ds
= \left| \int_{|y|\geq R-ct} 
\big(\Psi_0(x,y) (K^0_t \ast \sigma_2 )(y)
+\Psi_1(x,y) (K^1_t \ast \sigma_2)(y)\big) \ud y \right|
 \\ \qquad  \ds
\leq 
\|\Psi_0(x,\cdot)\|_{L^1(\complement B_y(0,ct-R))} 
\| K^0_t \ast \sigma_2 \|_{L^\infty(\R^3)}
\\ \qquad \qquad \ds
+  \|\Psi_1(x,\cdot)\|_{L^1(\complement B_y(0,ct-R))}
\| K^1_t \ast \sigma_2 \|_{L^\infty(\R^3)}
 \\ \qquad \ds
\leq  
(1+ct-R)^{-\epsilon}
\|(1+|\cdot|)^\epsilon \Psi_0(x,\cdot)\|_{L^1(\R^3)} 
(1+t)^{-1} c_0(\sigma_2,0)
\\ \qquad \qquad \ds
+  (1+ct-R)^{-\epsilon}
\|(1+|\cdot|)^\epsilon \Psi_1(x,\cdot)\|_{L^1(\R^3)} 
(1+t)^{-1} c_1(\sigma_2,0).
\end{array}
\end{equation}

For all $n\geq 3$, \eqref{ngeq4} and \eqref{negale3} allows us to 
find $C,\eta>0$ such that for all $t\geq 0$,
\[
 \left| \left| 
\int_{\R^n} \tilde{\Psi}(t,\cdot,y) \sigma_2(y) \ud y 
\right| \right|_{L^q(\R^d_x)} 
\leq 
\frac{C}{(1+t)^{1+\eta}}.
\]
Thanks to the Young inequality, for all $t\geq 0$ and all
$k \in \N$, we finally get
\[
\|\Phi(t)\|_{W^{k,\infty}(\R^d)}  
\leq \|\sigma_1\|_{W^{k,q^\ast}(\R^d)} \frac{C}{(1+t)^{1+\eta}}.
\]
Therefore $\Phi_0$ belongs to $L^1(\R_+,W^{1,\infty}(\R^d))$
and $\nabla \Phi_0 \in L^\infty(\R_+\times\R^d)$.

\end{ProofOf}

\subsection{Uniform estimates on the characteristic curves}

$V \in W^{2,\infty}_{loc}(\R^d)$ and $\Phi\in C(\R_+;W^{2,\infty}(\R^d))$ being fixed, we recall that the weak solutions of \eqref{kin} can be explicited thanks to the characteristic curves, solutions of the following ODE system

\begin{equation}
\label{char}
\left\lbrace \begin{array}{l}
\dot{X} = \xi \\ \dot{\xi} = -\nabla(V+\Phi(t))(X)
\end{array} \right. 
\qquad , \qquad
\left\lbrace \begin{array}{l}
X(\alpha) = X_{in} \\ \xi(\alpha) = \xi_{in}
\end{array} \right. 
\end{equation}
$\alpha \geq 0$ being fixed, we note $\varphi^t_\alpha (x,v) 
= \big(\varphi^t_{\alpha,x} (x,v),\varphi^t_{\alpha,v} (x,v)\big) $, the solution of \eqref{char} with initial data $(X_{in},\xi_{in}) =(x,v)$ at time $t$. $\varphi$ is well defined thanks to the Cauchy-Lipschitz theorem and the following bounds:

\begin{lemma}{Uniform estimates on the characteristic curves}
\label{char-est}

If $V \geq 0$ 
then for any $(X,\xi)$ solution of \eqref{char}, 
setting  $E_\alpha = V(X_{in}) + \frac{|\xi_{in}|^2}{2}$, we have for all $t \geq 0$
\[ \left\lbrace \begin{array}{l} \ds 
|\xi(t)| \leq (2 E_\alpha)^{1/2} 
+ \int_{\alpha \wedge t}^{\alpha \vee t} \|\nabla_x \Phi(s)\|_{L^\infty(\R^d)} \ud s,
\\ \ds
|X(t)-X_\alpha| \leq (2 E_0)^{1/2}|t -\alpha|
+ \int_{\alpha \wedge t}^{\alpha \vee t} + \int_{\alpha \wedge t}^s
\|\nabla_x \Phi(\tau)\|_{L^\infty(\R^d)} \ud \tau \ud s.
\end{array} \right.\]
\end{lemma}
Taking $f$ a solution of \eqref{kin}\eqref{CI}, 
$\chi_0 \in C_c^\infty(\R^d\times\R^d)$
and $T>0$, we set 
$\chi(t,x,v) = (\chi_0 \circ \varphi_t^T)(x,v)$ and we observe
\[ \frac{\ud}{\ud t} < f_t ; \chi(t) > 
\begin{array}[t]{l} \ds 
= < f_t ; v.\nabla_x \chi(t) - \nabla_x (V+\Phi(t)).\nabla_v \chi(t)> 
\\ \ds \qquad
-< f_t ; v.\nabla_x \chi(t) - \nabla_x (V+\Phi(t)).\nabla_v \chi(t)> 
\\ \ds 
=0.
\end{array}\]
In particular, for $t=0$ and $t=T$, we get for all solution of \eqref{kin} the identity
\begin{equation}
\label{rep-car}
\int_{\R^d\times\R^d} \chi(x,v) \ud f_T (x,v) 
= \int_{\R^d\times\R^d} (\chi \circ  \varphi_0^T)(x,v) 
\ud f_0 (x,v)
=: <{\varphi_0^T}_{\#} f_0 ; \chi > .
\end{equation}
It is easy to check that 
$t \mapsto {\varphi_0^t}_{\#} f_0 $ solves
\eqref{kin}\eqref{CI} as well. 
We recall that since the flow $\varphi$ is defined by a divergence free vector field ($\mbox{div}_{\R^{2d}} [(X,\xi)\mapsto 
(\xi,-\nabla_x(V+\Phi(t))(X))] = 0$ at any time), it preserves the measure: for any Borel subset $\mathcal{A}$ of 
$\R^d\times\R^d$ and any $\alpha,\beta \geq 0$,
$\left| \varphi_\alpha^\beta (\mathcal{A}) \right| =|\mathcal{A}|$.
When $f_0 \in L^1(\R^d\times\R^d)$, it allows us to transform
\eqref{rep-car} in $f_t = f_0 \circ \varphi_t^0$. We end by proving the estimates on $\varphi_\alpha^t$.

\begin{ProofOf}{Lemma \ref{char-est}}

By changing $\Phi(s)$ in $\Phi(\alpha-s)$ or 
$\Phi(\alpha+s)$, we can always suppose 
$0 = \alpha \leq t$. We take $(X,\xi)$ the solution of  \eqref{char}, we compute
\[ \left| \frac{\ud }{\ud t} \big( V(X(t)) + \frac{|\xi(t)|^2}{2} \big) \right|
= | -\xi(t).\nabla_x \Phi(t,X(t))| 
\leq \| \nabla_x \Phi(t) \|_{L^\infty(\R^d))}|\xi(t)|. \]
Integrating this expression, we get
\[ V(X(t)) + \frac{|\xi(t)|^2}{2} -E_0 
\leq \int_0^t \| \nabla_x \Phi(s)\|_{L^\infty(\R^d)} |\xi(s)| \ud s. \]
We now set 
$\ds u(t) = \int_0^t \|\nabla_x \Phi(s)\|_{L^\infty(\R^d)} |\xi(s)| \ud s$. 
Since $V \geq 0$, 
\begin{equation}
\label{char-est-int}
\frac{1}{2} \left(\frac{u'(t)}{|\|\nabla_x \Phi(t)\|_{L^\infty(\R^d)} } \right)^2
= \frac{|\xi(t)|^2}{2} 
\leq E_0 + u(t) .
\end{equation} 
It allows us to state
 $ \ds \frac{u'(t)}{(E_0 + u(t))^{1/2}} 
 \leq \sqrt{2}\| \nabla_x \Phi(t)\|_{L^\infty(\R^d)}  $. By integration  on $[0,t]$, we deduce
\[ 2 ( E_0 + u(t) )^{1/2} - 2 (E_0)^{1/2}
 \leq \sqrt{2} \int_0^t \|\nabla_x \Phi(s)\|_{L^\infty(\R^d)} \ud s, \]
and it can be recast as
\[ E_0 + u(t)  \leq \left( E_0^{1/2} + \frac{1}{\sqrt{2}}  
\int_0^t \|\nabla_x \Phi(s)\|_{L^\infty(\R^d)} \ud s \right)^2.\] 
Inserting this expression in \eqref{char-est-int},  we get exactly
\[ |\xi(t)| \leq (2 E_0)^{1/2} 
+ \int_0^t \|\nabla_x \Phi(s)\|_{L^\infty(\R^d)} \ud s. \]
 The estimate on $X$ is deduced directly by integrating this last inequality on $[0,t]$.
\end{ProofOf}

\subsection{Spatially-stationary states}

Let $\rho \in \mathcal{M}^1_+(\R^d)$,  we define the operator
$ \mathcal{S}^{\rho} $ acting from 
$\mathcal{M}^1_+(\R^d\times\R^d)$
to $C(\R ; \mathcal{M}^1_+(\R^d\times\R^d))$
such that $t \mapsto \mathcal{S}^{\rho}_t g_0 $ is the unique solution of 
the linear Vlasov equation:
\begin{equation}
\label{VWRho}
\left\lbrace \begin{array}{l}
\partial_t g + v.\nabla_x g 
= \nabla_x (V-\kappa \Sigma \ast \rho).\nabla_v g \\ \ds
g(0) = g_0.
\end{array} \right.
\end{equation}
Taking $\Phi(t,x) = -\kappa \Sigma \ast \rho(x)$, $\mathcal{S}^\rho$ is well defined when $\Sigma \in C_0^2(\R^d)$, 
$V \in W^{2,\infty}_{loc}(\R^d)$ and $V\geq 0$ thanks to Lemma 
\ref{char-est} and \eqref{rep-car}.

\begin{defin}{Spatially-stationary states}
\label{stat} 

Take $f$ in $\mathcal{M}^1_+(\R^d\times\R^d)$ 
and $\rho = \int_{\R^d_v} \ud f(\cdot,v)$, its spatial density.
\begin{enumerate}
\item We say that $f$ is a stationary state 
if for all $t\in \R$, 
$\mathcal{S}^{\rho}_t f \equiv f$
(in that case
$ v.\nabla_x f - \nabla_x(V-\kappa \Sigma \ast \rho).\nabla_v f = 0$).
\item We say that $f$ is a spatially-stationary state 
(and we write $f \in SEq$) if 
$\int_{\R^d_v} \ud (\mathcal{S}^{\rho}_t f)(v) \equiv \rho$ for all $t \in \R$.
\end{enumerate}
\end{defin}
For all $f \in \mathcal{M}^1_+(\R^d\times\R^d)$, we also 
define the linear energy
\[ \mathscr E (f) 
:= \int_{\R^d\times\R^d} \big( V(x) + \frac{|v|^2}{2} \big) 
\ud f(x,v).\]
We write $K(\mathfrak{m},M)$, the set of all positive measures $f$ of total mass $\mathfrak{m}$ and energy $\mathscr E (f)$ bounded by $M$,
 we set
$SEq(\mathfrak{m}, M)= SEq\cap K(\mathfrak{m},M)$.
A stationary state is always spatially-stationary but the the inclusion can be strict as we will see explicitly in the appendix. However, if one consider only the associated spatial density, the sets are the same in $K(\mathfrak{m},M)$:

\begin{proposition}
\label{semistat}
We suppose $V \in W^{2,\infty}_{loc}(\R^d)$, $V\geq 0$, $\Sigma \in C_0^2(\R^d)$ and $\mathfrak{m}, M>0$
\begin{enumerate}
\item If $f \in SEq(\mathfrak{m},M)$, then  we can find an equilibrium state $\tilde{f} \in K(\mathfrak{m},M)$ such that 
$\rho_{\tilde{f}} \equiv \rho_f$.
\item Take $\Phi \in C(\R_+, W^{2,\infty}(\R^d))$,
there exists a constant $\mathscr H (\mathfrak{m}, M)$ such that any  solution of \eqref{kin}\eqref{rho} $ f \in C_{W}(\R_+, K(\mathfrak{m}, M))$ satisfies
\[ W(\rho_t , \rho_s ) \leq \mathscr H (\mathfrak{m},M)
\left| \int_s^t (W ( f_\tau, SEq(\mathfrak{m}, M) ) )^{1/2}  \ud \tau \right|^{1/2}. \]
\end{enumerate}
\end{proposition}

\begin{ProofOf}{Proposition \ref{semistat}}

{\it Proof of $(1)$}

Take $f_0 \in SEq$ and $\rho_0$ its spatial density. We prove that we can find an equilibrium state $\tilde{f}$ such that $\int_{\R^d_v} \ud \tilde{f}(v) = \rho_0$.
For all $t >0$, we set $ f(t) =  \mathcal{S}^{\rho_0}_t f_0$
and we define the time averages $F_t = \frac{1}{t} \int_0^t f(s) \ud s$. According to the definition of $f$,  
$\int_{\R^d_v} \ud f_s(v) = \rho_0$ for all $s\geq 0$ and then $\int_{\R^d_v} \ud F_t(v) = \rho_0$ for all $t\geq 0$ as well.
It is clear that $(F_t)_{t>0}$ is a family of positive measures with total mass equal to $\mathfrak{m}$, we now prove it's also tight. Since 
$f(t) = \mathcal{S}_t^{\rho_0}f_0$, it solves \eqref{VWRho} and conserves its global energy:
\[\int_{\R^d \times \R^d} 
\big(V(x) +\frac{|v|^2}{2} - \kappa \Sigma \ast \rho_0(x)\big) \ud f_s(x,v)
= \mathscr E(f_0)- \kappa \int_{\R^d}  \Sigma \ast \rho_0(x)\ud \rho_0(x).
\]
Since the spatial density is also constant, we get
\[  
\int_{\R^d\times\R^d} |v|^2 \ud f_s(x,v) 
=  \int_{\R^d\times\R^d} |v|^2 \ud f_0(x,v) \leq 2 \mathscr E(f_0) .
\]
By averaging this identity, it is also satisfied by $F_t$ 
for all $t > 0$. For all $R>0$, we deduce
\[
\int_{|(x,v)|\geq R} \ud F_t(x,v) 
\begin{array}[t]{l} \ds
\leq \int_{|x|\geq R} \ud F_t(x,v) 
+ \int_{|v| \geq R}  \ud F_t(x,v) 
\\ \ds
\leq \int_{|x|\geq R} \ud \rho_0(x)
 + \frac{1}{R^2} \int_{|v| \geq R} |v|^2  \ud F_t(x,v)
\vspace{0.1cm}
\\ \ds
\leq \int_{|x|\geq R} \ud \rho_0(x) + \frac{2 \mathscr E (f_0)}{R^2}
\xrightarrow[ R \rightarrow \infty]{} 0.
\end{array}
\]
Hence $(F_t)_{t>0}$ is compact in 
$\mathcal{M}^1_+(\R^d \times \R^d)$ for the weak topology on 
$C_b(\R^d\times\R^d)$. 
It allows us to take a sequence 
$t_k \xrightarrow[k \rightarrow\infty]{} \infty$ and a limit measure $F$ such that
$\lim_{k\rightarrow \infty} W(F(t_k),F) = 0$.
The convergence is strong enough to ensure
\[\int_{\R^d_v} \ud F(v) 
= \lim_{k \rightarrow \infty} \int_{\R^d_v} \ud F_{t_k}(v) =\rho_0
\quad \mbox{and} \quad 
\mathscr E(F) \leq \lim_{k\rightarrow \infty} \mathscr E(F(t_k)) 
\leq \mathscr E(f_0). \]
We conclude by proving that $F$ is a stationary solution of \eqref{VWRho}.  
For all $\chi $ in $\mathcal{D}(\R^d\times\R^d)$ we have
\[ \begin{array}[t]{l} \ds \left\langle 
\nabla_x (V-\kappa \Sigma \ast \rho_0)
.\nabla_v F - v.\nabla_x F , \chi \right\rangle 
\vspace{0.1cm}
\\ \ds \qquad \qquad = 
- \left\langle 
F, \nabla_x (V-\kappa \Sigma \ast \rho_0 ).\nabla_v\chi  - v.\nabla_x \chi
\right\rangle 
\vspace{0.1cm}
\\ \ds \qquad \qquad = 
\lim_{k\rightarrow \infty} \frac{1}{t_k} \int_0^{t_k} \int_{\R^d\times \R^d}
\left(v.\nabla_x \chi-\nabla_x (V-\kappa \Sigma \ast \rho_0 ).\nabla_v \chi  \right) \ud f_s \ud s 
\vspace{0.1cm}
\\ \ds \qquad \qquad = 
\lim_{k\rightarrow \infty}  \frac{1}{t_k} \int_0^{t_k} 
\frac{\ud}{\ud s} \left\langle   f_s , \chi
\right\rangle  \ud s
\vspace{0.1cm}
\\ \ds \qquad \qquad = 
\lim_{k\rightarrow \infty}  \frac{1}{t_k}
\left(\left\langle \mathcal{S}^{\rho_0}_{t_k}f_0, \chi\right\rangle -  \left\langle f_0,\chi \right\rangle  \right)
\vspace{0.2cm}
\\ \ds \qquad \qquad = 0.
\end{array}\]
Then $F$ is also an equilibrium state. We now prove 
Proposition \ref{semistat}-$2)$, the proof is divided in two steps.

{\it First step: 
Controling $\partial_t \rho$ with $W(f(t),SEq)$:}

Take $f$, a solution of \eqref{kin} in 
$C_{W}(\R_+,K(\mathfrak{m},M))$. By integrating \eqref{kin} with respect to $v$, we get 
\begin{equation}
\label{cons-masse}
\partial_t \rho + \nabla . \left( \int_{\R^d} v \ud f(v) \right)
= 0.
\end{equation}
We set $\vartheta(v) =  1 \wedge [2-|v|]_+$ 
and $\vartheta_R(v) = \vartheta(v/R)$. For all 
$\chi$ in $W^{2,\infty}(\R^d)$ and all 
$g \in SEq(\mathfrak{m},M)$ we make the decomposition
\[ <\partial_t \rho, \chi> 
\begin{array}[t]{l} \ds
= \int_{\R^d\times\R^d} f \nabla_x \chi(x).v \ud x \ud v \\ \ds
= \int_{\R^d\times\R^d}(f -g)\nabla_x\chi(x).v \vartheta_R(v) 
\ud x \ud v 
\\ \ds \qquad
+ \int_{\R^d\times\R^d}(f -g)\nabla_x \chi(x).v (1-\vartheta_R(v))
\ud x \ud v 
\\ \ds
\leq \|f - g\|_{(W^{1,\infty}(\R^d\times \R^d))'} 
\|\vartheta_R \nabla_x\chi.v \|_{W^{1,\infty}(\R^d \times \R^d)}
\\ \ds \qquad
+ \frac{1}{R}  \|\nabla_x \chi\|_{L^\infty(\R^d)}
\int_{|v|\geq R} |f -g| |v|^2 \ud x \ud v 
\\ \ds
\leq  \|f - g\|_{(W^{1,\infty}(\R^d\times\R^d))'} 
(2R + 3) \|\nabla_x \chi\|_{W^{1,\infty}(\R^d)}
\\ \ds \qquad
+ \frac{4}{R} M \|\nabla_x \chi\|_{W^{1,\infty}(\R^d)}
\end{array}\]
Changing $\chi$ for $-\chi$ to get the absolute value and optimizing $R$, we get 
\[ | <\partial_t \rho, \chi> |
\leq \big(8 M^{1/2} + 3 (2\mathfrak{m})^{1/2}\big)
\|f - g \|_{(W^{1,\infty}(\R^d\times\R^d))'} ^{1/2}
\|\nabla_x \chi\|_{W^{1,\infty}(\R^d)}. \]
In the end, optimizing $g$ in $SEq(\mathfrak{m},M)$ and applying the definition of $W$,
we have established
\[ | <\partial_t \rho, \chi> |
\leq u_1 \big(W ( f, SEq(\mathfrak{m},M) ) \big)^{1/2} 
\|\nabla_x \chi\|_{W^{1,\infty}(\R^d)}\]
where we have set $u_1: = (8 M^{1/2} + 3 (2\mathfrak{m})^{1/2})$.

{\it Second step: interpolation between $(W^{2,\infty})'$ and 
$(C_b)'$:}

By integration, we deduce for all $t \geq s \geq 0$ and all 
$\chi$ in $W^{2,\infty}(\R^d)$,
\[ \left| \int_{\R^d} \chi(x) \ud(\rho_s-\rho_t) (x) \right|
\leq  u_1 \|\nabla_x \chi\|_{W^{1,\infty}(\R^d)}
\int_s^t (W ( f_\tau, SEq(\mathfrak{m},M ) )^{1/2}  \ud \tau.\]
We also have directly
$\left| \int_{\R^d} \chi(x) \ud(\rho_s-\rho_t) (x) \right|
\leq 2 \mathfrak{m} \|\chi\|_{L^\infty} $,
we now interpolate between those two inequality.
Take $\ds \theta(x) = \frac{d(d+1)}{|\mathbb{S}^{d-1}|} [1-|x|]_+$, one can check
\[ \int_{\R^d} \theta(x) \ud x = 1, \quad
\| \nabla \theta \|_{L^1(\R^d)} = d+1, \quad 
\int_{\R^d} |x| \theta(x) \ud x = \frac{d}{d+2}. \]
We set $\theta_\epsilon(x) = \frac{1}{\epsilon} \theta(\frac{x}{\epsilon} )$, For any $\chi$ in $W^{1,\infty}(\R^d)$ on the one hand we have
\[ \| \chi - \chi \ast \theta_\epsilon \|_{L^\infty(\R^d)}
\begin{array}[t]{l} \ds
= \sup_{x \in \R^d} \left| 
\int_{\R^d} (\chi(x) - \chi(x + \epsilon y )) \theta(y) \ud y \right| 
\\ \ds
\leq \epsilon \| \nabla \chi \|_{L^\infty(\R^d)} 
\int_{\R^d} |y| \theta(y) \ud y 
\\ \ds 
\leq \epsilon \| \chi \|_{W^{1,\infty}(\R^d)}.  
\end{array}\]
On the other hand, thanks to the Young inequality,
\[ \| \nabla( \chi \ast \theta_\epsilon) \|_{W^{1,\infty}(\R^d)} 
\leq  \| \nabla \theta_\epsilon \|_{L^1(\R^d)} 
\| \chi \|_{W^{1,\infty}(\R^d)}
\leq \frac{d+1}{\epsilon} \|\chi\|_{W^{1,\infty}(\R^d)} .\]
Splitting $\chi$ as 
$ (\chi - \chi \ast \theta_\epsilon) + (\chi \ast \theta_\epsilon) $, 
we deduce
\[ | < \rho_t - \rho_s , \chi > | 
\begin{array}[t]{l} \ds
\leq \frac{d+1}{\epsilon} 
\| \chi \|_{W^{1,\infty}(\R^d)}  
u_1
\int_s^t (W ( f_\tau, SEq^f ) )^{1/2}  \ud \tau 
\\ \ds \qquad
+ 2 \epsilon \mathfrak{m} \| \chi \|_{W^{1,\infty}(\R^d)}. 
\end{array} \]
In the end, we set 
$\mathscr H (\mathfrak{m},M) =\big( \mathfrak{m}(d+1)(8 M^{1/2} + 3 (2\mathfrak{m})^{1/2})\big)^{1/2}$,
taking the suppremum for all $\chi$ such that 
$\| \chi \|_{W^{1,\infty}(\R^d)} \leq 1$ and optimizing $\epsilon$, we get
\[ W(\rho_t , \rho_s ) \leq \mathscr H(\mathfrak{m},M)
\left( \int_s^t (W ( f_\tau, SEq(\mathfrak{m},M) ) )^{1/2}  \ud \tau \right)^{1/2}. \]

\end{ProofOf}

\section{Existence theory and explicit estimates}\label{sec2}

\subsection{Main result}

We now detail the assumptions we will suppose to be fullfield in all the paper.
For the initial data $f_0$, we assume
\[ {\bf (H1)} \quad f_0 \in \mathcal{M}^1_+(\R^d\times\R^d), \quad 
\int_{\R^d \times \R^d} \ud f_0 = \mathfrak{m}, \quad
\mathscr E (f_0) < \infty. \]
Our main motivation to work in this general framework is to get results describing both the solutions of 
\eqref{chm}-\eqref{chminit} and the solution of \eqref{kin}-\eqref{CI}  with an initial data $f_0 \in L^1(\R^d\times\R^d)$. For the well posedness of the characteristic curves, it is natural to ask:
\[ {\bf (H2)} \left\lbrace \begin{array}{l} \ds
\Phi_0 \in C(\R_+ W^{2,\infty} (\R^d)),
\quad p \in L^1_{loc}(\R), \quad \Sigma \in C_0^2(\R^d), 
\\ \ds
V \in W^{2,\infty}_{loc}(\R^d), \quad  V\geq 0.
\end{array} \right.\]
It is already enough to ensure the existence of solution for 
\eqref{kin}-\eqref{CI}. The dissipative behaviour we will describe is due to the following assumption:
\[ {\bf (H3)} \quad  \nabla_x \Phi_0 \in L^1(\R_+, L^{\infty}(\R^d)), \quad p=P',
\quad P \in L^1(\R), \quad \hat{P} \leq 0, \quad \hat{\Sigma} \geq 0.
\]
To our knowledge it is not enough to ensure the uniqueness of the solution. In order to give a uniqueness condition we introduce 
\begin{equation}
\label{Q}
Q(t) =  t \|\nabla_x \Phi_0 \|_{L^1(\R_+, L^\infty(\R^d))}
+ \mathfrak{m} \| \nabla \Sigma \|_{L^\infty} 
\int_0^t \int_0^{s_1} \int_0^{s_2} |p(s_3)| \ud s_3 \ud s_2 \ud s_1
\end{equation}
which is well defined under {\bf (H2)-(H3)} and the 
time-dependant weight
\begin{equation}
\label{w_t}
w_t(x,v) = \exp \left( 
\int_0^t \|\nabla_x^2 V \|_{L^\infty(B(x, t (|v|^2+2V(x))^{1/2} + Q(s))} 
\ud s \right).
\end{equation}
We now set the main result of this section:

\begin{theo}
\label{existence}
Assume {\bf(H1)-(H3)}.
\begin{enumerate}
\item
There exists a measure valued function $t \mapsto f_t$  solving  \eqref{kin}-\eqref{CI} in 

$C_{W}(\R_+; \mathcal{M}^1_+(\R^d\times\R^d))$.
\item
If $\int_{\R^d\times\R^d} w_T(x,v) \ud f_0(x,v) < \infty$, then 
$f$ is unique in $C_{W}(0,T; \mathcal{M}_+^1(\R^d\times\R^d))$
\item
There exists two constants $\mathcal{E}_1, \mathcal{E}_2$ depending on 
$\| \nabla_x \Phi_0\|_{L^1(\R_+, L^\infty(\R^d))}$,
$\mathfrak{m}$,\\ $\|P\|_{L^1(\R)}$, $\| \Sigma \|_{W^{2,\infty}(\R^d)}$ and $\mathscr E(f_0)$
and  such that the solution of \eqref{kin}-\eqref{CI} satisfies 
\end{enumerate}
\[\quad
\begin{array}{ll} \ds
(i) \quad
\int_{\R^d \times \R^d} \big( V(x) + \frac{|v|^2}{2} \big) \ud f_t (x,v) 
\leq \mathcal{E}_1
& \mbox{ for all } t \geq 0
\vspace{0.1cm} \\ \ds (ii) \quad
\int_\R \left| \left|  \int_0^\infty P(t-s) \partial_t (\Sigma \ast \rho) (s) \ud s \right| \right|_{W^{1,\infty}(\R^d)}^2 \ud t 
\leq \mathcal{E}_2 
& 
\end{array} \]
\end{theo}

The bound given by theorem \ref{existence}-$3ii)$ is the key estimate of the whole paper. In the next section it will allows us to prove that 
$ \ds \lim_{t \rightarrow \infty} 
\| \partial_t (\Sigma \ast \rho)(t) \|_{W^{1,\infty}(\R^d)} = 0$.
The rest of the current section is devoted to prove Theorem 
\ref{existence}.

\subsection{Cauchy theory}

We now prove the existence and uniqueness result.

\begin{ProofOf}{Theorem \ref{existence}-1,2)}

{\it Existence}

The proof of existence is classical and it uses some technics that we prefer to detail later in the paper, we only sketch the main points. 
From \cite{BGV}, there already exists a solution under the additional restriction
\begin{equation}
\label{cond-inutile}
 p \in C^1(\R_+),\quad \Sigma \in W^{3,\infty}(\R^d), 
\quad \Phi_0 \in C^1(\R_+; W^{2,\infty}(\R^d)).
\end{equation}
Taking $p,\Sigma$ and $\Phi_0$ satisfying  {\bf(H2)-(H3)} we can find three sequences $(p^k)_{k\geq 0}$, $(\Sigma^k)_{k\geq 0}$ and $(\Phi_0^k)_{k\geq 0}$ satisfying \eqref{cond-inutile} for all $k$ such that $p^k,\Sigma^k,\Phi_0^k$ converge to $p,\Sigma,\Phi_0$
in their respective spaces $L^1_{loc}(\R_+)$, $W^{2,\infty}(\R^d)$,
$C(\R_+ W^{2,\infty} (\R^d))$ when $k$ goes to infinity.
$k$ being fixed it allows us to find 
$f^k \in C_{W}(\R_+,\mathcal{M}_+^1(\R^d\times\R^d))$, a solution of \eqref{kin}-\eqref{CI} for the parameters $p^k,\Sigma^k$ and $\Phi_0^k$. For all $T \geq 0$,
\begin{equation}
\label{borne-Phi}
 \|\Phi^k \|_{C(0,T; W^{2,\infty}(\R^d))}
\leq 
\sup_{k\geq 0} \big( \|\Phi^k_0 \|_{C(0,T; W^{2,\infty}(\R^d))} 
+ \mathfrak{m} \|p^k\|_{L^1(0,T)} \|\Sigma^k\|_{ W^{2,\infty}(\R^d)} \big).
\end{equation}
It allows us to  extract a subsequence (still noted $f^k$) and to find a limit measured valued function $f$ such that 
$(f^k)_{k\geq 0}$ converges to $f$ in
$C_{w\ast}(\R_+,\mathcal{M}_+^1(\R^d\times\R^d))$ (see the beginning of the proof of Lemma \ref{ens-omega-limites}). Combining 
\eqref{borne-Phi}, Lemma \ref{char-est} and  \eqref{rep-car},
$(f^k)_{ k\geq 0 }$ is uniformly tight on $[0,T]$ for all $T\geq 0$. It allows us to deduce that the convergence also holds in 
$C_{W}(\R_+,\mathcal{M}_+^1(\R^d\times\R^d))$. From the definition of $W$, taking $\chi(y,v) = \nabla \Sigma(x-y)$, we get 
$\lim_{k\rightarrow \infty } 
\| \nabla \Sigma^k \ast \rho^k - \nabla \Sigma \ast \rho_f \|_{L^\infty([0,T] \times\R^d)}
= 0$ and it is enough to pass to the limit in \eqref{Phi}:
\[ \lim_{k\rightarrow \infty } \sup_{0\leq t \leq T}
\left| \left| \nabla_x \Phi^k(t) -\nabla_x \big(\Phi_0(t) 
- \int_0^t (\Sigma \ast \rho)(s) p(t-s) \ud s \big)
\right| \right|_{L^\infty(\R^d)} =0. \]
Passing to the limit in 
\eqref{kin} as well, we deduce that $f$ solves 
\eqref{kin}-\eqref{CI}.

{\it Uniqueness}

Take $f_1$ and $f_2$ two solutions of \eqref{kin}-\eqref{CI}. 
For $i=1,2$ we write $\rho_i = \int_{\R^d_v} \ud f_i(v)$,
$\Phi_i(t) = \Phi_0(t) - \int_0^t (\Sigma \ast \rho_i)(s) p(t-s) \ud s$, and $\varphi^{\Phi_i}$ the flow associated to \eqref{char} for $\Phi = \Phi_i$.
Setting $h = f_1-f_2$, $h$ solves the linear equation
\[ \left\lbrace \begin{array}{l} \ds
\partial_t h + v. \nabla_x h -\nabla_x(V+ \Phi_1).\nabla_v h 
=-\nabla_v f_2. \nabla_x 
\left( \int_0^t (\Sigma \ast \rho_h)(s) p(t-s) \ud s \right),
\\ \ds
h(0) = 0.
\end{array} \right.\]
Thanks to \eqref{rep-car} and the Duhamel formula, we get for all $T>0$,
\begin{equation}
\label{duhamel-h}
h(T) = -\int_0^T {\varphi_t^{\Phi_1,T}}_{\#} \left(\nabla_v f_2. \nabla_x 
\left( \int_0^t (\Sigma \ast \rho_h)(s) p(t-s) \ud s \right) \right) \ud t.
\end{equation}
The uniqueness of the solution will be provided by an estimates on $h$ for the strong norm on the dual of 
$W^{1,\infty}(\R^d\times\R^d)$. Picking $\chi$ a test function, 
and applying \eqref{rep-car} again, we get from \eqref{duhamel-h} 
\begin{equation}
\label{dual-unic}
\begin{array}{l} \ds 
\int_{\R^d\times\R^d} \chi(x,v) \ud h_T(x,v) 
\\ \ds 
=-\int_0^T \int_{\R^d\times\R^d} 
\nabla_v f_2. \nabla_x 
\left(\int_0^t (\Sigma \ast \rho_h)(s) p(t-s) \ud s \right) 
\chi \circ \varphi^{\Phi_1,T}_t \ud t
\\ \ds
= \int_0^T \int_{\R^d\times\R^d}
\nabla_x
\left(\int_0^t (\Sigma \ast \rho_h)(s) p(t-s) \ud s \right).
\nabla_v(\chi \circ \varphi^{\Phi_1,T}_t) 
\circ \varphi^{\Phi_2,t}_0 \ud f_0(x,v) \ud t.
\end{array}
\end{equation}
In order to conclude  the main point is to establish a satisfying estimate on
$u(t,T,x,v) = \nabla_v ( \chi \circ \varphi^{\Phi_1,T}_t) 
(\varphi^{\Phi_2,t}_0(x,v))$. We have directly
\[
|\nabla_v (\chi \circ \varphi^{\Phi_1,T}_t) (x,v)|
\leq \|\chi\|_{W^{1,\infty}(\R^d\times\R^d)} 
\left|D (\varphi^{\Phi_1,T}_t)(x,v) \right|\] 
where the differential of the flow 
$D \varphi^{\Phi_1,T}_{t}(x,v).(y,w)$ is given by the solution at time $\tau =T$ of the following ordinary differential equation
\[
\left\lbrace \begin{array}{l}
\ds \frac{\ud}{\ud \tau} \mathcal{Y}(\tau) = \mathcal{W} (\tau), \vspace{0.1cm}\\ \ds 
\frac{\ud}{\ud \tau} \mathcal{W} (\tau) = -\nabla^2_x(V+\Phi_1(\tau))(\varphi^{\Phi_1,\tau}_{t,x}(x,v)) \mathcal{Y}(\tau),
\end{array} \right. 
\qquad 
\left\lbrace \begin{array}{l}
\mathcal{Y}(t) = y, \vspace{0.6cm}\\ \mathcal{W} (t)=w.
\end{array} \right. 
\]
Thanks to the gronwall lemma, it allows us to get 
\[\left| D (\varphi^{\Phi_1,T}_{t})(x,v).(y,w) \right|
\leq |(y,w)| \exp \Big( \int_t^T (1+ 
|\nabla^2_x(V+\Phi_1(\tau))(\varphi^{\Phi_1,\tau}_{t,x}(x,v))|) \ud \tau
\Big).\]
From the expression of $u$, we deduce
\begin{equation}
\label{controle-u-1}
\frac{|u(t,T,x,v)|}{\|\chi\|_{W^{1,\infty}(\R^d\times\R^d)} }
\begin{array}[t]{l} \ds \leq
\exp \Big( \int_t^T (1+ 
|\nabla^2_x(V+\Phi_1(\tau))(\varphi^{\Phi_1,\tau}_{t,x}(\varphi^{\Phi_2,t}_0(x,v)))|) \ud \tau
\Big)
\\ \ds
\leq A_1(T) \exp \Big(  \int_t^T  
|\nabla^2_x V(\varphi^{\Phi_1,\tau}_{t,x}(\varphi^{\Phi_2,t}_0(x,v)))| \ud \tau
\Big)
\end{array}
\end{equation}
where we have set 
\[A_1(T) =  \exp \Big( 
\int_0^T (1+\|\nabla^2 \Phi_0(\tau) \|_{L^\infty(\R^d)}
+ \mathfrak{m} \| \nabla^2 \Sigma\|_{L^\infty(\R^d)}
\|p\|_{L^1(0,\tau)})\ud \tau  \Big). \]
Setting $\Phi(\tau) = \Phi_2(\tau)$  for $\tau <t$ and 
$\Phi(\tau) = \Phi_1(\tau)$  for $\tau \geq t$, by composition of the flow
$\varphi^{\Phi_1,\tau}_{t,x}(\varphi^{\Phi_2,t}_0(x,v))
= \varphi^{\Phi,\tau}_{0,x}(x,v)$ for all $\tau \geq t$.
Applying Lemma \ref{char-est} and remembering of \eqref{Q} we get for all $(x,v)$
\[|\varphi^{\Phi,\tau}_{0,x}(x,v)-x| 
\begin{array}[t]{l} \ds 
\leq
 (|v|^2+2V(x))^{1/2} \tau
+ \int_0^\tau \int_0^{\tau_1} \|\nabla_x \Phi(\tau_2)\|_{L^\infty(\R^d)} \ud \tau_2 \ud \tau_1
\\ \ds
\leq
 (|v|^2+2V(x))^{1/2} \tau
+ Q(\tau)
\end{array}
\]
Coming back to \eqref{controle-u-1} and \eqref{w_t}, we deduce
\begin{equation}
\label{controle-u}
\frac{|u(t,T,x,v)|}{\|\chi\|_{W^{1,\infty}(\R^d\times\R^d)} }
\begin{array}[t]{l} \ds 
\leq A_1(T) \exp \Big(  \int_t^T  
\|\nabla^2_x V\|_{L^\infty(B(x,(|v|^2+2V(x))^{1/2} \tau
+ Q(\tau))} \ud \tau
\Big)
\\ \ds
\leq A_1(T) w_T(x,v).
\end{array}
\end{equation}
Coming back to \eqref{dual-unic}, we get the simple controle
\[ \begin{array}{l} \ds
\|h(T)\|_{W^{1,\infty}(\R^d\times\R^d)'} 
\\ \ds \qquad
\leq A_1(T) \int_0^T \int_{\R^d\times\R^d}
\left| \left| 
\int_0^t \nabla_x (\Sigma \ast \rho_h)(s) p(t-s) \ud s 
\right| \right|_{L^\infty} w_T(x,v) \ud f_0(x,v) \ud t 
\\ \ds \qquad
\leq A_1(T) \Big(\int_{\R^d\times\R^d} w_T(x,v) \ud f_0(x,v) \Big)
\int_0^T  \| \nabla_x \Sigma \ast \rho_h\|_{L^\infty([0,t] \times\R^d)} \|p\|_{L^1(0,t)} \ud t 
\\ \ds \qquad
\leq A_2(T) \int_0^T 
\|h\|_{L^\infty \left([0,t] ; W^{1,\infty}(\R^d\times\R^d)'\right)}
\ud t. 
\end{array}
\]
where we have set
\[ A_2(T) = A_1(T) 
\left(\int_{\R^d\times\R^d} w_T(x,v) \ud f_0(x,v) \right) 
\|\nabla \Sigma\|_{W^{1,\infty}(\R^d)} \|p\|_{L^1(0,T)} \]
and we point out that $A_2(T)$ is well defined due to the additional integrability assumption on $f_0$. By construction, $A_2$ increases. For $0\leq t \leq T$, we deduce 
\[ 
\|h\|_{L^\infty \left([0,t] ; W^{1,\infty}(\R^d\times\R^d)'\right)}
\leq A_2(T) \int_0^t
\|h\|_{L^\infty \left([0,s] ; W^{1,\infty}(\R^d\times\R^d)'\right)}
\ud s.\]
By applying the gronwall lemma, we conclude that $h \equiv 0$ on
$[0,T]$.

\end{ProofOf}

\subsection{Explicit estimates}

In order to prove the second part of Theorem \ref{existence}, we introduce the following functional
\begin{equation}
\label{E}
\tilde{\mathscr E} (f) 
\begin{array}[t]{l} \ds
= \int_{\R^d\times\R^d} 
\big( V(x) - \frac{\kappa}{2} \Sigma \ast \rho_f(x) +\frac{|v|^2}{2} \big)
\ud f(x,v) \\ \ds
= \mathscr E(f) - \frac{\kappa}{2} \int_{\R^d} 
(\Sigma \ast \rho_f)(x) \ud \rho_f (x)
\end{array}
\end{equation}
which represents the energy for the Vlasov equation with the self potential $-\kappa \Sigma$. It  appears naturally in our analysis since in the end we will prove that $f_t$ gets closer and closer to a family of particular solutions of this equation when $t$ goes to infinity.
Thanks to the positivity of $V$ and the mass consevation, 
$\tilde{\mathscr E}$ is bounded from below :
\begin{equation}
\label{minE}
 \tilde{\mathscr E} (f(t))
+\frac{\kappa}{2} \mathfrak{m}^2 \|\Sigma\|_{L^\infty(\R^d)} 
\geq  
\int_{\R^d\times\R^d} \left( V(x) +\frac{|v|^2}{2} \right) \ud f_t(x,v)
\geq 0.
\end{equation}
Taking $h$ defined on $\R_+$, we also introduce the following notation:
\begin{equation}
\label{hT}
[h]_T(t) = \left\lbrace \begin{array}{ll} 
0 & t<0, \\
h(t) & 0\leq t \leq T, \\
0 & t>T.
\end{array}\right. 
\end{equation}
Before proving the estimates, we point out that an easier variant of the proof can be established with the restrictions $\Sigma = \sigma_1\ast \sigma_1$ and $P= -\check{q}_1\ast q_1$ satified when $p$ and $\Sigma$ are defined by \eqref{reduct-param}.
In that case, from computations which will be clear later we get simply:
\begin{equation}
\label{sdb}
\tilde{\mathscr E}(f(T)) - \tilde{\mathscr E}(f_0)
= \int_0^T \int_{\R^d} A_0 \partial_t \rho(t) \ud x \ud t
- \frac{1}{2} \int_{(0,T) \times \R^d}
\big(q_1 \ast [\partial_t \sigma_1 \ast \rho]_T\big)^2(t,x)\ud t \ud x  
\end{equation}
where $A_0 \in L^1(\R_+,W^{1,\infty}(\R^d))$. Our methods allow us to deduce that the three terms involved in \eqref{sdb} are uniformly bounded with respect to $T>0$. Theorem \ref{existence}-$3i)$ follows directly from \eqref{minE} and the bound on the right hand side of \eqref{sdb} allows us to find an explicit constant $\tilde{\mathcal{E}}_2$ such that
\[\int_\R \left| \left|  \int_0^\infty q_1(t-s) \partial_t (\Sigma \ast \rho) (s) \ud s \right| \right|_{W^{1,\infty}(\R^d)}^2 \ud t 
\leq \tilde{\mathcal{E}}_2 \]
which is slightly stronger than Theorem \ref{existence}-$3ii)$. Under our more general hypothesis, the right-hand side of \eqref{sdb} is replaced by a term which is still non negative but does not present like the integral of a square. The result holds anyway thanks to the following inequality

\begin{lemma}
\label{controle}
Take $k\geq 0$ and $\Sigma \in C_0^{2k}(\R^d)$ such that 
$\hat{\Sigma} \geq 0$. Then $(1+|\cdot |^{2k}) \hat{\Sigma}$ is a finite measure on $\R^d$. 
Setting $\Lambda_k= \sum_{j=0}^k (-\Delta)^j\Sigma(0)$,  
any  $u$ in $C^k_0(\R^d)'$ satisfies
\[ \|  \Sigma \ast u \|_{W^{k,\infty}(\R^d)}^2 
\leq \Lambda_k \big\langle u ; \Sigma \ast u 
\big\rangle_{(C_0^k)'\times C_0^k}
. \]
\end{lemma}


\begin{rmk}
\label{pot}
The arguments involved in the proof of Lemma \ref{controle} also allow us to deduce that $\hat{P} \in L^1(\R)$, therefore 
$\kappa := -P(0) = -\frac{1}{2\pi} \int_{-\infty}^{+\infty} \hat{P}(\omega) \ud \omega$ is positive. By the same way, from Lemma \ref{controle}, we also deduce 
$\Sigma(0) > \Sigma(x)$ for all $x\neq 0$. The self potential 
$-\kappa \Sigma$ involved in the limit equation \eqref{eq-lim} is never repulsive since
$-\kappa \Sigma(0)< -\kappa \Sigma(x_1-x_2)$ for all 
$x_1\neq x_2$ in $\R^d$.
\end{rmk}
We now end the proof of Theorem \ref{existence}, keeping the proof of Lemma \ref{controle} for the end.

\begin{ProofOf}{Theorem \ref{existence}-$3)$}

{\it First step: bounding $\tilde{\mathscr E}$}

We rewrite the expression of the self potential:
\begin{equation}
\label{devPhi}
\Phi(t,x) \begin{array}[t]{l} \ds
= \Phi_0(t,x) - \int_0^t p(t-s) \Sigma \ast \rho(s) \ud s \\ \ds
= \Phi_0(t,x) - P(t) \Sigma \ast \rho(0) + P(0) \Sigma \ast \rho(t) 
- \int_0^t P(t-s) \partial_t( \Sigma \ast \rho)(s) \ud s \\ \ds
= A_0(t,x) -\kappa \Sigma \ast \rho(t) 
- \int_0^t P(t-s) \partial_t( \Sigma \ast \rho)(s) \ud s
\end{array}
\end{equation}
where $A_0$ satisfies
$ \|\nabla_x A_0(t) \|_{L^\infty(\R^d)} 
\leq  \| \nabla_x \Phi_0(t) \|_{L^\infty(\R^d)} 
+ \mathfrak{m} \| \Sigma \|_{W^{1,\infty}(\R^d)} |P(t)| $, so that  setting
$K_1:=  \| \nabla_x \Phi_0 \|_{L^1(\R_+, L^\infty(\R^d))} 
+ \mathfrak{m} \| \Sigma \|_{W^{1,\infty}(\R^d)} \|P\|_{L^1(\R)} $, we have
\begin{equation}
\label{A0}
\int_0^\infty \| \nabla_x A_0(t) \|_{L^\infty(\R^d)}  \ud t \leq K_1 < \infty
\end{equation}

Using \eqref{devPhi}, we compute formally the time derivative of $\mathscr E$
\begin{equation}
\label{dtE1}
\frac{\ud}{\ud t} \tilde{\mathscr E} (f) 
\begin{array}[t]{l} \ds
=\int_{\R^d} \partial_t \rho ( V - \kappa  \Sigma \ast \rho_f ) \ud x 
- \int_{\R^d} \partial_t \rho ( V + \Phi(t)) \ud x \\ \ds
= \int_{\R^d} \partial_t \rho \left( -A_0(t) +
\int_0^t P(t-s) \partial_t( \Sigma \ast \rho)(s) \ud s \right) \ud x \\ \ds
= -\int_{\R^d}  A_0(t)\partial_t \rho \ud x
+\int_0^t \int_{\R^d} 
\partial_t\rho (t) \partial_t( \Sigma \ast \rho )(s)
P(t-s) \ud x \ud s.
\end{array}
\end{equation}
We start by estimating the first term,
thanks to \eqref{minE} and the Cauchy-Shwartz inequality:
\begin{equation}
\label{mom1}
\int_{\R^d\times \R^d} |v| \ud f_t (x,v) 
\begin{array}[t]{l} \ds
\leq\left(\int_{\R^d\times \R^d} \ud f_t(x,v) \right)^{1/2}
\left(\int_{\R^d\times \R^d} |v|^2 \ud f_t (x,v) \right)^{1/2} 
\\ \ds
\leq \mathfrak{m}^{1/2}  \left( \tilde{\mathscr E}(f(t))
 +\frac{\kappa}{2} \mathfrak{m}^2 \| \Sigma \|_{L^{\infty}(\R^d)}\right)^{1/2}.
 \end{array} 
\end{equation}
From \eqref{cons-masse}, we deduce
\begin{equation}
\label{est-init} 
\int_{\R^d} \partial_t \rho A_0(t) \ud x 
\leq \|\nabla_x A_0(t) \|_{L^\infty(\R^d)} 
\mathfrak{m}^{1/2}  \left( \tilde{\mathscr E}(f(t))
+\frac{\kappa}{2} \mathfrak{m}^2 \| \Sigma \|_{L^{\infty}(\R^d)} \right)^{1/2}.
\end{equation}
Integrating \eqref{dtE1} and using \eqref{A0}, we get
\begin{equation}
\label{dtE2}
\tilde{\mathscr E} (f(T)) -\tilde{\mathscr E} (f_0) 
\begin{array}[t]{l} \ds
\leq \int_0^T  \|\nabla_x A_0(t) \|_{L^\infty(\R^d)}
\mathfrak{m}^{1/2}  \left( \tilde{\mathscr E}(f(t))
+\frac{\kappa}{2} \mathfrak{m}^2 \| \Sigma \|_{L^{\infty}(\R^d)} \right)^{1/2} 
\ud t
\\ \ds \qquad +k(T) 
\\ \ds
\leq  K_1 \mathfrak{m}^{1/2} \sup_{0\leq t \leq T} 
\left(\tilde{\mathscr E}(f(t))+
\frac{\kappa}{2} \mathfrak{m}^2 \| \Sigma \|_{L^{\infty}(\R^d)}\right)^{1/2}  +k(T)
\end{array} 
\end{equation}
where $k$ is the integral of the right hand side of \eqref{dtE1}. Thanks to the parity of $P$, its expression can be transformed:
\begin{equation}
\label{k}
k(T) \begin{array}[t]{l} \ds
=\int_0^T \int_0^t \int_{\R^d} 
\partial_t\rho (t) \Sigma \ast \partial_t \rho (s)
P(t-s) \ud x \ud s \ud t \\ \ds
=\frac{1}{2} \int_{[0,T]^2\times \R^d} 
\partial_t\rho (t) \Sigma \ast \partial_t \rho (s)
P(t-s) \ud x \ud s \ud t \\ \ds
= \frac{1}{2} \int_{ \R \times \R^d} 
[\partial_t \rho ]_T(t) 
(P \ast \Sigma \ast [\partial_t \rho ]_T)(t)
\ud t \ud x 
\\ \ds
= \frac{1}{2(2\pi)^{d+1}} \int_{\R \times \R^d} 
|\mathcal{F}_{t,x}([\partial_t \rho ]_T)|^2(\omega,\xi) 
\hat{P}(\omega) \hat{\Sigma}(\xi) \ud \xi \ud \omega
\\ \ds
\leq 0.
\end{array} 
\end{equation}
We now set 
$\ds u(t) = \sup_{0\leq s \leq t} 
\left(\tilde{\mathscr E}(f(s))+\frac{\kappa}{2} \mathfrak{m}^2 
 \| \Sigma \|_{L^{\infty}(\R^d)}^2\right)^{1/2}$, 
from \eqref{dtE2} and \eqref{k}, we get
\[\tilde{\mathscr E}(f(T))+\frac{\kappa}{2} \mathfrak{m}^2 
 \| \Sigma \|_{L^{\infty}(\R^d)}
\leq \tilde{\mathscr E}(f_0)+\frac{\kappa}{2} \mathfrak{m}^2  
\| \Sigma \|_{L^{\infty}(\R^d)}
+ K_1 \mathfrak{m}^{1/2} u(T) \]
and since the right hand side increases, taking the supremum for all 
$0\leq t \leq T$, we get
\[u(T)^2 \leq u(0)^2 + K_1
\mathfrak{m}^{1/2} u(T)\]
and it allows us to deduce the uniform bound
\[ u(T) \leq \frac{1}{2} K_1 \mathfrak{m}^{1/2} +
\left( u(0)^2 + \frac{1}{4} K_1^2 \mathfrak{m} \right)^{1/2} =:K_2 .\]
Coming back to the definition of $u$, we have proved
\[ \tilde{\mathscr E}(f) (T) \leq 
 K_2^2 
- \frac{\kappa}{2} \mathfrak{m}^2  \| \Sigma \|_{L^{\infty}(\R^d)}. \]

{\it Second step: deducing $(i)$ and $(ii)$.}

From the definition of $\tilde{\mathscr E}$, it is clear that 
$(i)$ is satisfied for $\mathcal{E}_1 = K_2^2$. Coming back to \eqref{dtE2}, and \eqref{minE}, we also deduce for all $T\geq 0$,
\begin{equation}
\label{k2}
-k(T)
\leq K_1 \mathfrak{m}^{1/2} K_2 + \tilde{\mathscr E}(f_0) 
+  \frac{\kappa}{2} \mathfrak{m}^2  \| \Sigma \|_{L^{\infty}(\R^d)} =:K_3.
\end{equation}

We just have to deduce $(ii)$. Applying
Lemma \ref{controle}, we first have
\begin{equation}
\label{k3}
\begin{array}{l} \ds
\int_\R \| (P \ast \Sigma \ast \left[ \partial_t \rho \right]_T) (t)
\|_{W^{1,\infty}(\R^d)}^2\ud t 
\\ \qquad \qquad \qquad \ds
\leq \Lambda_1 \int_\R \int_{\R^d} 
(P\ast [\partial_t \rho ]_T) \Sigma \ast (P\ast [\partial_t \rho ]_T) 
\ud x \ud t 
\\ \qquad \qquad \qquad \ds
\leq \frac{\Lambda_1}{(2\pi)^{d+1}} \int_{\R \times\R^d}
|\mathcal{F}_{t,x}([\partial_t \rho ]_T)(\omega,\xi) |^2
|\hat{P}(\omega)|^2 \hat{\Sigma}(\xi) \ud \xi \ud \omega
\\ \qquad \qquad \qquad \ds
\leq -\frac{\Lambda_1}{(2\pi)^{d+1}} \|\hat{P} \|_{L^\infty(\R)}
\int_{\R \times\R^d} 
|\mathcal{F}_{t,x}([\partial_t \rho )]_T)(\omega,\xi) |^2
\hat{P}(\omega) \hat{\Sigma}(\xi) \ud \xi \ud \omega
\\ \qquad \qquad \qquad \ds
\leq 2 \Lambda_1 \|P\|_{L^1(\R)} K_3
\end{array}
\end{equation}
where we have recognized the last expression of $k$  given in \eqref{k} and applied \eqref{k2}.
The computations are rigourous because 
$|\cdot|^2 \hat{\Sigma}$ is a finite
measure thanks to Lemma \ref{controle} while 
$\widehat{\partial_t \rho}(t,\xi) = i\xi.\hat{J}(t,\xi)$ and 
$t \mapsto \|J(t)\|_{L^1(\R^d)}$ is bounded thanks to $(i)$
and \eqref{mom1}. 
In order to conclude, we just have to let $T$ go to infinity.
Applying Young inequality, thanks to \eqref{cons-masse}, \eqref{mom1} and $(i)$, for any $t \geq 0$, we have 
\begin{equation}
\label{dtsigmarho}
\| (\Sigma \ast \partial_t \rho)(t) \|_{W^{1,\infty}(\R^d)}
\leq \|\nabla \Sigma \|_{W^{1,\infty}(\R^d)} \|f_t v \|_{L^1(\R^d\times\R^d)}
\leq \|\Sigma \|_{W^{2,\infty}(\R^d)} 
(2 \mathfrak{m} \mathcal{E}_1)^{1/2}.
\end{equation}
Foll all $t,T\geq 0$, we deduce punctually
\[ \begin{array}{l} \ds
\left| \left| 
\left(P \ast [\partial_t \Sigma \ast (\rho )]_T \right)(t)
-\left(P \ast [\partial_t( \Sigma \ast \rho )]_\infty \right)(t)
\right| \right|_{W^{1,\infty}(\R^d)} \\ \ds \qquad
= \left| \left| 
\int_T^\infty \partial_t( \Sigma \ast \rho )(s) P(t-s) \ud s 
\right| \right|_{W^{1,\infty}(\R^d)} \\ \ds \qquad \qquad
\leq  
\|\Sigma \|_{W^{2,\infty}(\R^d)} 
(2 \mathfrak{m} \mathcal{E}_1)^{1/2}
\int_{-\infty}^{t-T} |P(s) | \ud s 
\xrightarrow[T \rightarrow +\infty]{} 0.
\end{array}
\]
When $T$ goes to infinity, it allows us to deduce from 
\eqref{k3}
\[\begin{array}{l} \ds
\int_\R \left| \left|  \int_0^\infty P(t-s) \partial_t (\Sigma \ast \rho) (s) \ud s \right| \right|_{W^{1,\infty}(\R^d)}^2 \ud t 
\\ \qquad \qquad \ds
=\int_\R \| P \ast [\partial_t( \Sigma \ast \rho )]_\infty(t)
\|_{W^{1,\infty}(\R^d)}^2 \ud t
\\ \qquad \qquad \ds
\leq   \liminf_{T\rightarrow +\infty} \int_{\R} 
\| P\ast [\partial_t( \Sigma \ast \rho ]_T(t) \|_{W^{1,\infty}(\R^d)}^2 \ud t
\\ \qquad \qquad \ds
\leq 2 \Lambda_1 \|P\|_{L^1(\R)} K_3 =:\mathcal{E}_2.
\end{array}
\]

\end{ProofOf}

We end this section by proving Lemma \ref{controle}.

\begin{ProofOf}{Lemma \ref{controle}}

{\it Integrability of $\hat{\Sigma}$}

Since $\hat{\Sigma} \geq 0$, $\hat{\Sigma}$ is already a locally finite non negative measure on $\R^d$. Taking $\theta$ such that 
$\hat{\theta} \in C_c^\infty(\R^d)$, $\hat{\theta} \geq 0$ and 
$\hat{\theta}(0) = 1$, we set 
$\theta_\epsilon(x) = \frac{1}{\epsilon^d} 
\theta \big( \frac{x}{\epsilon} \big)$ and 
$\Sigma_\epsilon = \Sigma \ast \theta_\epsilon$. 
According to our hypothesis $\theta \in \mathscr S (\R^d)$, $\int_{\R^d} \theta(x) \ud x = 1$ and $\Sigma \in C_0^{2k}(\R^d)$, therefore
$\lim_{\epsilon \rightarrow 0} 
\| \Sigma_\epsilon - \Sigma \|_{W^{2k,\infty}(\R^d)} = 0$.
For any $j \in \{0,..,k\}$, since $\hat{\theta}$ is compactly supported, the Fourier transform of 
$(-\Delta)^j\Sigma_\epsilon$ explicitely given by 
$\widehat{(-\Delta)^j\Sigma_\epsilon} = |\cdot|^{2j}\hat{\theta}(\epsilon \cdot) \hat{\Sigma}$ is 
a finite measure on $\R^d$. Applying the inverse Fourier transform, it allows us to set
\[ \frac{1}{(2\pi)^d}
\int_{\R^d} |\xi|^{2j} \hat{\theta}(\epsilon \xi) \ud \hat{\Sigma}(\xi)
= (-\Delta)^j \Sigma_\epsilon (0) 
\xrightarrow[\epsilon \rightarrow 0]{} (-\Delta)^j \Sigma (0).  \]
From the Fatou's lemma, we deduce that
$|\cdot|^{2j} \hat{\Sigma}$ is also a finite measure on $\R^d$. 
Since it is also the Fourier transform of $(-\Delta)^j \Sigma$,
applying the inverse Fourier transform again we have 
\begin{equation}
\label{mom-Sigma}
\int_{\R^d} |\xi|^{2j} \ud \hat{\Sigma}(\xi)
= (2\pi)^d (-\Delta)^j \Sigma (0) \quad \mbox{ for all } 
j \in \{0,..,k\}.
\end{equation}

{\it Proving the inequality for a whole sequence approximating $\Sigma$}

The key point is to approximate $\Sigma$ by a sequence 
$\sigma^\epsilon \ast \sigma^\epsilon$ that we define now.
We now take $\vartheta$ such that 
$\hat{\vartheta} \in C_c^\infty(\R^d)$, $\hat{\vartheta} \geq 0$
and $\int_{\R^d} \hat{\vartheta}(\xi) \ud \xi = (2\pi)^d$. We set $\vartheta_\epsilon (x) = \vartheta(\epsilon x)$
and $\Sigma^\epsilon = \Sigma \vartheta_\epsilon$. 
According to our hypothesis, $\vartheta \in \mathscr S(\R^d)$,
$\vartheta(0) =1$ and $\Sigma \in C^{2k}_0(\R^d)$, therefore 
$\lim_{\epsilon \rightarrow 0} 
\| \Sigma^\epsilon - \Sigma \|_{W^{2k,\infty}(\R^d)} = 0$. The Fourier transform of $\Sigma^\epsilon$ is now given by 
$\hat{\Sigma}^\epsilon
= \frac{1}{(2\pi)^d} \hat{\Sigma}  \ast \hat{\vartheta}_\epsilon$ therefore it is non negative and \eqref{mom-Sigma} allows us to deduce 
for any $j \in \{0,..,k\}$
\begin{equation}
\label{mom-Sigma-e}\int_{\R^d} |\xi|^{2j} (\hat{\Sigma}  \ast \hat{\vartheta}_\epsilon)(\xi) \ud \xi
= (2\pi)^d (-\Delta)^j \Sigma^\epsilon (0)
\xrightarrow[\epsilon \rightarrow 0]{} 
(2\pi)^d (-\Delta)^j \Sigma (0).  
\end{equation}
In particular for $j=0$, since 
$ \frac{1}{(2\pi)^d} \hat{\Sigma}  \ast \hat{\vartheta_\epsilon} \in L^1(\R^d)$, it allows us to define 
$\sigma^\epsilon$  in $L^2(\R^d)$ by its Fourier transform 
$\hat{\sigma^\epsilon} = (\frac{1}{(2\pi)^d} \hat{\Sigma}  \ast \vartheta_\epsilon)^{1/2}$.
For any $\epsilon$, 
$\Sigma^\epsilon = \sigma^\epsilon \ast \sigma^\epsilon$ and for any $j \in \{0,..,k\}$, the Plancherel theorem allows us to deduce from
\eqref{mom-Sigma-e} the uniform bound
\begin{equation}
\label{mom-sigma-e}
\lim_{\epsilon \rightarrow 0} \| \nabla^j \sigma^\epsilon \|_{L^2(\R^d)}^2
=\lim_{\epsilon \rightarrow 0} \frac{1}{(2\pi)^d}
\int_{\R^d} |\xi|^{2j} |\sigma^\epsilon(\xi)|^2 \ud \xi 
= (-\Delta)^j \Sigma (0).
\end{equation}
We first consider $u$ in $L^1(\R^d)$, for any $\epsilon>0$ and any $j \in \{1,..,k\}$, the Young inequality allows us to get
\[ \| \nabla^j \Sigma^\epsilon \ast u \|_{L^\infty(\R^d)}
=\| (\nabla^j \sigma^\epsilon) 
\ast (\sigma^\epsilon \ast u )\|_{L^\infty(\R^d)}
\leq \| \nabla^j \sigma^\epsilon \|_{L^2(\R^d)}
\|\sigma^\epsilon \ast u \|_{L^2(\R^d)}
\]
Since $\hat{\sigma^\epsilon} \geq 0$, $\sigma^\epsilon$ is even. It allows us to rewrite the right hand side in terms of $\Sigma^\epsilon$:
\[ \| \sigma^\epsilon \ast u \|_{L^2( \R^d)}^2
= \int_{\R^d} (\sigma^\epsilon \ast u )^2 \ud x
= \int_{\R^d} u (\sigma^\epsilon \ast \sigma^\epsilon )  \ast u \ud x 
= \int_{\R^d} u \Sigma^\epsilon \ast u \ud x .
\] 

Taking $\|h\|_{W^{k,\infty}(\R^d)} 
=(\|h\|_{L^\infty(\R^d)}^2 + \dots + \|\nabla^k h \|_{L^\infty(\R^d)}^2 )^{1/2}$ we have established
\[ \| \Sigma^\epsilon \ast u \|_{W^{k,\infty}(\R^d)}^2
\leq 
\big( \sum_{j=0}^k \| \nabla^j \sigma^\epsilon \|_{L^2(\R^d)}^2 \big)
\int_{\R^d} u \Sigma^\epsilon \ast u \ud x .
\] 
    
{\it Conclusion}

We now let $\epsilon$ go to $0$. Since 
$\Sigma^\epsilon$ goes to $\Sigma$ in $W^{2,\infty}(\R^d)$, 
we already have
$\lim_{\epsilon \rightarrow 0 }  
\| \Sigma^\epsilon \ast u \|_{W^{k,\infty}(\R^d)} 
= \| \Sigma \ast u \|_{W^{k,\infty}(\R^d)}  $
while according to \eqref{mom-sigma-e}, the constant does not explode:
$ \lim_{\epsilon \rightarrow 0 }
\sum_{j=0}^k \| \nabla^j \sigma^\epsilon \|_{L^2(\R^d)}^2 
= \sum_{j=0}^k (-\Delta)^j \Sigma (0)=:\Lambda_k$. For all
$u \in L^1(\R^d)$ we have proved 
\begin{equation}
\label{presque}
\| \Sigma \ast u \|_{W^{k,\infty}(\R^d)}^2
\leq \Lambda_k \int_{\R^d} u \Sigma \ast u \ud x
= \Lambda_k \big\langle u ; \Sigma \ast u \big\rangle_{(C_0^k)'\times C_0^k} .
\end{equation}
Taking now $u \in C^k_0(\R^d)'$, 
$\Sigma \ast u(x): =  < u ; \Sigma(x-\cdot) >$ is well defined
for all $x$ and it is clear that
$\| \Sigma \ast u \|_{W^{k,\infty}(\R^d)} 
\leq \|\Sigma \|_{W^{2k,\infty}(\R^d)} \|u \|_{C_0^k(\R^d)'}$.
At least when $u\in L^1(\R^d)$, it is also clear that 
$\lim_{|x|\rightarrow \infty} \nabla^\ell \Sigma \ast u(x) =0$
for all $\ell \in \{1,..,k\}$.
Since both sides of \eqref{presque} are continuous in $u \in C^k_0(\R^d)'$, the inequality also holds in that space by density.
\end{ProofOf}

\section{Long time behaviour}\label{sec3}

In order to get a precise description of the long time behaviour of $f$, 
we slightly ristrict our hypothesis. We will now ask
\[ 
\begin{array}{ll}
{\bf(H4)}  &
\quad \nabla_x \Phi_0 \in L^\infty(\R_+ \times \R^d),
\quad P \in L^\infty(\R),
\quad |\nabla V| \leq \alpha + \beta V, 
\vspace{0.2cm}
\\ \ds
{\bf(H5)} &  \quad \hat{P}(\omega) \neq 0 \mbox{ for all } \omega \neq 0, 
\quad supp(\hat{\Sigma}) = \R^d. 
\end{array}\]
We can now set our main result:

\begin{theo}{Long time behaviour}
\label{main}

Assume {\bf(H1)-(H5)} and take $f$, a solution of \eqref{kin}-\eqref{CI}
\begin{enumerate}
\item
If $\ds \lim_{|x| \rightarrow \infty} V(x) = +\infty$ or more generally if $(f_t)_{t\geq 0}$ is tight then
\[\begin{array}{llll}
(i) & \ds 
\lim_{T\rightarrow \infty} |\rho|_{C^{1/2}_{W}(T,+\infty)} = 0, 
& (iii) & \ds 
\lim_{t\rightarrow \infty} 
\|\partial_t \Sigma\ast \rho \|_{W^{1,\infty}(\R^d)} =0,
\vspace{0.2cm}
\\ (ii) & \ds
\lim_{t\rightarrow \infty} W (f(t),SEq) = 0, 
& (iv) & \ds 
\lim_{u \rightarrow \infty} \sup_{|t-u|\leq T} 
W(f(t), \mathcal{S}^{\rho(u)}_{t-u} f(u)) = 0.
\end{array}\]
\item
Without any additional assumption $(ii)$-$(iv)$ still hold if one replace $W$
by any metric $\mathrm{d}$ adapted to the weak star topology on $C_0$.
Instead of $(i)$, we always have
$\ds \lim_{u \rightarrow \infty} \sup_{|t-u|\leq T} 
\mathrm{d}(\rho(u),\rho(t)) = 0$.
\end{enumerate}
\end{theo}

\begin{rmk}
\label{W-W1}
When $f_0$ is a probability measure, the distance $W$ can be replaced by  the 
Kantorowich-Rubinstein distance $W_1$ in Theorem \ref{main}-$1)$ thanks to \eqref{KR}. 
\end{rmk}

\begin{rmk}
Under {\bf (A1)}, when $p$ is defined by \eqref{reduct-param},
we always have $p(0)=0$, we point out that 
{\bf (H2)-(H5)} offers more possibility. For example, one can fix $\delta \in (0,1)$, $\varsigma >0$ and take 
$p(t) = \frac{\sgn(t)}{|t|^\delta} 
e^{-\frac{t^2}{2 \varsigma^2} }$. It is clear that
$P: t\mapsto \int_{-\infty}^t p(s) \ud s $ belongs to
$(L^1\cap L^\infty)(\R)$ and one can compute
\[ \hat{P}(\omega) = -\frac{\mathcal{Z}}{\omega} \int_{-\infty}^{+\infty}
e^{-\varsigma^2 (\omega-s)^2/2} \frac{\sgn(s)}{|s|^{1-\delta}} \ud s
= -2 \mathcal{Z} \int_0^\infty \frac{\sinh (\varsigma^2 \omega s)}{\omega} 
e^{-\varsigma^2 (s^2+\omega^2)/2}
\frac{\ud s }{s^{1-\delta}}\]
where we have set 
$\mathcal{Z} := \frac{2 \varsigma}{ \sqrt{2\pi}} 
\int_0^\infty \frac{\sin(t)}{t^\delta} \ud t $.
$\hat{P}(\omega)$ is negative for all $\omega \in \R$. With such parameters $p$, \eqref{kin}-\eqref{CI} is formally closer to the Vlasov equation \eqref{VWRho}.
\end{rmk}

\begin{rmk}
\label{non-conf-simp}
We point out that an additional assumption is required to ensure the tightness of $(f(t))_{t\geq 0}$ even when $V=0$. If one takes for example $d=3$, $\Phi_0 =  V = 0$, $\Sigma$ radially symmetric such that $\supp(\Sigma) \subset B(0,R)$ and 
$p\in L^1(\R_+)$, then setting 
$\mathcal{K} = 2 \|\Sigma\|_{W^{1,\infty}(\R^d)} R^2 \|p\|_{L^1(\R_+)}$,
we fix $r_0, r_1 > 0$ such that $r_0 > 2R$ while
$\mathcal{U}_0 := \frac{r_1^2}{2} - \frac{2\mathcal{K}}{r_0} >0$ and we set
$\ds f_0 = \frac{1}{4\pi(r_0^2 + r_1^2)} 
\delta_{\{ (r_0 u ,r_1 u) \mid u \in \mathbb{S}^2\}}$.
By symetry and mass conservation, the unique solution of 
\eqref{kin}-\eqref{CI} is given at any time by
\begin{equation}
\label{contrex}
f_t = \frac{1}{4\pi (r(t)^2+ \dot{r}(t)^2)} 
\delta_{\{ (r(t) u ,\dot{r}(t) u) \mid u \in \mathbb{S}^2\}}
\end{equation}
while due to \eqref{rep-car}, $r$ solves
\begin{equation}
\label{non-conf-edo}
\left\lbrace 
\begin{array}{l} \ds
\ddot{r}(t) = \int_0^t p(t-s) 
\frac{1}{4\pi r(s)^2} \int_{|y|=r(s)} e_1.\nabla_x \Sigma(r(t) e_1 - y ) \ud S(y) \ud s 
\vspace{0.1cm}
\\ \ds
r(0) = r_0,\quad  \dot{r}(0) = r_1.
\end{array}
\right.
\end{equation}
We now prove that $\lim_{t\rightarrow \infty} r(t) = \infty$.
Since $\Sigma$ is compactly supported, we have
\[ |\ddot{r}(t)| \leq \frac{\|\Sigma\|_{W^{1,\infty}(\R^d)}}{4\pi}
\int_0^t
\frac{|p(t-s)|}{r(s)^2}
\Big(\int_{\begin{subarray}{l} |y|=r(s) \\ |y-r(t)e_1|\leq R \end{subarray}}  \ud S(y) \Big)\ud s. \]
From geometrical considerations, we control the second integral:
\[\int_{\begin{subarray}{l} |y|=r(s) \\ |y-r(t)e_1|\leq R \end{subarray}}  \ud S(y) \leq 4\pi R^2 \mathds{1}_{\{|r(t)-r(s)| \leq R\}}.\]
As long as $r(t) >2 R$ and $\dot{r}(t) >0$, we deduce first
\[ |\ddot{r}(t)| \leq 
\|\Sigma\|_{W^{1,\infty}(\R^d)} R^2 \|p\|_{L^1(\R_+)} \frac{1}{(r(t)-R)^2}\leq \frac{\mathcal{K}}{r(t)^2}.\]
and multiplying by $\dot{r}$, we get simply
\[ \left( \ddot{r} \dot{r} + \mathcal{K}  \frac{\dot{r}}{r^2} \right)(t) \geq 0 \]
By integration we have 
\[ \frac{1}{2} \dot{r}(t)^2 - \frac{\mathcal{K}}{r(t)} \geq 
\frac{1}{2} r_1^2- \frac{\mathcal{K}}{r_0} =:\mathcal{U}_0 \]
We deduce that $\dot{r}$ never vanishes and then that 
$r(t) \geq r_0 + t (2\mathcal{U}_0)^{1/2} $.
Remembering of \eqref{contrex}, we have found $f_0$ such that the unique solution of \eqref{kin}-\eqref{CI} with initial data $f_0$ goes to $0$ for the weak star topology on $C_0(\R^d \times \R^d)$.
\end{rmk}
The proof of Theorem \ref{main} is divided in three part. First, we will prove that Theorem \ref{existence}-$3ii)$ allows us to deduce Theorem \ref{main}-$1iii)$. Second, we will deduce that the omega limit set of $(f(t))_{t\geq 0}$ is contained in 
$SEq$. The other claims of Theorem \ref{main} follow easily in the third part.

\subsection{Relaxation of $\partial_t \Sigma \ast \rho$ to $0$}

The first point leading to Theorem \ref{main} is the following 
result

\begin{lemma}
\label{relax}
Assume {\bf (H1)-(H5)}, take $f$ a solution of 
\eqref{kin}-\eqref{CI} and $\rho$ its spatial density, then 
$ \ds \lim_{t \rightarrow \infty} 
\| \partial_t (\Sigma \ast \rho)(t) \|_{W^{1,\infty}(\R^d)} = 0$.
\end{lemma}

\begin{ProofOf}{Lemma \ref{relax}}

We reminds the notation $[h]_T$ defined by \eqref{hT}.

{\it Step 1: Proving 
$\lim_{t \rightarrow \infty} 
\| P \ast [\partial_t( \Sigma \ast \rho )]_\infty(t)
\|_{L^{\infty}(\R^d)}
= 0 $.}

We set 
\[ u(t) = P \ast [\partial_t( \Sigma \ast \rho )]_\infty(t)
= \int_{-\infty}^t P(s) \partial_t (\Sigma \ast \rho)(t-s) \ud s. \]
According to Theorem \ref{existence}-$3ii)$, we already know that 
$t \mapsto \|u(t)\|_{L^\infty(\R^d)}$ belongs to $L^2(\R)$.
We now prove the time uniform continuity of $u$, derivating its expression, we first get
\begin{equation}
\label{dtu}
\partial_t u(t) = P(t) \partial_t (\Sigma \ast \rho)(0)
+ \int_{-\infty}^t P(s) \partial_t^2 (\Sigma \ast \rho)(t-s) \ud s, 
\end{equation}
From \eqref{cons-masse}, the expression of $\partial_t^2 \rho$ is given by 
\begin{equation}
\label{dt2rho}
\partial_t^2 \rho 
=\nabla_x^2 . \left( \int_{\R^d_v} v \otimes v \ud f \right)
+\nabla_x. (\rho\nabla_x (V+\Phi(t)))
\end{equation}  
According to \eqref{devPhi} and \eqref{dtsigmarho}, for all $t\geq 0$, we have directly
\begin{equation}
\label{nablaphi}
\|\nabla_x \Phi \|_{L^\infty(\R_+ \times \R^d)}
\begin{array}[t]{l} \ds
\leq 
\|\nabla_x \Phi_0 \|_{L^\infty(\R_+ \times \R^d)}
+2 \mathfrak{m} \|P \|_{L^\infty(\R)} \|\nabla_x \Sigma \|_{L^\infty(\R^d)}
\\ \ds \qquad \quad
+  \|P \|_{L^1(\R)} 
 \|\Sigma \|_{W^{2,\infty}(\R^d)} 
(2 \mathfrak{m} \mathcal{E}_1)^{1/2}.
\end{array}
\end{equation}

From {\bf (H4)} and Theorem \ref{existence}-$3i)$, it leads to the uniform bound
\begin{equation}
\label{dt2sigmarho}
\| \partial_t^2 (\Sigma \ast \rho)(t) \|_{L^\infty(\R^d)}
\begin{array}[t]{l} \ds
\leq d^2 \| \nabla^2 \Sigma \|_{L^\infty(\R^d)} 
\int_{\R^d\times \R^d} |v|^2 \ud f_t 
\\ \ds \qquad
+ \| \nabla \Sigma \|_{L^\infty(\R^d)} \int_{\R^d} 
\left(\alpha + \beta V 
+ \|\nabla_x \Phi \|_{L^\infty(\R_+ \times \R^d)}\right) \ud \rho_t
\vspace{0.2cm} \\ \ds
\leq (2 d^2 + \beta) \| \Sigma \|_{W^{2,\infty}(\R^d)} \mathcal{E}_1
\\ \ds \qquad
+ \mathfrak{m} \| \Sigma \|_{W^{2,\infty}(\R^d)} 
\left(\alpha + \|\nabla_x \Phi \|_{L^\infty(\R_+ \times \R^d)}\right).
\end{array}
\end{equation}
Finally coming back to \eqref{dtu} and using \eqref{dtsigmarho} we have 
\[ \| \partial_t u  \|_{L^\infty(\R_+ \times \R^d)}
\begin{array}[t]{l} \ds
\leq \| P \|_{L^\infty(\R)} \|\Sigma \|_{W^{2,\infty}(\R^d)} (2 \mathfrak{m} \mathcal{E}_1)^{1/2}
\\ \ds \qquad
+ \|P \|_{L^1(\R^d)} 
\| \partial_t^2 (\Sigma \ast \rho) \|_{L^\infty(\R_+\times \R^d)}.
\end{array}
\]
Then the function $t \mapsto \| u(t) \|_{L^\infty(\R^d)}$ is uniformly continuous on $\R$. Since it also belongs to $L^2(\R)$, from its expression we deduce
\begin{equation}
\label{cavapourP}
\lim_{t\rightarrow \infty}  \| 
P \ast [\partial_t (\Sigma \ast \rho)]_\infty(t) \|_{L^\infty(\R^d)} = 0. 
\end{equation}

{\it Step 2: Deducing
$\lim_{t \rightarrow \infty} 
\| g \ast [\partial_t( \Sigma \ast \rho )]_\infty(t)
\|_{L^{\infty}(\R^d)}
= 0 $ when $\int_{\R} g(t) \ud t =0$. }

We now set 
\[ Y = \left\lbrace g \in L^1(\R) \left|  \lim_{t\rightarrow \infty}  \| 
g \ast [\partial_t (\Sigma \ast \rho)]_\infty(t) \|_{L^\infty(\R^d)} = 0
\right\rbrace \right. ,
\] 
We have proved that $Y$ contains $P$ and if one note $\tau_h$ the translation operator defined by $(\tau_h g)(x) = g(x+h)$, since 
$(\tau_h P) \ast g = \tau_h(P \ast g)$, $Y$ also contains all the translations of $P$ and their linear combinations, setting
\[ X:= \overline{\mbox{span} \big( \{ \tau_h P \mid h \in \R \} \big)}^{L^1(\R)} \]
$Y$ is closed in $L^1(\R)$ thanks to the Young inequality, therefore it also contains $X$. 
Since 
$\widehat{\tau_h P}(\omega) = e^{ih\omega} \hat{P}(\omega)$, setting $Z(X):= \cap_{g \in X} \hat{g}^{-1}(\{0\})$, 
$Z(X) =  \hat{P}^{-1}(\{0\})$.
$X$ is a closed translation-invariant subspace of 
$L^1(\R)$, it is enough to deduce that
$X= \{ g\in L^1(\R) \mid \hat{g}(\omega) = 0 
\quad \forall \omega \in Z(X) \}$ if 
$\partial Z(X)$ contains no perfect set
(see Corollary $7.2.4$-$(b)$ in \cite{Rud}). 
Since $ \hat{P}^{-1}(\{0\}) \subset \{0\}$ by {\bf (H5)},
we deduce $ \{ g \in L^1(\R) \mid \int_{\R} g(t) \ud t = 0 \}
\subset X \subset Y$.
Coming back to the definition of $Y$, for all 
$g \in L^1(\R)$ such that $\int_{\R} g(t) \ud t = 0$, we have
proved
\begin{equation}
\label{wien}
\lim_{t\rightarrow \infty} 
\| g \ast [\partial_t (\Sigma \ast \rho)]_\infty (t) \|_{L^\infty(\R^d)} = 0.
\end{equation}

{\it Step 3: Deducing $\lim_{t \rightarrow \infty} 
\| \partial_t( \Sigma \ast \rho )(t)
\|_{L^{\infty}(\R^d)} = 0 $.}

We now take $\theta \in \mathscr D( \R)$ such that $\theta \geq 0$, 
$\int_{\R} \theta(t) \ud t =1$ and for all $r,\epsilon >0$, we set
$g_\epsilon(t) = \frac{1}{\epsilon} \theta(\frac{t}{\epsilon}) - \frac{1}{\epsilon} \theta(\frac{t+r}{\epsilon})$. 
On the one hand $g_\epsilon \in L^1(\R)$ and $\int_{\R} g_\epsilon(t) \ud t = 0$, therefore thanks to \eqref{wien},
$\lim_{t \rightarrow \infty}
\| g_\epsilon \ast [\partial_t (\Sigma \ast \rho)]_\infty(t) \|_{L^\infty(\R^d)} = 0$.
On the other hand, for all $r,t,\epsilon >0$, 
\[ \begin{array}{l} \ds
\| g_\epsilon \ast [\partial_t (\Sigma \ast \rho)]_\infty (t)
-\partial_t (\Sigma\ast \rho)(t) + \partial_t (\Sigma\ast \rho)(t+r)
\|_{L^\infty(\R^d)} 
\\ \ds \qquad 
\leq \int_{-\infty}^{t/\epsilon} \theta(s) 
\|\partial_t (\Sigma\ast \rho)(t-\epsilon s) - \partial_t(\Sigma\ast \rho)(t)\|_{L^\infty(\R^d)} \ud s
\\ \ds
\qquad \qquad   
+ \| \partial_t(\Sigma\ast \rho)(t)\|_{L^\infty(\R^d)} 
\int^{+\infty}_{t/\epsilon} \theta(s) \ud s
\\ \ds
\qquad \qquad   
+ \int_{-\infty}^{(t+r)/\epsilon} \theta(s) 
 \|\partial_t (\Sigma\ast \rho)(t+ r -\epsilon s) - \partial_t(\Sigma\ast \rho)(t+r)\|_{L^\infty(\R^d)} \ud s
\\ \ds
\qquad \qquad   
+ \| \partial_t(\Sigma\ast \rho)(t+r)\|_{L^\infty(\R^d)} 
\int^{+\infty}_{(t+r)/\epsilon} \theta(s) \ud s
\\ \ds 
\qquad
\leq 2 \epsilon 
\|\partial_t^2 (\Sigma \ast \rho) \|_{L^\infty(\R_+ \times \R^d)}
\int_{\R} \theta(s) |s| \ud s
\\ \ds
\qquad \qquad  
+ 2 \|\partial_t (\Sigma \ast \rho) \|_{L^\infty(\R_+ \times \R^d)}
\int^{+\infty}_{t/\epsilon} \theta(s) \ud s
\\ \ds 
\qquad \quad
\xrightarrow[ \epsilon \rightarrow 0]{} 0.
\end{array}
\]
Thanks to \eqref{wien}, for all $r\geq 0$, it allows us to deduce first
\[ \lim_{t \rightarrow \infty } 
\| \partial_t (\Sigma\ast \rho)(t+r) - \partial_t (\Sigma\ast \rho)(t)
\|_{L^\infty(\R^d)} =0. \]
By integration, we get
$ \lim_{t \rightarrow +\infty} 
\| \Sigma\ast \rho(t) + s \partial_t(\Sigma\ast \rho)(t)- \Sigma\ast \rho(t+s) \|_{L^\infty(\R^d)} = 0 $ for all
$s> 0$.
Thanks to the mass conservation, we deduce
\[ \limsup_{t \rightarrow \infty} 
\| \partial_t \Sigma \ast \rho(t) \|_{L^\infty(\R^d)} \leq 
2 \frac{\mathfrak{m} \|\Sigma\|_{L^\infty(\R^d)}}{s} .\]
Letting $s$ go to infinity, we have proved
$\lim_{t \rightarrow \infty} 
\| \partial_t \Sigma \ast \rho(t) \|_{L^\infty(\R^d)} =0$.

{\it Step 4: Enlarging the convergence}

It just remains to be proved that the convergence also holds in 
$W^{1,\infty}(\R^d)$. For all $t\geq 0$ and all $x,y \in \R^d$, we have
\[ \begin{array}{l} \ds
|\nabla (\partial_t \Sigma \ast \rho)(t,x) 
- \nabla (\partial_t \Sigma \ast \rho)(t,y) |
\\ \ds
\qquad \qquad
= \left| \int_{\R^d\times\R^d} 
(\nabla^2 \Sigma(x-z)-\nabla^2 \Sigma(y-z))v \ud f_t(z,v) \right|
\\ \ds
\qquad \qquad
\leq (2 \mathfrak{m} \mathcal{E}_1)^{1/2} \|\nabla^2\Sigma - \tau_{x-y}(\nabla^2\Sigma) \|_{L^\infty(\R^d)}
\end{array}
\]

Since $\Sigma \in C^2_0(\R^d)$, 
it allows us to find a continuous modulus $k$ such that $k(0) = 0$ and
$ |\nabla (\partial_t \Sigma \ast \rho)(t,x) 
- \nabla (\partial_t \Sigma \ast \rho)(t,y) |
 \leq k(|x-y|)$.
For all $t,x,y$, we have
\[ \begin{array}{l} \ds
|\partial_t \Sigma \ast \rho(t,x) 
- \partial_t \Sigma \ast \rho(t,y)
- \nabla (\partial_t \Sigma \ast \rho)(t,x).(y-x) |
\\ \ds \qquad \qquad
= \left| \int_0^1 \big(\nabla (\partial_t \Sigma \ast \rho)(t,x+s(y-x) ) 
- \nabla(\partial_t \Sigma \ast \rho)(t,x) \big).(y-x) \ud s \right|
\\ \ds \qquad \qquad
\leq |x-y| \int_0^1 k(s|x-y|)\ud s.
\end{array}\]
Fixing $\epsilon>0$, we first take $\eta$ such that 
$\int_0^1 k(s\eta) \ud s \leq \epsilon$, and then we choose
$y = x + \eta \frac{\nabla(\partial_t \Sigma \ast \rho)(t,x)}{|\nabla(\partial_t \Sigma \ast \rho)(t,x)|}$ we get
\[ |\nabla(\partial_t \Sigma \ast \rho)(t,x)| \leq \epsilon 
+\frac{2}{\eta} \|\partial_t \Sigma \ast \rho(t)\|_{L^\infty(\R^d)}. \]
Since the right hand side can be taken as small as desired when $t$ goes to infinity, we have proved
\[ \lim_{t \rightarrow \infty} 
\| \partial_t( \Sigma \ast \rho )(t) \|_{W^{1,\infty}(\R^d)} = 0 .
\]

\end{ProofOf}

\begin{rmk}
\label{Pitt}
The steps 2-3 establish a variant of Pitt's extension of the Wiener's tauberian theorem (see \cite[Chapter~9]{Rud91}). If  $\hat{P}$ does not cancel on $\R$, this variant is not necessary and
one can deduce directly $\lim_{t \rightarrow \infty} 
\| \partial_t \Sigma \ast \rho(t) \|_{L^\infty} =0$ from 
\eqref{cavapourP} and \eqref{dt2sigmarho} by using
\cite[Theorem~9.7-$(b)$]{Rud91}.

\end{rmk}

\subsection{Characterization of the omega-limit set}

\begin{lemma}{Omega limit set}

\label{ens-omega-limites}
Assume {\bf(H1)-(H5)} and take $f$ a solution of 
\eqref{kin}-\eqref{CI}.
For any sequence $(t_k)_{k\geq 0}$ such that 
$\lim_{k\rightarrow \infty} t_k = +\infty$, 
if $(f(t_k))_{k\geq 0}$ converges to $f^\ast$
for the weak star topology on $C_0(\R^d\times\R^d)$, 
then $f^\ast \in SEq$ and 
$f(t_k+\cdot)$ goes to $\mathcal{S}^{\rho^\ast} f^\ast$ 
in $C_{w\ast}([-T,T],\mathcal{M}_+^1(\R^d\times\R^d))$
for all $T\geq 0$.
\end{lemma}

\begin{ProofOf}{Lemma \ref{ens-omega-limites}}

Take $t_k  \underset{k \rightarrow \infty}{\longrightarrow}\infty$ and 
$f^\ast$ such that 
$f(t_k) \underset{n \rightarrow \infty}{\rightharpoonup} f^\ast$. Fixing 
$T>0$, we define the sequence $(f^k)_{k\geq 0}$ in 
$C_W([-T,T],\mathcal{M}_+^1(\R^d\times\R^d))$ by setting 
$f^k(t) := f(t+t_k)$.

{\it First step: compactness of $(f^k)_k$ }

Pick $\chi$ in $C_c^\infty(\R^d\times\R^d)$ and $t\in [-T,T]$, on
the one hand we have the uniform bound
\begin{equation}
\label{masse}
\left| \int_{\R^d\times\R^d} \chi(x,v) \ud f^k(t) (x,v) \right|
\begin{array}[t]{l} \ds  \leq 
\| f(t+t_k) \|_{\mathcal{M}_+^1(\R^d\times\R^d)} 
\|\chi\|_{L^\infty(\R^d\times\R^d)} \\ \ds
\leq \mathfrak{m}
\|\chi\|_{L^\infty(\R^d\times\R^d)}.
\end{array}
\end{equation}
On the other hand, thanks to \eqref{nablaphi}, $\nabla_x \Phi$ is bounded on 
$\R_+\times\R^d$, therefore
\[
\left| \frac{\ud}{\ud t}\int_{\R^d\times\R^d} \chi(x,v) 
\ud f^k_t(x,v) \right| \begin{array}[t]{l} \ds
= \left| \int_{\R^d\times\R^d}
\left( v.\nabla_x \chi - \nabla_x(V+\Phi(t+t_k)). \nabla_v \chi \right) \ud f^k_t \right|
\\ \ds 
\leq 
\big( \| v.\nabla_x \chi - \nabla_x V.\nabla_v \chi \|_{L^\infty(\R^d\times\R^d)}
\\ \ds 
\qquad \qquad + \|\nabla_x \Phi \|_{L^\infty(\R_+\times \R^d)} 
\|\nabla_v \chi \|_{L^\infty(\R^d\times \R^d)} \big) \mathfrak{m}.
\end{array}
\]
Finally the set 
\[\left\lbrace \left.
t \mapsto \int_{\R^d\times\R^d} \chi(x,v) \ud f^k_t \right| k\in \N
\right\rbrace\]
is equibounded and equicontinuous. Thanks to the Arzela-Ascoli's theorem, we deduce it is compact in $C([-T,T])$. Going back to \eqref{masse}, a simple approximation argument allows us to extend the conclusion to any function $\chi$ in $C_0(\R^d\times\R^d)$. Since this space is separable, using a diagonal argument, we can extract a subsequence and find a measure valued function 
$\tilde{f}$ in $C_{w\ast}([-T,T],\mathcal{M}_+^1(\R^d\times\R^d))$
such that for all $\chi$ in $C_0(\R^d\times\R^d)$,
\begin{equation}
\label{limfk}
\lim_{k\rightarrow\infty} 
\int_{\R^d\times\R^d} \chi(x,v) \ud f^k_t
= \int_{\R^d\times\R^d} \chi(x,v) \ud \tilde{f}_t 
\end{equation}
holds uniformly on $[-T,T]$. 

{\it Second step: limit equation satisfied by $\tilde{f}$.}

Setting $\Phi^k(t) = \Phi(t+t_k)$ for any $k\geq 0$, $f^k$ solves
\begin{equation}
\label{edpfk}
\left\lbrace \begin{array}{l}
\partial_t f^k + v.\nabla_x f^k = \nabla_x(V+\Phi^k(t)). \nabla_v f^k
\\
f^k(0) = f(t_k).
\end{array} \right.
\end{equation}
By definition of $(t_k)_k$ and $f^\ast$, 
we already get $\tilde{f}(0) = f^\ast$ from \eqref{limfk} and it is
also clear that 
\begin{equation}
\label{faible}
\partial_t f^k + v.\nabla_x f^k - \nabla_x V. \nabla_v f^k
\underset{k \rightarrow \infty}{\rightharpoonup}
\partial_t \tilde{f} + v.\nabla_x \tilde{f} - \nabla_x V. \nabla_v \tilde{f} 
\end{equation}
in $\mathscr D'((-T,T)\times \R^d \times \R^d)$. In order to pass to the limit in the last term of \eqref{edpfk}, we now prove that
$\nabla_x \Phi^k$ converges strongly to 
$-\kappa \nabla \Sigma \ast \rho^\ast$. First, according to \eqref{devPhi},
we get
\[\begin{array}{l} \ds
\int_{-T}^T\| \nabla_x( \Phi^k(t)+\kappa\Sigma \ast \rho^k(t))
\|_{L^\infty(\R^d)} \ud t
\\ \qquad \ds
\leq  \int_{-T}^T
 \big(\| \nabla_x \Phi_0(t+t_k) \|_{L^\infty(\R^d)} + |P(t+t_k)| 
 \|\nabla_x \Sigma \ast \rho(0) \|_{L^\infty(\R^d)} \big) \ud t
\\ \qquad \qquad \ds
+ \int_{-T}^T \int_0^{t+t_k} |P(t+t_k-s)| 
\| \nabla_x(\partial_t \Sigma \ast \rho)(s) \|_{L^\infty(\R^d)} \ud s \ud t
\\ \qquad \ds
\leq \int_{t_k-T}^\infty 
\big(\| \nabla_x \Phi_0(t) \|_{L^\infty(\R^d)} 
+ \mathfrak{m} \| \Sigma \|_{W^{1,\infty}(\R^d)} |P(t)| \big) \ud t
\\ \qquad \qquad \ds
+ 2T \sup_{-T\leq t \leq T} 
\int_0^{t +t_k} |P(t+t_k-s)| \| \partial_t (\Sigma \ast \rho)(s) 
\|_{W^{1,\infty}(\R^d)} \ud s
\end{array}
\]
According to {\bf (H3)}, the first of those last two terms goes to $0$ when $k$ goes to infinity.
For the second one we split $[0,t+t_k]$ in $[0,\frac{t+t_k}{2}]$ and 
$[\frac{t+t_k}{2}, t+t_k]$, setting
$u(s) = \| \partial_t (\Sigma \ast \rho)(s) \|_{W^{1,\infty}(\R^d)}$
we get 
\[\int_0^{t +t_k} |P(t+t_k-s)| u(s) \ud s 
\begin{array}[t]{l} \ds
\leq \|u\|_{L^\infty} \int_{\frac{t +t_k}{2}}^{t +t_k} |P(t)|\ud t 
+ \|P\|_{L^1(\R_+)}
\| u \|_{L^\infty(\frac{t +t_k}{2},t +t_k)}
\vspace{0.1cm}
\\ \ds
\leq \|u\|_{L^\infty} \int_{\frac{t_k-T}{2}}^\infty |P(t)|\ud t 
+ \|P\|_{L^1(\R_+)}
\| u \|_{L^\infty(\frac{t_k-T}{2},\infty)}.
\end{array}
\]
Since  $\lim_{s\rightarrow \infty}u(s) =0$ thanks to
Lemma \ref{relax}, we have already proved
\begin{equation}
\label{convphi1}
\lim_{k \rightarrow \infty} 
\int_{-T}^T\| \nabla_x( \Phi^k(t)+\kappa\Sigma \ast \rho^k(t))
\|_{L^\infty(\R^d)} \ud t = 0 .
\end{equation}
It remains to establish the convergence of 
$\nabla_x \Sigma \ast \rho^k$ when $k$ goes to infinity. 
According to  Theorem \ref{existence}-$3i)$, 
$\int_{\R^d\times\R^d} |v|^2 \ud f_t^k (x,v) \leq 2\mathcal{E}_1$ while 
$f^k_t$ converges weakly to to $\tilde{f}_t$ on $[-T,T]$ by \eqref{limfk}. 
A standard approximation argument allows us to deduce that
$\rho^k(t) \underset{n \rightarrow \infty}{\rightharpoonup} 
\tilde{\rho}(t)$ for any $t$ in $[-T,T]$.
On the one hand, since $\Sigma \in C^2_0(\R^d)$, 
$\nabla_x \Sigma (x-\cdot) \in C_0(\R^d,\R^d)$ for all $x \in \R^d$, it is enough to get punctually
$\lim_{k\rightarrow \infty} (\nabla_x \Sigma\ast \rho^k)(t,x)
=(\nabla_x \Sigma \ast \tilde{\rho})(t,x) $ 
for all $(t,x)\in [-T,T]\times\R^d$.
On the other hand, the following estimates
$\|\nabla_x (\nabla_x \Sigma \ast \rho^k(t) ) \|_{L^\infty(\R^d)}
\leq \mathfrak{m} \|\Sigma \|_{W^{2,\infty}(\R^d)}$
and 
$\|\partial_t (\nabla_x \Sigma \ast \rho^k(t) )\|_{L^\infty(\R^d)} \leq (2 \mathfrak{m} \mathcal{E}_1)^{1/2} \|\Sigma \|_{W^{2,\infty}(\R^d)}$
holds for all $k\geq 0$. It allows us to deduce that the convergence is uniform on all compact set of $[-T,T]\times \R^d$.
For all $R>0$, we have already established
\begin{equation}
\label{convrhok1}
\lim_{k \rightarrow \infty} \| \nabla_x \Sigma\ast \rho^k
-\nabla_x \Sigma \ast \tilde{\rho} 
\|_{L^\infty([-T,T]\times B(0,R))} = 0.
\end{equation}
By using the following decomposition
\[ \nabla_x \Sigma \ast \rho^k(t) - \nabla_x \Sigma \ast \rho^\ast
\begin{array}[t]{l} \ds
= \big(\nabla_x \Sigma \ast \rho^k(0) 
- \nabla_x \Sigma \ast \tilde{\rho}(0) \big) 
\\ \qquad \qquad \qquad \ds 
+ \int_{t_k}^{t_k+t} 
\nabla_x (\partial_t \Sigma\ast \rho)(s) \ud s,
\end{array}\]
Lemma \ref{relax} and \eqref{convrhok1} allows us to deduce 
\begin{equation}
\label{convrhok2}
\lim_{k \rightarrow \infty} 
\|\nabla_x \Sigma \ast \rho^k- \nabla_x \Sigma \ast \rho^\ast 
\|_{L^\infty([-T,T] \times B(0,R))} = 0.
\end{equation}
Coming back to \eqref{convphi1}, we get 

\begin{equation}
\label{convphi3}
\lim_{k \rightarrow \infty} 
\int_{-T}^T\| \nabla_x \Phi^k(t)+\kappa \nabla_x \Sigma \ast \rho^\ast
\|_{L^\infty(B(0,R))} \ud t = 0 .
\end{equation}

To conclude that step, we now take $\chi $
in $\mathscr D((-T,T)\times \R^d \times \R^d)$ and $R$ such that 
$ supp(\chi) \subset [-T,T] \times B(0,R)  \times B(0,R)$. For any 
$k\geq 0$, we have
\[ \begin{array}{l} \ds
\left|\int_{(-T,T) \times \R^d\times\R^d} 
\nabla_v \chi.\nabla_x \Phi^k \ud f^k_t \ud t 
+ \kappa \int_{(-T,T) \times \R^d\times\R^d} 
\nabla_v \chi . \nabla_x (\Sigma \ast \rho^\ast) \ud \tilde{f}_t \ud t 
\right|
\\ \ds \qquad
\leq \int_{-T}^T \left( \int_{\R^d\times\R^d} 
\| \nabla_v \chi. \nabla_x (\Phi^k + \kappa (\Sigma \ast \rho^\ast)) 
\|_{L^\infty(B(0,R))} \ud f_t^k \right)\ud t 
\\ \ds \qquad \qquad
+ \kappa \int_{-T}^T \left| \int_{\R^d\times\R^d} 
\nabla_v \chi . \nabla_x (\Sigma \ast \rho^\ast) 
(\ud \tilde{f}_t -\ud f^k_t) \right| \ud t
\end{array}\]
The first term is controled by 
$ \mathfrak{m} \int_{-T}^T\| \nabla_x \Phi^k(t)+\kappa \nabla_x \Sigma \ast \rho^\ast \|_{L^\infty(B(0,R))} \ud t$ which goes to $0$ by
\eqref{convphi3}. The second term goes to $0$ as well by 
\eqref{limfk}. Finally, we have established  $\nabla_x \Phi^k. \nabla_v f^k 
\rightharpoonup -\kappa (\nabla_x \Sigma \ast \rho^\ast). \nabla_v \tilde{f}$ in $\mathscr D'((-T,T)\times \R^d \times \R^d)$.

Coming back to \eqref{edpfk} and \eqref{faible},
$\tilde{f}$ solves the linear transport equation 
\begin{equation}
\label{edpftilde}
\left\lbrace \begin{array}{l}
\partial_t \tilde{f} + v.\nabla_x \tilde{f} 
= \nabla_x(V-\kappa \Sigma \ast \rho^\ast). \nabla_v \tilde{f}
\\
\tilde{f}(0) = f^\ast.
\end{array} \right.
\end{equation}

{\it Third Step: Conclusion.}

According to the definition of 
$\mathcal{S}$, 
$\tilde{f}(t) = \mathcal{S}^{\rho^\ast}_t f^\ast$ for all $t$ in $[-T,T]$ by \eqref{edpftilde}.
The sequence $\{f^k\}_{k\geq 0}$ is compact and has a unique acumulation point, thus it converges. Coming back to the definition of $f^k$, we have proved 
 \[f(t_k+t) \underset{k \rightarrow \infty}{\rightharpoonup} 
\mathcal{S}^{\rho^\ast}_t f^\ast\] 
uniformy in $t \in [-T,T]$ for all $T>0$.
Thanks to \eqref{convrhok1} and \eqref{convrhok2}, 
$\nabla \Sigma \ast \tilde{\rho}(t) = \nabla \Sigma \ast \rho^\ast$. Since $supp(\hat{\Sigma}) = \R^d$ the mapping 
$\rho \mapsto \nabla \Sigma \ast \rho $ is injective in $\mathcal{M}_+^1(\R^d)$, 
therefore $\tilde{\rho}(t) = \rho^\ast$.
Taking $T$ as big as desired we have establised
$\int_{\R^d_v} \ud \mathcal{S}^{\rho^\ast}_t f^\ast = \rho^\ast$
for all $t\in \R$.

\end{ProofOf}

\subsection{Proof of Theorem \ref{main}}

The proof of Theorem \ref{main} is short. We have already established all the main arguments it requires.

\begin{ProofOf}{Theorem \ref{main}}

Assuming {\bf(H1)-(H5)},
we point out that $(i)$ is already proved by Lemma \ref{relax}. 

{\it Proof of 2)}

We take $\mathrm{d}$ a distance adapted to the weak star topology on $C_0$.
It is clear that $\{f(t)\}_{t\geq 0}$ is compact for this topology. For any sequence $(f(t_k))_k$ such that
$\lim_{k\rightarrow \infty}t_k = +\infty$, we can extract a subsequence (still noted $t_k$) and find an accumulation point 
$f^\ast$ such that 
$\lim_{k\rightarrow \infty} \mathrm{d}(f(t_k),f^\ast)=0 $.
Thanks to Lemma \ref{ens-omega-limites}, $f^\ast \in SEq$ and 
\begin{equation}
\label{int-th-1}
\lim_{k\rightarrow \infty} 
\mathrm{d}(f(t_k + t ), \mathcal{S}^{\rho^\ast}_t f^\ast)=0
\end{equation}
uniformly in $t$ on all compact set of $\R$. We now set 
$g_k (t) = \mathcal{S}^{\rho(t_k)}_t f(t_k)$. By definition of 
$\mathcal{S}$, $g_k$ is the unique solution of 
\[ \left\lbrace \begin{array}{l}
\partial_t g_k + v.\nabla_x g_k 
= \nabla_x (V-\kappa \Sigma \ast \rho(t_k)).\nabla_v g_k 
\\
g_k(0) = f(t_k).
\end{array}\right.\]
For all $T>0$, from the same method we used to prove Lemma 
\ref{ens-omega-limites}, $\{g_k\}_k$ is compact 
in $C_{w\ast}([-T,T],\mathcal{M}_+^1(\R^d\times\R^d))$ and it has a unique accumulation point:  
$t \mapsto \mathcal{S}^{\rho^\ast}_t f^\ast$, thus it converges.
according to \eqref{int-th-1}, we already have
\begin{equation}
\label{int-th-2}
\lim_{k\rightarrow \infty} \sup_{-T \leq t \leq T}
\mathrm{d}(f(t_k+t), \mathcal{S}^{\rho(t_k)}_t f(t_k) )=0,
\qquad
\lim_{k\rightarrow \infty} \mathrm{d}(f(t_k),SEq)=0.
\end{equation}
Since $f^\ast \in SEq$, 
$\int_{\R^d_v } \ud (\mathcal{S}^{\rho^\ast}_t f^\ast) \equiv \rho^\ast$ for all $t$, then the bound of kinetic energy 
$\int_{\R^d\times\R^d} |v|^2\ud f_t(x,v) \leq 2 \mathcal{E}_1$ given by Theorem \ref{existence}-$2i)$ and \eqref{int-th-1} allows us to deduce 
\begin{equation}
\label{int-th-3}
\lim_{k\rightarrow \infty} \sup_{-T \leq t \leq T} 
\mathrm{d}(\rho(t_k+t) ,\rho^\ast) 
= 
\lim_{k\rightarrow \infty} \sup_{-T \leq t \leq T} 
\mathrm{d}(\rho(t_k+t),\rho(t_k))
= 0.
\end{equation}
For any sequence $(t_k)_k$ such that $\lim_{k \rightarrow \infty} t_k = + \infty$,
we have proved that \eqref{int-th-2} and \eqref{int-th-3} are at least satisfied for a subsequence of $(t_k)_k$. It allows us to deduce 
\begin{equation}
\label{int-th-4}
\lim_{u\rightarrow \infty} \sup_{-T \leq t \leq T}
\mathrm{d}(f(u+t), \mathcal{S}^{\rho(u)}_t f(u) )=0,
\quad \lim_{t\rightarrow \infty} \mathrm{d}(f(t),SEq)=0
\end{equation}
as well as
\[\lim_{u\rightarrow \infty} \sup_{-T \leq t \leq T} 
\mathrm{d}(\rho(u+t),\rho(u)) =0.\]

{\it Proof of 1)}

If $\ds \lim_{x\rightarrow \infty} V(x) = + \infty$ then 
$(f(u))_{u\geq 0}$ is tight thanks to Theorem \ref{existence}-$3i)$.
Thanks to Lemma \ref{char-est},
$\ds  \left( \mathcal{S}^{\rho(u)}_t f(u)
\right)_{\begin{subarray}{l} u\geq 0 \\ -T\leq t \leq T  \end{subarray} }$
is also tight by \eqref{rep-car}. It allows us to improve directly \eqref{int-th-4} as
\[\lim_{u\rightarrow \infty} \sup_{-T \leq t \leq T}
W(f(u+t), \mathcal{S}^{\rho(u)}_t f(u) )=0,
\quad \lim_{t\rightarrow \infty} W(f(t),SEq)=0.\]
We just have to prove the Hölder continuity estimate on $\rho$. Take $g \in SEq$ such that for a sequence $t_k$ going to infinity,
$\lim_{k\rightarrow \infty} W(f(t_k),g)=0$. A standard approximation argument allows us to state
\[ \int_{\R^d\times \R^d} \big(V(x) + \frac{|v|^2}{2}\big) \ud g(x,v)
\leq
\liminf_{k\rightarrow \infty}
\int_{\R^d\times \R^d} \big(V(x) + \frac{|v|^2}{2}\big) \ud f_{t_k}(x,v)
\leq 
\mathcal{E}_1.\]
while $\int_{\R^d\times\R^d} \ud g (x,v) = \mathfrak{m}$.
By compactness of $\{f(t)\}_{t\geq 0}$, we deduce
\begin{equation}
\label{la-derr-th}
\lim_{t\rightarrow \infty} W(f(t),SEq(\mathfrak{m},\mathcal{E}_1))=0.
\end{equation}
According to Proposition \ref{stat} for all $s,t\geq T$,
\[ \frac{W(\rho_t,\rho_s)}{|t-s|^{1/2}} 
\begin{array}[t]{l} \ds
\leq  \mathscr H (\mathfrak{m},\mathcal{E}_1)
\left| \frac{1}{|t-s|}
 \int_s^t (W ( f_\tau, SEq(\mathfrak{m}, \mathcal{E}_1) ) )^{1/2}  \ud \tau \right|^{1/2}
\\ \ds
\leq \mathscr H (\mathfrak{m},\mathcal{E}_1) 
\sup_{t\geq T}  \big(W( f_t, SEq(\mathfrak{m}, \mathcal{E}_1) ) \big)^{1/4}
\end{array}\]
Then
$\ds \lim_{t\rightarrow \infty} |\rho|_{C^{1/2}_{W}(t,+\infty)} = 0$ by \eqref{la-derr-th}.
\end{ProofOf}

\section{Application: the mean-field equations}\label{sec4}

\subsection{General property}

We now consider the $N$-particle system $(q_k,p_k)_{1\leq k \leq N}$ defined by \eqref{chm}-\eqref{chminit}.
We define the emprical density $\hat{f}^N$ and its associated spatial empirical density $\hat{\rho}^N$ by 
\[ \hat{f}^N_t =  \frac{1}{N} \sum_{k=1}^N \delta_{(q_k(t),p_k(t))},
\quad 
\hat{\rho}^N_t =  \int_{\R^d_v} \ud \hat{f}^N_t(v)  
=\frac{1}{N} \sum_{k=1}^N \delta_{q_k(t)}. \]
We point out that the evolution equation of $\Psi$ 
in \eqref{chmbdb} can be recast as
\[ (\partial_{tt}^2 \Psi- c^2 \Delta_z \Psi)(t,x,y) = 
-\sigma_2(y)\sigma_1 \ast \hat{\rho}^N_t(x). \]
Under {\bf (A1)-(A2)}, Lemma \ref{resolution} allows us to deduce that \eqref{chm} is an equivalent formulation of \eqref{chmbdb}
completed by the initial data $(\Psi_0,\Psi_1)$.
The system \eqref{chm}-\eqref{chminit} is linked with our previous considerations thanks to the following fundamental observation:
\begin{lemma}{Assume {\bf (H2)-(H3)},}
\label{equiv-kyn-chm}
\begin{itemize}
\item[(i)] 
If $(q_k,p_k)_{1 \leq k \leq N}$ solves \eqref{chm}, then 
$\hat{f}^N$ solves \eqref{kin}-\eqref{rho}.
\item[(ii)] 
Reciproquely if $f$ solves \eqref{kin}-\eqref{CI} 
with the initial data 
$f_0 = \frac{1}{N} \sum_{k=1}^N \delta_{(q_{k,0},p_{k,0})}$, then we can find $(q_k,p_k)_{1 \leq k \leq N}$ solving 
\eqref{chm}-\eqref{chminit} such that for all $t\geq 0$,
$f(t) = \frac{1}{N} \sum_{k=1}^N \delta_{(q_{k}(t),p_{k}(t))}$.
\end{itemize}
\end{lemma}

The proof of Lemma \ref{equiv-kyn-chm} is classical. We refer the reader to \cite{GV} for the details.  $(i)$ is obtained by computing the time derivative of 
$t\mapsto < \hat{f}^N_t \chi >$ for any test function 
$\chi \in C_c^{\infty}(\R^d\times\R^d)$. $(ii)$ is a direct concequence of \eqref{rep-car}. It allows us to deduce first

\begin{proposition}
\label{exist-chm}
Assume {\bf (H2)-(H3)}, for all family of initial datas \\
$(q_{k,0},p_{k,0})_{1 \leq k \leq N}$ there exists a unique global solution 
$(q_k,p_k)_{1\leq k \leq N}$ of \eqref{chm}-\eqref{chminit}.
\end{proposition}

\begin{ProofOf}{Proposition \ref{exist-chm}}

The existence is a direct concequence of Theorem \ref{existence} and Lemma
\ref{equiv-kyn-chm}: for all initial data $(q_{k,0},p_{k,0})_{1 \leq k \leq N}$, we define the empirical initial distribution 
$ \hat{f}^N_0 =  \frac{1}{N} \sum_{k=1}^N \delta_{(q_{k,0},p_{k,0})}$, 
take the solution given by Theorem \ref{existence}-$1)$ and apply Lemma 
\ref{equiv-kyn-chm}-$ii)$. For the uniqueness we have to be a little more cautious. 

If $(q^1_k,p^1_k)_{1 \leq k \leq N}$ 
and $(q^2_k,p^2_k)_{1 \leq k \leq N}$
both solve \eqref{chm}-\eqref{chminit}, then according to Lemma \ref{equiv-kyn-chm}-$i)$ and Theorem \ref{existence}-$2)$ (the integrability condition is trivially satisfied on $\R_+$), we have the identity
\[ \frac{1}{N} \sum_{k=1}^N \delta_{(q_k^1(t),p^1_k(t))}
=  \frac{1}{N} \sum_{k=1}^N \delta_{(q_k^2(t),p_k^2(t))} \]
(but it only gives directly 
$(q^2_k,p^2_k)_{1 \leq k \leq N} 
= (q^1_{\varsigma(k)},p^1_{\varsigma(k)})_{1 \leq k \leq N}$ for some 
$\varsigma$ in $\mathfrak{S}_N$). For all $k$, it allows us to deduce that 
 $(q^1_k,p^1_k)$ and $(q^2_k,p^2_k)$ both solve the same ordinary differential equation
\[ \left\lbrace \begin{array}{l}
\ds \dot{q} =p  \\ \ds \dot{p}= -\nabla_x(V+\Phi(t))(q)
\end{array}\right. \]
with the same initial data. Uniqueness follows directly from 
the Cauchy-Lipschitz theorem.
\end{ProofOf}

\subsection{Long time behaviour}

Since the empirical density associated with the solutions of 
\eqref{chm} satisfies \eqref{kin}-\eqref{rho}, we can already apply Theorem \ref{main} but in this case due to our restriction to a very particular kind of initial data, the consequences are stronger. Under the following additional assumption
\[ {\bf (M)} \quad \lim_{x\rightarrow \infty} V(x) = + \infty, 
\quad x.\nabla V(x) \geq 0, \quad x.\nabla \Sigma(x) \leq 0, \]
("M" stands for monotony) we can even prove the convergence of all the particles.

\begin{theo}
\label{monotonie}
Assume {\bf (H1)-(H5)} and {\bf(M)}, if 
$(q_k,p_k)_{1\leq k \leq N}$ solves \eqref{chm} then  
$\ds \lim_{t\rightarrow \infty} p_k(t) = 0$ 
for all $k\in \{1,..,N\}$. Moreover,
\begin{enumerate}
\item 
if $x.\nabla V(x) >0$ for all $x\neq 0$, or if $\Sigma$ is radially symmetric and $\nabla V$ do not cancel on $\R^d\setminus \{0\}$, then
\[ \lim_{t\rightarrow \infty} (q_k,p_k) = (0,0) 
\quad \mbox{ for all } k\in \{1,..,N\}.\]
\item
if $\Sigma$ is radially symmetric and if $\nabla\Sigma$ does not cancel on 
$B(0,r)\setminus \{0\}$ for some $r>0$, then 
we can find $m \leq N$ 
and a family of disjoint sets $ (I_k)_{1 \leq k \leq m}$ such that $\ds \cup_{k=1}^m I_k = \{1,..,N\}$ and for all $k\neq \ell$ in $\{1,..,m\}$
and all $i,j \in \{1,..,N\}$,
\[\begin{array}{llll} \ds 
(i) & \ds 
\lim_{t \rightarrow \infty} \max_{i,j \in I_k} |q_i(t) - q_j(t) | = 0,
& (iii) & \ds
\liminf_{t \rightarrow \infty} 
\min_{(i,j) \in I_k\times I_\ell} |q_i(t) - q_j(t) | \geq r,
\\ \ds (ii) & \ds
\lim_{t\rightarrow \infty} 
\nabla \Sigma(q_i(t)-q_j(t)) = 0,
& (iv) & \ds
\lim_{t\rightarrow \infty} 
D\big(q_i(t); \nabla V^{-1}(\{0\})\big) = 0.
\end{array}\]
Moreover, if $\nabla V(x) \neq 0$ for all $x\in \R^d\setminus B(0,R)$ then
$m \leq \big( \frac{2R+ r}{r} \big)^d$.
\end{enumerate}
\end{theo}

\begin{rmk}{Comments}

\begin{enumerate}
\item
We point out that {\bf (M)} allows us to take 
$V(x)= [|x|-R]_+^2$. If one can prove that all the particules stay in a compact set of $\R^d$ during the whole evolution when $V=0$, then the conclusions of Theorem \ref{monotonie}-$2)$ would also hold in this situation.
\item
In the next sections, we will see that the speed of convergence of those quantity goes to $0$ when $N$ goes to infinity.
\end{enumerate}
\end{rmk}
The proof of Theorem \ref{monotonie} can be splitted in two intermediate results. First, defining $Eq^N$, the set of equilibrium of \eqref{chm}
\begin{equation}
\label{EqN} 
Eq^N:= \big\{ (x_k,v_k)_{1\leq k \leq N} \big|
\nabla_x V(x_k) -\kappa \sum_{j=1}^N \nabla_x \Sigma (x_k -x_j)
= v_k = 0 \mbox{ for all } k \big\}, 
\end{equation}
Theorem \ref{main}-$1)$ can be simply traduced in our new framwork by
\begin{proposition}
\label{kin-part}
Assume {\bf (H1)-(H5)} and take $(q_k,p_k)_{1\leq k \leq N}$ a solution 
of \eqref{chm}. 
If $\ds \lim_{t\rightarrow \infty} V(x) = + \infty$ 
or more generaly if
$(q_k,p_k)_{1\leq k \leq N} \in C_b(\R_+, (\R^d\times\R^d)^N)$ 
then
\[ \lim_{t\rightarrow \infty} 
D\big((q_k(t),p_k(t))_{1\leq k \leq N}, Eq^N \big) = 0. \]
\end{proposition}
The second step is to prove that under {\bf(M)} and the other restrictions we added in Theorem \ref{monotonie}, the set $Eq^N$ is exactly
\[ \big\{ (x_k,v_k)_{1\leq k \leq N} \big|
\nabla_x V(x_k) = \nabla_x \Sigma (x_k -x_\ell)
= v_k = 0 \mbox{ for all } k,\ell \big\}. \]
It is a direct concequence of the following lemma:

\begin{lemma}
\label{vir}
Assume {\bf (M)} and suppose $\Sigma $ is even, then for any $(x_k,v_k)_{1\leq k \leq N} \in Eq^N$,
\[x_k . \nabla V(x_k) = 0 
\mbox{ and }
(x_k-x_j). \nabla \Sigma (x_k-x_j) =0
\mbox{ for all } k,j \in \{1,..,N\}.\]
\end{lemma}
We now prove those claims.

\begin{ProofOf}{Proposition \ref{kin-part}}

Take $(q_k,p_k)_{1\leq k \leq N}$ a  solution of \eqref{chm},
we start by setting the concequences of our previous results. 
From Lemma \ref{equiv-kyn-chm}, the empirical density 
$\hat{f}^N_t =  \frac{1}{N} \sum_{k=1}^N \delta_{(q_k(t),p_k(t))}$
solves \eqref{kin}-\eqref{rho}. According to Theorem \ref{existence}, we can find $\mathcal{E}_1$ such that 
\[ \frac{1}{N} \sum_{k=1}^N \frac{|p_k|^2}{2} + V(q_k) \leq \mathcal{E}_1. \]
Therefore if $\lim_{x\rightarrow \infty} V(x) = + \infty$ or if
$(q_k,p_k)_{1\leq k \leq N} \in C_b(\R_+, (\R^d\times\R^d)^N)$,
 we can find a compact set $K$ such that for any $t\geq 0$ and any 
$k \in \{1,..,N\}$, $q_k(t) \in K$. From Theorem \ref{main}-$1)$, we get
\begin{equation}
\label{CVemp0}
\lim_{t \rightarrow \infty} W_1(\hat{f}^N_t, SEq) =0.
\end{equation}
We now set 
\[ \begin{array}{lcll} 
\mathscr J: & \big( \ds K^N , |\cdot | \big) 
& \longrightarrow & \ds 
\big( \mathcal{M}^1_+(\R^d\times \R^d) , W_1 \big) 
\\ \ds 
\mathscr J: & \ds (x_k,v_k)_{1 \leq k \leq N} & \longmapsto &
\ds \frac{1}{N} \sum_{k=1}^N \delta_{(x_k,v_k)}.
\end{array} \] 
and we explain our strategy: in a first step, we prove that
\eqref{CVemp0} allows us to deduce 
$\lim_{t \rightarrow \infty} 
D \big((q_k(t),p_k(t))_{1,\leq k \leq N}, \mathscr J^{-1}(SEq) \big)$, 
we conclude in a second step by proving 
$ \mathscr J^{-1}(SEq)  = Eq^N$.

{\it Step 1: topological considerations}

It is clear that  $\mathscr J$ is continuous, therefore any  accumulation points of $(\hat{f}^N(t))_{t\geq 0}$ belongs to $\mathscr J(K^N)$ by compactness. Thanks to \eqref{CVemp0}, we already deduce 
\begin{equation}
\label{CVemp1} 
\lim_{t \rightarrow\infty} 
W_1 \left( \mathscr J ((q_k(t),p_k(t))_{1 \leq k \leq N} ), 
\mathscr J (K^N) \cap SEq \right) = 0. 
\end{equation}
Taking $X =(x_k,v_k)_{1 \leq k \leq N}$ in $K^N$ and $\varsigma$ in 
$\mathfrak{S}_N$, we define the natural action 
of $\mathfrak{S}_N$ on $K^N$ by setting 
$X^\varsigma := (x_{\varsigma(k)},v_{\varsigma(k)})_{1 \leq k \leq N}$, we
note $\bar{X} \in K^N / \mathfrak{S}_N$ the element 
$\{X^\varsigma\mid \varsigma \in \mathfrak{S}_N\}$. 
It is clear that $\mathscr J(X) = \mathscr J(Y)$ 
if and only if we can find  $\varsigma$ such that $Y = X^\varsigma$ 
therefore $\mathscr{J}$ can be restricted in a continuous bijection
\[ \begin{array}{lcll} 
\bar{\mathscr J}: & \big( \ds K^N / \mathfrak{S}_N, \bar{D} \big) 
& \longrightarrow & \ds 
\big( \mathscr J(K^N) , W_1 \big) 
\\ \ds 
\bar{\mathscr J}: & \ds \overline{(x_k,v_k)}_{1 \leq k \leq N} & \longmapsto &
\ds \frac{1}{N} \sum_{k=1}^N \delta_{(x_k,v_k)}
\end{array} \]
(for the quotient distance 
$\bar{D}(\bar{X},\bar{Y}) := \min_{X \in \bar{X},Y\in \bar{Y}} D(X,Y)$).
Since $K^N$ is compact, $K^N / \mathfrak{S}_N$ is also compact. Hence
$ \bar{\mathscr J}^{-1}$ is continuous and we get from \eqref{CVemp1}
\begin{equation}
\label{CVemp1-2} \lim_{t \rightarrow \infty} 
\bar{D} \Big( \overline{(q_k(t),p_k(t))}_{1,\leq k \leq N}, 
\bar{\mathscr J}^{-1} 
\left( \mathscr J (K^N) \cap SEq \right) \Big)= 0.
\end{equation}
Obviously, for any subset $A$ of $\mathscr J (K^N)$, 
$\cup_{\mu \in A} \bar{\mathscr J}^{-1}(\mu) = \mathscr J^{-1}(A)$ 
as well as 
$ \big( \mathscr J^{-1}(A) \big)^{\varsigma} = \mathscr J^{-1}(A)$ 
for all $\varsigma$ in $\mathfrak{S}_N$. It allows us to transform
\eqref{CVemp1-2} in
\begin{equation}
\label{CVemp2}
\lim_{t \rightarrow \infty} 
D \left((q_k(t),p_k(t))_{1,\leq k \leq N}, 
\mathscr J^{-1}(SEq) \right)= 0
\end{equation}
(for the usual distance $D(x,y) = |x-y|$ on $(\R^d\times\R^d)^N$).

{\it Step 2: conclusion}

In order to conclude, we now prove that
 $\mathscr J^{-1}(SEq)$ is the set of equilibrium of \eqref{chm}.
Take $(x_{k,0},v_{k,0})_{1\leq k \leq N}$ in 
$\mathscr J^{-1} (SEq)$. Coming back to definition \ref{stat}, 
and solving \eqref{VWRho} for $\rho = \frac{1}{N} \sum_{k=1}^N \delta_{x_{k,0}}$ with the initial data 
$ g_0 = \frac{1}{N} \sum_{k=1}^N \delta_{(x_{k,0},v_{k,0})}$,
it means that the unique solution $(x_k,v_k)_{1\leq k \leq N}$ of 
 \begin{equation}
\label{chm-stat}
\left\lbrace \begin{array}{l} \ds 
\dot{x_k} = v_k, \\  \ds \dot{v_k} = -\nabla_x V (x_k)
+ \frac{\kappa}{N} \sum_{j=1}^N \nabla_x \Sigma(x_k-x_{0,j})
\end{array} \right.
\qquad 
\left\lbrace \begin{array}{l} \ds 
x_k(0) = x_{k,0} \vspace{0.5cm} \\ v_k(0) = v_{k,0}
\end{array} \right.
\end{equation}
satisfies for any $t\geq 0$,
\[
\frac{1}{N} \sum_{k=1}^N \delta_{x_k(t)} 
= \frac{1}{N} \sum_{k=1}^N \delta_{x_{k,0}}.
\]
For any $t\geq 0$, it allows us to deduce that $ \ds (x_k(t))_{1,\leq k \leq N} \in 
\{ (x_{0,\varsigma (k)})_{1\leq k \leq N}\mid \varsigma \in \mathfrak{S}_N \}$ and since $\{ (x_k(t),v_k(t))_{1\leq k \leq N} \mid t \geq 0 \}$
is connex $ x_k(t) = x_{k,0} $ for all $k,t$.
Thanks to \eqref{chm-stat}, we deduce first $v_k(t)=0$ for all $k,t$ and then
$\frac{\kappa}{N} \sum_{j=1}^N \nabla_x \Sigma(x_k(t) -x_{0,j}) 
-\nabla_x V(x_k(t)) = 0$. Finally,  $(x_k,v_k)_{1\leq k \leq N} 
\in \mathscr J^{-1} ( SEq )$
if and only if 
\[
\left\lbrace \begin{array}{l}
\ds v_k = 0 \vspace{0.2cm}
\\ \ds
\frac{\kappa}{N} \sum_{j=1}^N \nabla_x \Sigma(x_k -x_j) 
-\nabla_x V(x_k) = 0
\end{array} \right.
\mbox{ for all } k\in \{1,..,N\}.
\]
Therefore,  
$\mathscr J^{-1} (SEq) = Eq^N$ and the result follows by \eqref{CVemp2}.

\end{ProofOf}

\begin{ProofOf}{Theorem \ref{monotonie}}

{\it Radial case}

We first simplify the expression of $Eq^N$ by assuming that 
$\Sigma$ is radially symmetric.
Take $(x_k,v_k)_{1\leq k \leq N} \in Eq^N$,
from Lemma \ref{vir}, we get 
$\nabla \Sigma (x_k-x_j)=0$ for all $k,j \in \{1,..,N\}$.
Since for all $k \in \{1,..,N\}$ we also have 
$\nabla_x V(x_k) =\frac{\kappa}{N} \sum_{j=1}^N \nabla_x \Sigma(x_k -x_j) $ by \eqref{EqN}, we get $\nabla V(x_k)=0$ for all $k$ as well. We have proved
\begin{equation}
\label{EqN-rad-M}
Eq^N = \big\{ (x_k,v_k)_{1\leq k \leq N} \big|
\nabla_x V(x_k) = \nabla_x \Sigma (x_k -x_\ell)
= v_k = 0 \mbox{ for all } k,\ell \big\}.
\end{equation}

{\it Proof of $(1)$: }
 
If $x.\nabla V(x) >0$ for all $x \neq 0$, then 
according to Lemma \ref{vir}, $Eq^N = \{(0,0)^{\otimes N}\}$.
If $\Sigma$ is radially symmetric and if $\nabla V$ does not cancel on $\R^d \setminus \{0\}$, then the same conclusion holds by\eqref{EqN-rad-M}. According to Proposition \ref{kin-part}, we have proved
\[ \lim_{t\rightarrow \infty} (q_k(t),p_k(t)) = (0,0) 
\quad \mbox{ for all } k\in \{1,..,N\}.\]

{\it Proof of $(2)(i)$-$(iv)$:}

We now assume that $\Sigma$ is radially symmetric such that
\begin{equation}
\label{r}
\nabla \Sigma(x) \neq 0 \mbox{ for all } x \in B(0,r)\setminus \{0\}
\end{equation}
for some $r>0$. Take $(x_k,v_k)_{1\leq k \leq N} \in Eq^N$,
for all $k,\ell \in \{1,..,N\}$ according to \eqref{EqN-rad-M} and \eqref{r}, 
$|x_k-x_\ell| \in \{0\} \cup [r,+\infty)$.
Then, thanks to Proposition \ref{kin-part}, 
\[ \lim_{t\rightarrow \infty} 
D\big(|q_k(t)-q_\ell(t)|, \{0\} \cup [r,+\infty)\big) = 0.\]
Since $t\mapsto |q_k(t)-q_\ell(t)|$ is continuous, we get
\[ \lim_{t\rightarrow \infty} |q_k(t)-q_\ell(t)| = 0 
\quad \mbox{or} \quad 
\liminf_{t\rightarrow \infty} |q_k(t)-q_\ell(t)| \geq r
\quad \mbox{ for all } k,\ell \in \{1,..,N\}.\]
It allows us to find $m$ and $(I_k)_{1\leq k \leq m}$ such that 
for all $k \neq \ell$ in $\{1,..,m\}$
\[
 \lim_{t \rightarrow \infty} 
\max_{i,j \in I_k} |q_i(t) - q_j(t) | = 0, 
\quad \mbox{and} \quad
\liminf_{t \rightarrow \infty} 
\min_{(i,j) \in I_k\times I_\ell} |q_i(t) - q_j(t) | \geq r.
\]
From \eqref{EqN-rad-M} and Proposition \ref{kin-part}, we have proved $(i)$-$(iv)$. 

{\it Uniform bound on $m$}

We now suppose that $\nabla V^{-1}(\{0\}) \subset B(0,R)$.
For all $k \in \{1,..,m\}$, we pick $i_k \in I_k$. 
Fixing $\epsilon>0$, we can find $t \geq 0$ such that 
for all $k \neq l$ in $\{1,..,m\}$,
\[B\big(q_{i_k}(t),\frac{r-\epsilon}{2} \big) 
\cap B\big(q_{i_l}(t),\frac{r-\epsilon}{2} \big) 
= \emptyset 
\quad \mbox{ and } \quad 
B \big(q_{i_k}(t),\frac{r-\epsilon}{2} \big) \subset B \big(0, R+\frac{r}{2} \big).\]
By considering the measure of those sets,  we get
$ m \omega_d \big( \frac{r-\epsilon}{2} \big)^d 
\leq \omega_d \big( R + \frac{r}{2} \big)^d$.
Taking $\epsilon$ as small as desired, it can be recast as 
$m \leq \big( \frac{2R+r}{r} \big)^d$.

\end{ProofOf}

\begin{ProofOf}{Lemma \ref{vir}}

Assume {\bf(M)} and take
$(x_k,v_k)_{1\leq k \leq N} \in Eq^N$, for all $k \in \{1,..,N\}$ we have
\begin{equation}
\label{EqN-scal}
x_k.\nabla_x V(x_k)
-\frac{\kappa}{N} \sum_{j=1}^N x_k.\nabla_x \Sigma(x_k -x_j) 
= 0.
\end{equation}
In view of summing those identity for all $k \in \{1,..,N\}$,  we make the following computation:
\[ \sum_{k=1}^N  \sum_{j=1}^N x_k. \nabla \Sigma (x_k-x_j)
\begin{array}[t]{l} \ds
= \sum_{k=1}^N  \sum_{j=1}^N (x_k-x_j). \nabla \Sigma (x_k-x_j)
+ \sum_{k=1}^N  \sum_{j=1}^N x_j. \nabla \Sigma (x_k-x_j)
\\ \ds 
= \sum_{k=1}^N  \sum_{j=1}^N (x_k-x_j). \nabla \Sigma (x_k-x_j)
+ \sum_{j=1}^N  \sum_{k=1}^N x_k. \nabla \Sigma (x_j-x_k)
\end{array}\]
Since $\Sigma$ is even, we get
\begin{equation}
\label{viriel}
\sum_{k=1}^N  \sum_{j=1}^N x_k. \nabla \Sigma (x_k-x_j)
= \frac{1}{2}  \sum_{1\leq k,j \leq N} (x_k-x_j). \nabla \Sigma (x_k-x_j).
\end{equation}
Summing \eqref{EqN-scal} for all $k \in \{1,..,N\}$, \eqref{viriel} allows us to deduce
\begin{equation}
\label{viriel2}
\sum_{k=1}^N x_k.\nabla V(x_k) 
- \frac{\kappa}{2N} 
\sum_{1\leq k,j \leq N} (x_k-x_j). \nabla \Sigma (x_k-x_j) =0.
\end{equation}
Thanks to {\bf (M)}, all the terms involved in \eqref{viriel2} are non negative, by assuming  $(x_k,v_k)_{1\leq k \leq N} \in Eq^N$, we have established
\[
x_k . \nabla V(x_k) = 0 
\mbox{ and }
(x_k-x_j). \nabla \Sigma (x_k-x_j) =0
\mbox{ for all } k,j \in \{1,..,N\}.
\]
\end{ProofOf}

\subsection{Mean-field limit}

We now introduce some probability on the intial distribution of 
particle. Our motivation is to exhibit some lower bounds on the convergences proved in the last subsection in a "typical initial configuration", the meaning of this expression is given by the following hypothesis  
\[ {\bf (H6)} \quad \mbox{
$(q_{k,0},p_{k,0})_{k\in \N}$ is identically distributed according to } f_0  \in \mathcal{P}(\R^d\times\R^d). \]
Once again for all $N>0$, we consider $(q_k^N,p_k^N)_{1\leq k \leq N}$, the solution of \eqref{chm}-\eqref{chminit} and the empirical density 
$\hat{f}^N_t = \frac{1}{N} \sum_{k=1}^N \delta_{(q_k^N(t),p_k^N(t))}$
and we warm the reader that those object are now random variable depending on the initial configuration. The system is now described by the measure 
$f^{(N,N)}_t$ acting on $(\R^d \times\R^d)^N$ such that 
\[ f^{(N,N)}_t ( \mathcal{A}) = \mathbb{P}
[ (q_1^N(t),p_1^N(t),...,q_N^N(t),p_N^N(t)) \in \mathcal{A} ]
.\]
Thanks to the symetry in \eqref{chm}, at any time $t$ the particles 
stay identically distributed according to
the probability law  of one particle $f^{(1,N)}_t$, the first marginal of $f^{(N,N)}_t$:
\[ f^{(1,N)}_t(\mathcal{A}) := \int_{A\times(\R^d\times\R^d)^{N-1}}
\ud f^{(N,N)}_t = \mathbb{P}[(q_1^N(t),p^N_1(t)) \in \mathcal{A}]
.\]
The measure $\mathbb{E}[\hat{f}^N_t]$ is identically equal to $f^{(1,N)}_t$ since for all Borel subset $\mathcal{A}$ of $\R^d\times\R^d$,
\[ \mathbb{E}[\hat{f}^N_t(\mathcal{A}) ]
= \frac{1}{N} \sum_{k=1}^N \mathbb{E}[\delta_{(q_k^N(t),p^N_k(t))}(\mathcal{A})]
= \frac{1}{N} \sum_{k=1}^N \mathbb{P}[(q_k^N(t),p^N_k(t))\in \mathcal{A}] 
= f^{(1,N)}_t(\mathcal{A}).\]
When $N$ goes to infinity we have proved in \cite{GV} that $\hat{f}^N$ and eventually $f^{(1,N)}$ converges to a solution of \eqref{kin}-\eqref{CI}. From
Proposition \ref{kin-part}, we know that those measures get concentrated on $\R^d_x\times\{0\}$ in large time. 
However when $V$ is a confining potential, we will see that this localisation of the measure is impossible for any solution of \eqref{kin}-\eqref{CI} with an initial data $f_0\in L^1(\R^d\times\R^d)$. We summarize those informations in the following result:

\begin{theo}{Mean-field limit}
\label{lim-chm}

Assume {\bf (H1)-(H6)},
\begin{enumerate}
\item 
Almost surely, from any subsequence of $(\hat{f}^N)_{N\geq 1}$we can extract a subsequence $(\hat{f}^{N_k})_{k\geq 1}$ and find $f$, a solution of \eqref{kin}-\eqref{CI} such that 
\[ \lim_{k\rightarrow \infty} \sup_{0\leq t \leq T} 
W_1(\hat{f}^{N_k}_t, f_t) = 0 \mbox{ for all } T>0. \]
\item 
If $\nabla^2 V \in L^\infty(\R^d)$, then for all $t\geq 0$ we can find 
$C(t)$ explicit such that
\begin{itemize}
\item[(i)] 
$\ds W_1(f^{(1,N)}_t, f_t)  \leq \frac{C(t)}{\sqrt{N}}.$
\item[(ii)]
$\ds \mathbb{E}\left[ \left| 
\int_{\R^d\times\R^d} \chi \ud \widehat{f}^N_t 
-\int_{\R^d\times\R^d} \chi \ud f_t \right| \right]
\leq \frac{C(t)\|\nabla \chi \|_{L^\infty} + \|\chi \|_{L^\infty} }{\sqrt{N}}$ for any test function  
$\chi \in W^{1,\infty}(\R^d\times\R^d)$.
\end{itemize}
\item
If $\lim_{t\rightarrow \infty} V(x) = +\infty$ and 
$f_0\in L^1(\R^d\times\R^d)$ then there exists $\Theta>0$ 
depending on $f_0$, $V$, $\|P\|_{L^1(\R)}$, 
$\| \Sigma \|_{W^{2,\infty}(\R^d)}$ and
$\| \nabla_x \Phi_0\|_{L^1(\R_+, L^\infty(\R^d))}$ such that for all solution of \eqref{kin}-\eqref{CI} with initial data $f_0$,
\begin{itemize}
\item[(i)] $\ds \liminf_{t\rightarrow\infty} 
W_1(\hat{f}^N(t), f(t)) \geq \Theta$ surely,
\item[(ii)]
$\ds \liminf_{t\rightarrow\infty} 
W_1(f^{(1,N)}(t), f(t)) \geq \Theta$.
\end{itemize}
\end{enumerate}
\end{theo}
We point out that the proof of Theorem \ref{lim-chm} allows us to explicit the constant $\Theta$. However, the values available by this method would be far to be optimal in most of the practical situations.
With the restrictions of Theorem \ref{monotonie}, one can get
 
\begin{proposition}
\label{Theta}
Assume {\bf(M)} and suppose that $\Sigma$ is radially symmetric,
$f_0 \in L^\infty(\R^d\times\R^d)$ and 
$(\nabla V)^{-1}(\{0\}) \subset B(0,R)$ for some $R\geq 0$ then one can take
\begin{itemize}
\item 
$ \Theta = \frac{d}{d+1} 
\big[ (\omega_{d}^2 (d+1) (R+1)^d 
\|f_0\|_{L^\infty(\R^d\times\R^d)})^{-\frac{1}{d}} \wedge 1 \big]$ when $R>0$. 
\item 
$ \Theta = \frac{2d}{2d+1} 
\big[ (\omega_{2d} (2d+1) 
\|f_0\|_{L^\infty(\R^d\times\R^d)})^{-\frac{1}{2d}} \wedge 1 \big]$ when $R=0$. 
\end{itemize}
\end{proposition}

\begin{rmk}{Comments}

\begin{itemize}
\item
The two first points of Theorem \ref{lim-chm} are classical and mostly already established 
(both in a different probabilistic context) in \cite{GV} where they allowed us to justify the physical validity of \eqref{kin}-\eqref{CI} with a continuous density $f_0 \in L^1(\R^d\times\R^d)$. 
\item
When $\lim_{|x|\rightarrow \infty} V(x) = +\infty$, $\{\hat{f}^N(t)\}_{t\geq 0}$ and $\{f(t)\}_{t\geq 0}$ are both compact for the topology defined by the Kantorowich-Rubinstein distance $W_1$. Theorem \ref{lim-chm}-$3)$ allows us to deduce that there is no orbital stability to expect for any distance weak enough to allow a sequence of finite sum of Dirac measures to approximate a function $f \in L^1(\R^d\times\R^d)$.
\end{itemize}
\end{rmk}

\begin{rmk}
\label{C-explicite}
Following the proof of Theorem \ref{lim-chm}-$2)$ steps by steps, one can get
\[\begin{array}{l} \ds
C(T) = 3 \|\nabla \Sigma \|_{L^\infty(\R^d)} 
\int_0^T \|p\|_{L^1(0,t)}
\Big( 1 + \int_t^T 
\exp \big( \int_s^T r_2(\tau) \ud \tau \big) \ud s  \Big) \ud t,
\\ \ds 
r_2(t): = 
 1+\|\nabla^2 V\|_{L^\infty(\R^d)} 
+\|\nabla^2 \Phi_0\|_{L^\infty([0,t]\times \R^d)}
+2 \|\nabla^2 \Sigma\|_{L^\infty(\R^d)} \|p\|_{L^1(0,t)} .
\end{array}\]
When $\lambda := \|r_2 \|_{L^\infty(\R_+)}$ is well defined, 
it leads to the simpler estimate
\[ C(T) \leq 
  3 \|\nabla \Sigma \|_{L^\infty(\R^d)} \|p\|_{L^1(\R_+)}
\big( T+ \frac{e^{\lambda T} -1 -\lambda T}{\lambda^2}  \big).\]
\end{rmk}

\begin{ProofOf}{Theorem \ref{lim-chm}}

{\it Proof of (1)}

Thanks to {\bf(H6)}, for any 
$\varphi \in C_b(\R^d\times\R^d)$
and any $k$ in $\{1,...,N\}$, we already have
$\mathbb{E}[\varphi(q_{k,0},p_{k,0})] 
= \int_{\R^d\times\R^d} \varphi(x,v) \ud f_0(x,v) $. 
The strong law of large numbers and  the separability of 
$C_0(\R^d\times\R^d)$  allows us to deduce almost surely, for all $\varphi$ in $C_0(\R^d\times\R^d)$
\[
\int_{\R^d\times\R^d} \varphi(x,v) \ud \hat{f}^N_0
= \frac{1}{N} \sum_{k=1}^N \varphi(q_{k,0},p_{k,0}) 
\xrightarrow[N \rightarrow \infty]{}
\int_{\R^d\times\R^d} \varphi(x,v) \ud f_0
\]
while for all $N\geq 0$, 
$\int_{\R^d\times\R^d} \ud \hat{f}^N_0
=\int_{\R^d\times\R^d} \ud f_0 =1 $. 
It is enough to ensure that $f^N_0$ converges to $f_0$ for the weak topology on $C_b(\R^d\times\R^d)$.
Since the Kantorowich-Rubinstein distance metrizes this topology, we deduce
\begin{equation}
\label{limchm1}
 \lim_{N \rightarrow \infty} W_1(\hat{f}_0^N,f_0) = 0 
\mbox{ almost surely. }
\end{equation}
For all $(q_{k,0},p_{k,0})_{1\leq k \leq N}$, thanks to Lemma 
\ref{equiv-kyn-chm}, $(\hat{f}^N)_{N\geq 1}$ is a sequence of solution of \eqref{kin}-\eqref{rho}. 
From any subsequence of $(\hat{f}^N)_{N\geq 1}$, the arguments of the first step of the proof of 
Lemma \ref{ens-omega-limites} allows us to extract a subsequence $N_k$
and to find a measure valued function $f$ such that for all  
$\chi$ in $C_0(\R^d\times\R^d)$,
\begin{equation}
\label{limchm2}
\lim_{k\rightarrow \infty} 
\int_{\R^d\times\R^d} \chi (x,v) \ud \hat{f}^{N_k}_t(x,v) = 
\int_{\R^d\times\R^d} \chi (x,v) \ud f_t(x,v)
\end{equation}
holds uniformly on all compact set of $[0,+\infty)$. From now the proof is classical, we just summarize the main arguments by sake of completeness. 
According to \eqref{limchm1}, $(\hat{f}^N_0)_{N\geq 1}$ is tight, from \eqref{rep-car} and Lemma \ref{char-est}, we deduce that  
$(\hat{f}^{N_k}_t)_k$ is time-uniformly tight on all compact set of 
$[0,+\infty)$. It allows us to improve \eqref{limchm2} in
\begin{equation}
\label{limchm3}
\lim_{k \rightarrow \infty} \sup_{0\leq t \leq T} 
W_1(\hat{f}^{N_k}(t),f(t))) = 0 \quad \mbox{ for all }T>0.
\end{equation}
As in the second step of the proof of Lemma 
\ref{ens-omega-limites} all the linear term of \eqref{kin} pass through
the limit when $k$ goes to infinity
\[
\partial_t \hat{f}^{N_k} + v.\nabla_x \hat{f}^{N_k} 
- \nabla_x (V+\Phi_0). \nabla_v \hat{f}^{N_k}
\underset{k \rightarrow \infty}{\rightharpoonup}
\partial_t f + v.\nabla_x f - \nabla_x (V+\Phi_0). \nabla_v f,
\]
while thanks to \eqref{limchm3} and \eqref{KR}, 
\[ \lim_{k \rightarrow \infty} 
\|\nabla_x \Sigma \ast \hat{\rho}^{N_k}- \nabla_x \Sigma \ast \rho^f 
\|_{L^\infty([0,T] \times \R^d)} = 0. \]
For all $\chi \in \mathscr D((0,T) \times \R^d\times\R^d)$, it is enough to prove the weak convergence of the non linear term:
\[ \begin{array}{l} \ds
\lim_{k\rightarrow \infty}  
\int_{(0,T)\times\R^d\times\R^d} 
\nabla_v \chi(t).\big(\int_0^t 
\nabla_x(\Sigma \ast \hat{\rho}^{N_k}(s)) p(t-s)\ud s\big) \ud \hat{f}^N_t(x,v) \ud t 
\\ \ds \qquad
= \int_{(0,T)\times\R^d\times\R^d} 
\nabla_v \chi(t).\big(\int_0^t 
\nabla_x(\Sigma \ast \rho^f(s)) p(t-s)\ud s\big) \ud f_t(x,v) \ud t.
\end{array}\]
Therefore, $f$ solves \eqref{kin}-\eqref{CI}.

{\it Proof of (2)-(i)}

The result has already been proved in
\cite{GV} where an understanding of stochastic PDEs is required
due to the presence of an additionnal diffusion term in \eqref{kin} (which can be taken equal to $0$). We just introduce the objects (simpler here) involved in the previous proof before inviting the reader to skip in \cite{GV}.

We now suppose $\nabla^2V \in L^\infty(\R^d)$ and we point out that the solution of \eqref{kin}-\eqref{CI} is unique thanks to Theorem \ref{existence}-$2)$. In order to simplify the expressions 
we set 
\[ \mathcal{L}(h)(t,x) = \int_0^t p(t-s) \Sigma \ast \rho^h(s) \ud s.\]
$(q_{k,0},p_{k,0})_{k\geq 1}$ still defined by 
{\bf (H6)}, we introduce  a new sequence of random variable 
$(\tilde{q}_k,\tilde{p}_k)_{k\geq 1}$ defined by the following ordinary differential equation:
\begin{equation}
\label{chm-couplage}
\left\lbrace \begin{array}{l} \ds 
\frac{\ud}{\ud t}\tilde{q}_k = \tilde{p}_k, 
\vspace{0.2cm}\\  \ds 
\frac{\ud}{\ud t}\tilde{p}_k = 
-\nabla_x (V+\Phi_0-\mathcal{L}(f)) (\tilde{q}_k)
\end{array} \right. 
\quad
\left\lbrace \begin{array}{l} 
\tilde{q}_k(0) = q_{k,0},
\vspace{0.4cm}\\ 
\tilde{p}_k(0) = p_{k,0}. \end{array} \right.
\end{equation}
By construction $\tilde{q}_k(0)= q_k^N(0)$ and $\tilde{p}_k(0)= p_k^N(0)$.
From {\bf (H6)}, for any $t\geq 0$, 
$(\tilde{q}_k(t),\tilde{p}_k(t))_{k\geq 1}$ is a family of random variable iid. We temporary write $\mu_t$, their common law.
Still according to to {\bf (H6)} and the initial condition 
in \eqref{chm-couplage}, 
$\mu_0 = f_0$.
For any $\varphi$ in $\mathcal{D}(\R^d\times\R^d)$, thanks to 
\eqref{chm-couplage},
\[ \frac{\ud}{\ud t}\varphi(\tilde{q}_1,\tilde{p}_1)
= \tilde{p}_1.\nabla_x \varphi(\tilde{q}_1,\tilde{p}_1) 
- \nabla_x (V+\Phi_0-\mathcal{L}(f))(\tilde{q}_1). 
\nabla_v \varphi(\tilde{q}_1,\tilde{p}_1).\]
We integrate this expression on $[0,T]$ and take the expected value, we get for $\mu$:
\[\begin{array}[t]{l}\ds 
\int_{\R^d\times\R^d} \varphi  \ud \mu_T 
- \int_{\R^d\times\R^d} \varphi  \ud f_0
\\ \ds \qquad
= \int_0^T\int_{\R^d\times\R^d} (v.\nabla_x \varphi(x,v) 
- \nabla_x (V +\Phi_0-\mathcal{L}(f)) (x).\nabla_v \varphi(x,v) ) \ud \mu_t \ud t
\end{array}\]
which is a weak formulation of 
\[ \left\lbrace \begin{array}{l}
\partial_t \mu + v.\nabla_x \mu 
= \nabla_x (V+\Phi_0-\mathcal{L}(f)).\nabla_v \mu
\\ \ds 
\mu(0) = f_0.
\end{array}\right.\]
solved by $f$. Hence by uniqueness $\mu \equiv f$.

We now estimate the difference between 
$(q_k^N,p_k^N)_{1\leq k \leq N}$ and 
$(\tilde{q}_k,\tilde{p}_k)_{1\leq k\leq N}$. From 
\eqref{chm} and \eqref{chm-couplage}, we get for all $k$
\[ 
\left\lbrace \begin{array}{l} \ds
\frac{\ud}{\ud t}(q_k^N - \tilde{q}_k ) = p_k^N - \tilde{p}_k
\vspace{0.2cm}
\\ \ds
\frac{\ud}{\ud t}(p_k^N - \tilde{p}_k ) = 
\big(\nabla_x V(\tilde{q}_k) - \nabla_x V(q_k^N) \big)
+\big(\nabla_x \Phi_0(\tilde{q}_k) - \nabla_x \Phi_0(q_k^N) \big)
\\ \ds \qquad \qquad \qquad
+ \big( \nabla_x \mathcal{L}(\hat{f}^N)(q_k^N)
- \nabla_x \mathcal{L}(\hat{f}^N)(\tilde{q}_k) \big)
+ \nabla_x \mathcal{L}(f-\hat{f}^N) (\tilde{q}_k).
\end{array} \right.\]
Setting
$z_k^N= (q_k^N,p_k^N)$, $\tilde{z}_k = (\tilde{q}_k,\tilde{p}_k)$ 
and
\[ r_1(t) = 1+\|\nabla^2 V\|_{L^\infty(\R^d)} 
+\|\nabla^2 \Phi_0(t)\|_{L^\infty(\R^d)}
+\|\nabla^2 \Sigma\|_{L^\infty(\R^d)} \|p\|_{L^1(0,t)} \]
we get
\[ \frac{\ud}{\ud t} |z_k^N - \tilde{z}| 
\leq r_1(t)  |z_k^N - \tilde{z}| 
+ |\nabla_x \mathcal{L}(f-\hat{f}^N) (\tilde{q}_k)|.\]
As announced we now invite the reader to skip in the proof of 
Theorem 5.3 in \cite{GV} to get directly for all $i$
\begin{equation}
\label{GV}
\mathbb{E}\left[ \sup_{0\leq t \leq T} |z_i^N- \tilde{z}_i|(t)\right]
\leq \frac{C(T)}{\sqrt{N}}.
\end{equation}
Since $f^{(1,N)}_t$ and $f_t$ are the respective laws of $z_1(t)$ and $\tilde{z_1}(t)$, the law $\pi_t$ of $(z_1(t),\tilde{z_1}(t))$
is a coupling of $(f^{(1,N)}_t,f_t)$.
From the definition of the Kantorowich-Rubinstein distance we get
\[ W_1\left(f_t,f^{(1,N)}_t\right) 
\leq \int_{(\R^d\times\R^d)^2} |z-\tilde{z}|\wedge 1 
\ud \pi_t(z,\tilde{z})
= \mathbb{E}\left[|z_i^N(t)- \tilde{z}_i(t)| \wedge 1 \right]
\leq \frac{C(T)}{\sqrt{N}}. \]

{\it Proof of (2)-(ii)}

The second estimate is also a direct concequence of \eqref{GV}:
\[ \begin{array}[t]{l} \ds
\mathbb{E}\left[ \left| 
\int_{\R^d\times\R^d} \chi \ud \widehat{f}^N_t 
-\int_{\R^d\times\R^d} \chi \ud f_t \right| \right]
= \mathbb{E}\left[ \left| \frac{1}{N}\sum_{i=1}^N \chi(z_i(t)) -\int_{\R^d\times\R^d} \chi \ud f_t \right| \right] 
\\ \quad \qquad \qquad \qquad \ds
\leq \mathbb{E}[|\chi(z_1(t)) - \chi(\tilde{z}_1(t))|] +
\mathbb{E}\left[ \left| 
\frac{1}{N} \sum_{i=1}^N \chi(\tilde{z}_i(t))
-\int_{\R^d\times\R^d} \chi \ud f_t \right| \right] 
\\ \quad \qquad \qquad \qquad \ds
\leq \|\nabla \chi \|_{L^\infty} \frac{C(t)}{\sqrt{N}} 
+\mathbb{E}\left[ \left| 
\frac{1}{N} \sum_{i=1}^N \chi(\tilde{z}_i(t))
- \mathbb{E}[\chi(\tilde{z}_1(t))] \right| \right] 
\\ \quad \qquad \qquad \qquad \ds
\leq \|\nabla \chi \|_{L^\infty} \frac{C(t)}{\sqrt{N}} + \frac{(\mbox{Var}[\chi(\tilde{z_1}(t))])^{1/2}}{\sqrt{N}} 
\\ \quad \qquad \qquad \qquad \ds
\leq \frac{C(t)\|\nabla \chi \|_{L^\infty} + \|\chi \|_{L^\infty} }{\sqrt{N}},
\end{array}\]
where we have just applied  \eqref{GV} to estimate the first term and  the law of large numbers to the family $(\chi(\tilde{z}_i(t))_{1\leq i \leq N}$ to deal with the second term.

{\it Proof of (3):}

We recall the preservation of the measure by the flow $\varphi$ defined by \eqref{char}. For any $R,\eta >0$, we estimate
\[ \int_{\R^d_x \times B(0,\eta)}  f(t,x,v) \ud x \ud v
\begin{array}[t]{l} \ds
\leq\int_{B(0,R) \times B(0,\eta)} f_0 (\varphi_t^0(x,v)) \ud x \ud v
\\ \quad \ds 
+ \big(\sup_{|x|\geq R} \frac{1}{V(x)}\big) 
\int_{(\complement B(0,R)) \times B(0,\eta)} f(t,x,v)V(x) \ud x \ud v
\\ \ds
\leq \sup_{|A| \leq (\omega_d)^2 R^d\eta^d} \int_A f_0 \ud x \ud v
+ \mathcal{E}_1 \sup_{|x|\geq R} \frac{1}{V(x)} 
\end{array}\]
Fixing $\epsilon >0$, $R$ such that 
$\mathcal{E}_1  \sup_{|x|\geq R} \frac{1}{V(x)} \leq \epsilon/2$ 
and then $\eta$ small enough to ensure
$ \sup_{|A| \leq (\omega_d)^2 R^d\eta^d} \ud f_0(A) \leq \epsilon/2$,
we get at any time
\[
\int_{\R^d_x \times B(0,\eta)}  f(t,x,v) \ud x \ud v \leq \epsilon
\]
Thanks to \eqref{KR}, we get first
\[ W_1(f_t,\hat{f}^N_t) 
\begin{array}[t]{l} \ds
\geq  \int_{\R^d\times\R^d} (1\wedge |v|)f(t,x,v) \ud x \ud v 
- \int_{\R^d\times\R^d} (1\wedge |v|)\ud \hat{f}_t^N (x,v) 
\\ \ds
\geq  (1\wedge \eta) \int_{|v|\geq \eta} f(t,x,v) \ud x \ud v 
-\frac{1}{N}\sum_{k=1}^N (1 \wedge |p_k^N(t)|) 
\\ \ds
\geq (1\wedge \eta)(1-\epsilon) 
-\frac{1}{N} \sum_{k=1}^N (1 \wedge |p_k^N(t)|).
\end{array}\]
The second term goes to $0$ by Proposition \ref{kin-part}, taking 
$\Theta = (1\wedge \eta)(1-\epsilon)$, we deduce $(i)$.
By the same way for $f^{(1,N)}$, we have
 \[ W_1(f_t,f^{(1,N)}_t) 
\begin{array}[t]{l} \ds
\geq  \int_{\R^d\times\R^d} (1\wedge |v|)f(t,x,v) \ud x \ud v 
- \int_{\R^d\times\R^d} (1\wedge |v|)\ud f^{(1,N)}_t (x,v) 
\\ \ds
\geq (1\wedge \eta)(1-\epsilon) -\mathbb{E}[1 \wedge |p_1^N(t)|].
\end{array}\]
and Proposition \ref{kin-part} allows us to get $(ii)$ as well
thanks to the dominated convergence theorem.
\end{ProofOf}
We now exhibit some better estimates on $\Theta$.

\begin{ProofOf}{Proposition \ref{Theta}}

Thanks to \eqref{rep-car},
$\|f(t)\|_{L^\infty(\R^d\times\R^d)} 
= \|f_0\|_{L^\infty(\R^d\times\R^d)}$ for all $t\geq 0$. 
We set 
\begin{equation}
\label{chi-R} \chi_R(x,v) := d\big((x,v),B_{\R^d_x}(0,R)\times\{0_{\R^d_v}\}\big)\wedge 1= ( [|x|-R]_+^2 +|v|^2 )^{1/2} \wedge 1. 
\end{equation}
Combining the definitions of $\hat{f}^N$, $f^{(1,N)}$ and \eqref{KR}, we have 
\begin{equation}
\label{Kant-R}
\begin{array}{l} \ds
W_1(f_t,\hat{f}^N_t) 
\geq \int_{\R^d\times\R^d} \chi_R(x,v) f_t(x,v) \ud x \ud v
- \frac{1}{N} \sum_{k=1}^N \chi_R(q_k^N(t),p_k^N(t))
\\ \ds
W_1(f_t,f^{(1,N)}) 
\geq \int_{\R^d\times\R^d} \chi_R(x,v) f_t(x,v) \ud x \ud v
- \mathbb{E} \big[ \chi_R(q_1^N(t),p_1^N(t))\big]
\end{array}
\end{equation}
While thanks to Theorem \ref{monotonie}, 
\begin{equation}
\label{test-0}
\lim_{t \rightarrow \infty} \frac{1}{N} \sum_{k=1}^N \chi_R(q_k^N(t),p_k^N(t))=0, \quad 
\lim_{t \rightarrow \infty} 
\mathbb{E} \big[ \chi_R(q_1^N(t),p_1^N(t))\big] = 0.
\end{equation}
In order to conclude we just have to find a time uniform lower bound on 

$\int_{\R^d\times\R^d} \chi_R(x,v) f_t(x,v) \ud x \ud v$.

{\it Estimate when $R=0$:}

For all $\eta \in [0,1]$ and all $t \geq 0$, we have
\begin{equation}
\label{R=0}
 \int_{\R^d\times\R^d} \chi_0(x,v) f_t(x,v)  \ud x \ud v
\begin{array}[t]{l} \ds
= \int_{\R^d\times\R^d} (|(x,v)|\wedge 1) f_t(x,v)  \ud x \ud v
\\ \ds
\geq \eta \int_{|(x,v)|\geq \eta} f_t(x,v)  \ud x \ud v
\\ \ds
\geq \eta \big( 1 - \int_{|(x,v)|\leq \eta} f_t(x,v)  \ud x \ud v ) 
\\ \ds
\geq \eta \big( 1 - \omega_{2d} \eta^{2d} \|f_0\|_{L^\infty(\R^d\times \R^d)} \big) = : g(\eta).
\end{array}
\end{equation}
By computing the derivative of $g$, we get easily
\[ \max_{0\leq \eta \leq 1}g(\eta) = \left\lbrace 
\begin{array}{ll} 
\frac{2d}{2d+1} \big((2d+1)\omega_{2d} 
\|f_0\|_{L^\infty(\R^d\times\R^d)} \big)^{\frac{-1}{2d}}
&  \mbox{if } 
\omega_{2d} \|f_0\|_{L^\infty(\R^d\times\R^d)} \geq \frac{1}{2d+1},
\\ 
1-\omega_{2d} \|f_0\|_{L^\infty(\R^d\times\R^d)}
& \mbox{if } 
\omega_{2d} \|f_0\|_{L^\infty(\R^d\times\R^d)} \leq \frac{1}{2d+1}.
\end{array}\right.\]
In both cases 
$\ds \max_{0\leq \eta \leq 1}g(\eta) \geq 
\frac{2d}{2d+1} 
\big[ (\omega_{2d} (2d+1) 
\|f_0\|_{L^\infty(\R^d\times\R^d)})^{-\frac{1}{2d}} \wedge 1 \big]$.

{\it Estimate when $R>0$:}

When $R>0$ we proceed as in \eqref{R=0} in order to get for all 
$t\geq 0$ and all $\eta \in [0,1]$,
\begin{equation}
\label{R>0}
\int_{\R^d\times\R^d} \chi_R(x,v) f_t(x,v)  \ud x \ud v
\geq \eta \big(1-|K_{\eta,R} | \|f_0\|_{L^\infty(\R^d\times\R^d)}\big)
\end{equation}
where we have set
\[ K_{R,\eta} = 
\left\lbrace (x,v) \mid \chi_R(x,v) \leq \eta \right\rbrace 
\begin{array}[t]{l} \ds
=\left\lbrace (x,v) \left| |v| \leq \eta, 
|x| \leq R+(\eta^2-|v|^2)^{1/2} \right\rbrace\right.
\\ \ds
\subset B_{\R^d_x}(0,R+1)\times B_{\R^d_v}(0,\eta).
\end{array}
\]
Inserting the bound $|K_{R,\eta}|\leq \omega_d^2 (R+1)^d \eta^d$ into \eqref{R>0}, we get as above
\[\int_{\R^d\times\R^d} \chi_0(x,v) f_t(x,v)  \ud x \ud v
\geq \frac{d}{d+1} 
\big[ (\omega_{d}^2 (d+1) (R+1)^d 
\|f\|_{L^\infty(\R^d\times\R^d)})^{-\frac{1}{d}} \wedge 1 \big].\]
In both cases, we have established
\begin{equation}
\label{borninf-theta}
\int_{\R^d\times\R^d} \chi_R(x,v) f_t(x,v)  \ud x \ud v
\geq \Theta \quad \mbox{ for all } t\geq 0.
\end{equation}
Combining \eqref{Kant-R}, \eqref{test-0} and \eqref{borninf-theta},
we have proved Proposition \ref{Theta}.
\end{ProofOf}

\subsection{Lower bounds for the convergences and conclusion}

For any $Z_1 = (x_k^1,v^1_k)_{1\leq k\leq N}$ and 
$Z_2 = (x_k^2,v^2_k)_{1\leq k\leq N}$
in $(\R^d\times\R^d)^N$, we define the distance
$D^N(Z_1, Z_2) := \frac{1}{N} 
\sum_{k=1}^N (|x_k^1- x_k^2|^2 + |v_k^1- v_k^2|^2)^{1/2} \wedge 1$. 
We end by proving that the convergences established by Theorem \ref{monotonie} and Proposition \ref{kin-part} turn to be slower and slower when $N$ goes to infinity.

\begin{proposition}
\label{bornes-infs}
Assume {\bf (H1)-(H6)}, {\bf(M)} and suppose that $\Sigma$ is radially symmetric,
$f_0 \in L^\infty(\R^d\times\R^d)$ and 
$(\nabla V)^{-1}(\{0\}) \subset B(0,R)$ for some $R\geq 0$.
Taking $\Theta$ defined by Proposition \ref{Theta}, we have
\begin{enumerate}
\item Almost surely for all $t\geq 0$, 
\[ \liminf_{N\rightarrow \infty} D^N\big((p_k(t),q_k(t))_{1\leq k \leq N}, Eq^N\big) \geq \Theta. \]
\item If $\nabla^2 V \in L^\infty(\R^d)$, then taking 
$C(t)$ given by Remark \ref{C-explicite}),
for all $N,t\geq 0$ and $\epsilon \in (0,\Theta)$, we have explicitely
\begin{itemize}
\item[(i)]
$\ds \mathbb{E}\left[ D^N\big((p_k(t),q_k(t))_{1\leq k \leq N}, Eq^N\big) \right] 
\geq \Theta - \frac{C(t)}{\sqrt{N}} $
\item[(ii)]
$\ds \mathbb{P}\left[ D^N\big((p_k(t),q_k(t))_{1\leq k \leq N}, Eq^N\big)  
\leq \Theta -\epsilon \right] 
\leq \frac{1+C(t)}{\epsilon \sqrt{N}} $ 
\end{itemize}
\end{enumerate}
\end{proposition}


Before proving Proposition \ref{bornes-infs}, we now explain the division in four periods announced in the introduction. Assuming {\bf (H1)-(H6)} and 
$\lim_{|x| \rightarrow \infty} V(x) = + \infty$, we take $f_0 \in L^1(\R^d\times\R^d)$ such that
\eqref{kin}-\eqref{CI} has a unique solution (we remind the reader that a uniqueness condition is given in 
Theorem \ref{existence}-$2)$. To fix the ideas, we take
$0<\epsilon <<\Theta$ and we define the following times
\[\begin{array}{l}
T_{1,\epsilon} = \max \{t\geq 0 \mid W_1(f(t), SEq) \geq \epsilon \}
\\
T_{2,\epsilon}^N = 
\min \{t\geq 0 \mid W_1(f(t), \hat{f}^N(t)) \geq \epsilon \} 
\\
T_{3,\epsilon}^N = 
\max \{t\geq 0 \mid D^N((q_k^N(t),p_k^N(t))_{1\leq k \leq N}, Eq^N ) \geq \epsilon \} 
\end{array}
\]
$T_{1,\epsilon}$ and $T_{3,\epsilon}^N$ are well defined thanks to Theorem \ref{main} and Proposition \ref{kin-part}.
$T_{2,\epsilon}^N$ is almost surely well defined thanks to Theorem \ref{lim-chm}-$1)$.  According to Theorem \ref{lim-chm}-$1)$,
$\lim_{N \rightarrow \infty} T_{2,\epsilon}^N = +\infty$
almost surely; if $N$ is large enough then $T_{2,\epsilon}^N\geq T_{1,\epsilon}$. By construction of $\Theta$, we always have 
$T_{3,\epsilon}^N\geq T_{2,\epsilon}^N$. By definition of those times we deduce the following behaviour:
\begin{enumerate}
\item
On $[0,T_{1,\epsilon}]$ the empirical density is close to $f$, it moves from its initial value to get in an $2\epsilon$-neighbourhood of $SEq^{f_0}$, the set of spatially-stationary states dynamiquely accessible by a solution of  
\eqref{kin}-\eqref{CI}. 
\item
On $[T_{1,\epsilon},T_{2,\epsilon}^N]$, the empirical density is still close to $f$ and located in a $2\epsilon$-neighbourhood of $SEq^{f_0}$. The empirical spatial density is quasi-stationary due to Proposition \ref{semistat}-$2)$. 
\item 
On $[T_{2,\epsilon}^N,T_{3,\epsilon}^N]$, $f$ is no longer a good approximation of $\hat{f}^N$. The empirical density leaves the neighbourhood of $SEq^{f_0}$ to to get in a $\epsilon$-neighbourhood of $\mathscr J(Eq^N)$.
\item
On $[T_{3,\epsilon}^N,+\infty)$, the particles are close to $Eq^N$, the set of equilibrium of the system. Their speeds go to $0$ while their positions are at least quasi-stationary and possibly converge 
(surely under  {\bf (M)} and $x.\nabla V(x) >0$ on $\R^d\setminus \{0\}$
thanks to Theorem \ref{monotonie}).
\end{enumerate}

We end by proving Proposition \ref{bornes-infs}.

\begin{ProofOf}{Proposition \ref{bornes-infs}}

According our assumptions and thanks to \eqref{EqN-rad-M} and , 
$Eq^N \subset (B_{\R^d_q}(0,R)\times\{0_{\R^d_p}\})^N$.
Taking $\chi_R$ defined by \eqref{chi-R}, it allows us to deduce for any $N,t\geq 0$
\begin{equation}
\label{DN-chiR} 
D^N(p_k(t),q_k(t))_{1\leq k \leq N}, Eq^N) 
\begin{array}[t]{l} \ds
\geq \frac{1}{N} \sum_{k=1}^N \chi_R (p_k(t),q_k(t)) 
\\ \qquad \ds
=\int_{\R^d\times\R^d} \chi_R(x,v) \ud \hat{f}^N_t (x,v).
\end{array}
\end{equation}

{\it Proof of $(1)$}

Taking $f$ a solution of \eqref{kin}-\eqref{CI} and splitting 
$\hat{f}^N_t = f_t + (\hat{f}^N_t - f_t)$ from 
\eqref{borninf-theta} and \eqref{KR}, we get
\begin{equation}
\label{borninf1}
D^N((q_k^N,p_k^N)_{1\leq k\leq N}, Eq^N) \geq 
\Theta - W_1(\hat{f}^N(t), f(t)).
\end{equation}
We now take a subsequence $N_\ell$ such that 
\[\lim_{\ell \rightarrow \infty} 
D^{N_\ell}((q_k^{N_\ell},p_k^{N_\ell})_{1\leq k\leq N_\ell}, Eq^{N_\ell}) 
= \liminf_{N \rightarrow \infty}  
D^N(p_k(t),q_k(t))_{1\leq k \leq N}, Eq^N)\]
Almost surely from Theorem \ref{lim-chm}-$1)$, we can extract a subsequence 
(still noted $\hat{f}^{N_\ell}$) and find $f$ a solution of \eqref{kin}-\eqref{CI}  such that 
$\lim_{\ell \rightarrow \infty} W_1(\hat{f}^{N_\ell}(t), f(t) ) = 0$. coming back to \eqref{borninf1}, we have proved
\[ \liminf_{N \rightarrow \infty}  
D^N\big((p_k(t),q_k(t))_{1\leq k \leq N}, Eq^N\big) \geq \Theta. \]

{\it Proof of $(2i)$}

The proof is very similar to the previous one. For all 
$t \geq 0$, taking the expected value in \eqref{DN-chiR} we get as well
\[ \mathbb{E}[D^N((q_k^N,p_k^N)_{1\leq k\leq N}, Eq^N)]
\geq  
\int_{\R^d\times\R^d} \chi_R(x,v) \ud f^{(1,N)}_t(x,v).
\]
Splitting $f^{(1,N)}_t = f_t + (f_t - f^{(1,N)}_t)$, from \eqref{borninf-theta} 
and \eqref{KR}, we deduce
\[ \mathbb{E}[D^N((q_k^N,p_k^N)_{1\leq k\leq N}, Eq^N)] \geq 
\Theta - W_1(f^{(1,N)}_t, f(t))\]
The result follows directly from Theorem \ref{lim-chm}-$2i)$.

{\it Proof of $(2ii)$}

Applying successively \eqref{DN-chiR}, \eqref{borninf-theta} and the Markov inequality, we get 
\[  \begin{array}{l} \ds
\mathbb{P}\left[ D^N(p_k(t),q_k(t))_{1\leq k \leq N}, Eq^N)  
\leq \Theta -\epsilon \right] 
\\ \ds \qquad \qquad
\leq \mathbb{P}\left[ \frac{1}{N} \sum_{k=1}^N \chi_R (p_k(t),q_k(t)) \leq \Theta -\epsilon \right]
\\ \ds \qquad \qquad
\leq \mathbb{P}\left[ \frac{1}{N} \sum_{k=1}^N \chi_R (p_k(t),q_k(t)) 
\leq \int_{\R^d\times\R^d} \chi_R(x,v) \ud f_t(x,v) - \epsilon 
\right]
\\ \ds \qquad \qquad
\leq \mathbb{P}\left[ \left| 
\int_{\R^d\times\R^d} \chi_R(x,v) \ud \hat{f}^N_t(x,v)
-\int_{\R^d\times\R^d} \chi_R(x,v) \ud f_t(x,v)
\right| \geq \epsilon  \right]
\\ \ds \qquad \qquad
\leq \frac{1}{\epsilon}
\mathbb{E}\left[ \left| 
\int_{\R^d\times\R^d} \chi_R(x,v) \ud \hat{f}^N_t(x,v)
-\int_{\R^d\times\R^d} \chi_R(x,v) \ud f_t(x,v)
\right| \right].
\end{array} \]
And the result follows directly from Theorem \ref{lim-chm}-$2ii)$.
\end{ProofOf}

\section{Appendix}

\subsection*{A) Non trivial spatially-stationary states}

We now gives some examples of spatially-stationary states which are not stationary. In that purpose, we will use some simple property of the the linear Vlasov equation for the harmonic potential $ \mathscr V(x)=\frac{|x|^2}{2}$:
\begin{equation}
\label{harm}
\partial_t f + v.\nabla_x f = x.\nabla_v f. \\
\end{equation}
Setting
\[ A f = \Delta_v \Delta_x f - \mbox{div}_v (\nabla_x (\mbox{div}_v (\nabla_x f))), \]
our examples can easily be constructed thanks to the following observation
\begin{lemma}
Assume $d\geq 2$ and take  $f \in W^{4,1}(\R^d\times\R^d)$.
\label{patos}
\begin{enumerate}
\item[(i)] If $f$ solves \eqref{harm} for the initial data $f_0$ then 
$Af$ solves \eqref{harm} for the initial data $A f_0$ .
\item[(ii)] $f_0$ is a stationary solution of \eqref{harm} if and only if 
$Af_0$ is.
\end{enumerate}
\end{lemma}
$d\geq 2$, $\Sigma$ and $p$ being fixed, we pick $h$ and $\chi$ such that 
\begin{itemize}
\item $h \in L^1(\R^d\times\R^d)$ 
and $h$ is a stationary state of \eqref{harm},
\item $\chi \in W^{4,1}(\R^d\times\R^d)$  
and $\chi$ is not a stationary state of \eqref{harm},
\item $f_0 :=h+A\chi \geq 0$.
\end{itemize}
$\rho_h$ being the spatial density of $h$, we choose the external potential $V$ such that \eqref{VWRho} for 
$\rho =\rho_h$ can be simplified as \eqref{harm}:
\[
V(x) := \frac{|x|^2}{2} + \kappa (\Sigma \ast \rho_h )(x)
\]
By construction, $\mathcal{S}^{\rho_h} g$ solves \eqref{harm} for all $g\in \mathcal{M}_+^1(\R^d\times\R^d)$. 
On the one hand 
$\mathcal{S}^{\rho_h} f_0 
=  h  + \mathcal{S}^{\rho_h} (A \chi)$
and Lemma \ref{patos}-$ii)$ allows us to deduce that $f_0$ is not a stationary state of \eqref{harm}. On the other hand for any 
$t\in \R$, thanks to Lemma \ref{patos}-$i)$,
\[ \int_{\R^d_v} \mathcal{S}^{\rho_h}_tf_0(x,v) \ud v
= \rho_h   + \int_{\R^d_v} A \mathcal{S}^{\rho_h}_t \chi(x,v) \ud v
= \rho_h = \rho_0.\]
(the two last identity are obtained by integration by part thanks to the expression of $A$). Hence, $f_0 \in SEq$ while
$t\mapsto \mathcal{S}_t^{\rho_0} f_0$ is not stationary.

\begin{rmk}{We also point out that:}
\label{non-amortissement}

\begin{itemize}
\item $f_0$ can be chosen as smooth
as desired, even analytic or compactly supported.
\item The mapping $t\mapsto \mathcal{S}_t^{\rho_0} f_0$ is $2\pi$-periodic as all the solutions of \eqref{harm}, no damping can be expected in large time for the distribution function.
\item 
$t\mapsto \mathcal{S}_t^{\rho_0} f_0$ also solves 
\eqref{kin}-\eqref{CI} for $\Phi_0(t) = P(t) \Sigma \ast \rho_0$
and $V= \mathscr V + \kappa \Sigma \ast \rho_h$.
\end{itemize}
\end{rmk}

\begin{ProofOf}{Lemma \ref{patos}}

We consider the partial Fourier tranform on the phase space
\[\hat{f}(k,\xi) = \int_{\R^d\times\R^d} f(x,v) e^{-i(k.x+\xi.v)} 
\ud x \ud v.\]
If $f$ solves \eqref{harm}, then for $\hat{f}$ we get
\[ \partial_t \hat{f} \begin{array}[t]{l}
=- (-i) \nabla_\xi . (-i k \hat{f}) 
+ (-i) \nabla_k . (-i\xi \hat{f}) \\ \ds
= -\xi.\nabla_k \hat{f} + k \nabla_\xi \hat{f}. 
\end{array}\]
Therefore, $f$ solves \eqref{harm} if and only if $\hat{f}$ does too. 
One can check that for $a(k,\xi)=|k|^2|\xi|^2-(k.\xi)^2$, we have
$\widehat{Af}=a\hat{f}$ and that $a$ is (at least formaly) a stationary solution of \eqref{harm}. For all $f$ solution of \eqref{harm} with initial data $f_0$, since $\hat{f}$ also solves \eqref{harm}, it follows that $a\hat{f}$ solves \eqref{harm} as well for the initial data $a \hat{f_0}$. Applying the inverse Fourier transform, we deduce that $Af$ is a solution of \eqref{harm} for the initial data $Af_0$. 

By injectivity of the Fourier transform, $f_0$ is a stationary solution of 
\eqref{harm} if and only if $\hat{f_0}$ is. Since $a$ is a non zero polynomial function in $(k,\xi)$, we have from \cite{Mit}
\[ \left| \left\lbrace (k,\xi)\in \R^d\times\R^d \left| 
a(k,\xi)=0 \right.\right\rbrace \right| = 0. \]
Therefore, $\hat{f_0}$ is a stationary solution of \eqref{harm} if and only if $a \hat{f_0}$ is too. 

\end{ProofOf}
\subsection*{B) Property of $\sigma_1 \ast \sigma_1$}

According to Lemma \ref{resolution}, the parameter $\Sigma$
involved in the simplification of the Vlasov-Wave system 
\eqref{kin},\eqref{rho}-\eqref{VW} is given by 
$\Sigma = \sigma_1 \ast \sigma_1$. 
Under {\bf (A1)}, 
if one also suppose that $x.\nabla \sigma_1(x) \leq 0$ for all $x \in \R^d$, then $x.\nabla \Sigma(x) \leq 0$ as well thanks to the following lemma:

\begin{lemma}
\label{monot}
Take $\varphi_1 \in C_c^1(\R_+)$ 
and $\varphi_2 : \R_+ \rightarrow \R$ such that 
$\varphi_2 (|\cdot|) \in L^1_{loc}(\R^d)$.
We define
$u_1,u_2,u_3 : \R^d \rightarrow \R$ by
$u_1(x) = \varphi_1(|x|)$, $u_2(x) = \varphi_2(|x|)$ and
$u_3 = u_1\ast u_2$. 
If $\varphi_1$ and $\varphi_2$ are both non increasing then 
$x.\nabla u_3(x) \leq 0$ for all $x\in \R^d$.
\end{lemma}
When $\Sigma=\sigma_1\ast \sigma_1$, we also point out that 
  $\nabla^2 \Sigma(0) 
= -\frac{1}{d}\|\nabla \sigma_1 \|_{L^2(\R^d)}^2 Id$ 
as well as $\nabla \Sigma(0) =0$
under {\bf(A1)}. 
The condition
$\nabla \Sigma \neq 0$ on $B(0,r)\setminus \{0\}$ 
involved in Theorem \ref{monotonie} is always satisfied for some $r>0$.

\begin{ProofOf}{Lemma \ref{monot}}

Since $u_1$ and $u_2$ are both radially symmetric, $u_3$ is radially symmetric as well. Setting $\varphi_3(t) = u_3(te_1)$, we just have to prove that $\varphi_3' \leq 0$ on $\R_+$.
From the following expression
\[\varphi_3(t) = \int_{\R \times \R^{d-1} }
\varphi_1( [ (t-y_1)^2 + |y^\perp|^2]^{1/2} )
\varphi_2( [ y_1^2 + |y^\perp|^2]^{1/2} ) \ud y_1 \ud y^\perp,
\]
we get directly
\[ \varphi_3'(t) 
\begin{array}[t]{l} \ds
= \int_{\R \times \R^{d-1} }
\frac{y_1}{[ y_1^2 + |y^\perp|^2]^{1/2}}
\varphi_1'( [ y_1^2 + |y^\perp|^2]^{1/2} )
\varphi_2( [ (t-y_1)^2 + |y^\perp|^2]^{1/2} ) \ud y_1 
\ud y^\perp
\vspace{0.2cm}
\\ \ds
=\int_{\R^{d-1}} \int_0^\infty 
\frac{y_1}{[ y_1^2 + |y^\perp|^2]^{1/2}}
\varphi_1'( [ y_1^2 + |y^\perp|^2]^{1/2} ) 
v(t,y_1, |y^\perp|) \ud y_1 \ud y^\perp
\end{array}
\]
where we have set 
\[ v(t,y_1,r) = \varphi_2( [ (t-y_1)^2 + r^2]^{1/2} )
-\varphi_2( [ (t+y_1)^2 + r^2]^{1/2} ). \]
For all $t,y_1\geq 0$, $|t-y_1| \leq |t+y_1|$, 
since $\varphi_2$ is non increasing, it is enough to deduce that 
$v(t,y_1,r) \geq 0$ for all $t,y_1,r\geq 0$. 
We also know that $\varphi_1' \leq 0$, therefore 
$\varphi_3'(t) \leq 0$ for all $t\geq 0$.

\end{ProofOf}

\section*{Acknowledgements}
This research is supported by the Basque Government through the BERC 2018-2021 program and by Spanish Ministry of Science, innovation and Universities (Agencia Estatal de Investigación) through BCAM Severo Ochoa excellence accreditation SEV-2017-0718.

\bibliographystyle{siam} 
\footnotesize

\bibliography{LT-VW.bib}

\end{document}